\documentclass[11pt,reqno]{amsart}
 \usepackage[english]{babel}
 \usepackage[utf8]{inputenc} %This pack allows to use UTF8 codification for special characters.
\usepackage[T1]{fontenc} 
  \usepackage{textcomp}% required to get special symbols  

\usepackage[margin=3cm, marginparwidth=0.85in]{geometry}
 \usepackage{mathtools,mathrsfs}
\mathtoolsset{showonlyrefs}% use \eqref and \refeq to refer equations

\usepackage{graphicx,enumerate}
\usepackage{hyperref,color}
\usepackage[dvipsnames]{xcolor}
\hypersetup{
	colorlinks=true,
	pdfpagemode=UseNone,
    citecolor=OliveGreen,
    linkcolor=RoyalBlue,
    urlcolor=black,
	pdfstartview=FitW
}
 
\usepackage{subcaption}
\usepackage{tikz}  
\usetikzlibrary{calc, intersections, math, arrows}
\usetikzlibrary{decorations.pathreplacing,positioning}
\tikzset{
  dot/.style={
    circle, inner sep=0pt, 
    minimum size=1mm, fill=gray
  },
  solid node/.style={circle, draw, inner sep=.8, fill=black},
  red node/.style={circle, draw, inner sep=.8, fill=red},
  blue node/.style={circle, draw, inner sep=.8, fill=blue},
} 
\usepackage{comment}
\usepackage{appendix}
\usepackage{multicol} %Multiple columns
  \usepackage[lining]{libertinus}  
  \usepackage[varqu,varl]{zi4}% inconsolata for mono, not LibertineMono
 
  \usepackage[amsthm]{libertinust1math} 
\usepackage{microtype} %for fixing some overfull box  
\newcommand{\ind}[1]{\mathbf{1}_{\left\{ #1 \right\}}}
\theoremstyle{plain}
\newtheorem{theorem}{Theorem}[section]
\newtheorem{proposition}[theorem]{Proposition}
\newtheorem*{proposition*}{Proposition}
\newtheorem*{theorem*}{Theorem}
\newtheorem{lemma}[theorem]{Lemma}
\newtheorem*{lemma*}{Lemma}
\newtheorem{corollary}[theorem]{Corollary}
\theoremstyle{definition}

\newtheorem{remark}[theorem]{Remark}
\newtheorem{question}{Question}
\newcommand*{\dif}{\ensuremath{\mathop{}\!\mathrm{d}\hspace{1pt}}}

\def\Z{\mathbb{Z}}

\def\E{\mathbf{E}}
\def\P{\mathbf{P}}

\numberwithin{equation}{section}

\renewcommand{\bar}[1]{\mkern 1.5mu\overline{\mkern-1.5mu#1\mkern-1.5mu}\mkern 1.5mu}
 
\DeclareMathOperator*\diam{diam}
\DeclareMathOperator*\dist{dist}

\title{Multifractal spectrum of branching random walks on free groups}

\author{Shuwen Lai}
\address{Shuwen Lai: School of Statistics and Data Science, Nankai University, 94 Weijin Road, Nankai District, Tianjin, P.R. China} \email{1120220052(at)mail(dot)nankai(dot)edu(dot)cn}

\author{Heng Ma}
\address{Heng Ma: School of Mathematical Sciences, Peking University, Beijing, 100871, P.R. China. }\email{hengmamath(at)gmail(dot)com} \urladdr{\url{hengmamath.github.io}} 

\author{Longmin Wang}
\address{Longmin Wang : School of Statistics and Data Science \& LPMC, Nankai University, 94 Weijin Road, Nankai District, Tianjin, P.R. China } \email{wanglm(at)nankai(dot)edu(dot)cn}

\begin{document}

\begin{abstract}
Consider a symmetric branching random walk on a free group $\mathbb{F}$
in the transient regime $1<r\leq R$, where $r$ is the mean offspring
number and $R$ is the reciprocal of the spectral radius of the
underlying random walk.
 The limit set $\Lambda_r$—consisting of all ends of $\mathbb{F}$ to which the BRW’s particle trajectories converge—is a proper random subset of the boundary $\partial \mathbb{F}$.  Hueter and Lalley (2000)   determined the Hausdorff dimension of $\Lambda_{r}$, and proved that $ \dim_{\mathrm{H}} \Lambda_{r} \leq  \frac{1}{2} \dim_{\mathrm{H}} \partial \mathbb{F}$ with equality possible only when $r = R$.

We further extend this study by conducting a multifractal analysis of the limit set $\Lambda_r$.
We obtain the Hausdorff dimensions of the sub-fractals $\Lambda_{r}(\alpha) \subset \Lambda_{r}$ which consist of all ends of $\mathbb{F}$ approached by particle trajectories escaping at the rate $\alpha \in [0,1]$.  Notably, there exists a unique $\alpha(r) \in [0,1]$ such that
 \begin{equation}
  \dim_{\mathrm{H}} \Lambda_{r} = \dim_{\mathrm{H}} \Lambda_{r}( \alpha(r) ).
 \end{equation}
  Moreover, the maximizing speed exhibits a phase transition: $\alpha(r)>0$ for $1<r<R$, whereas $\alpha(R)=0$.
    \end{abstract}

% \begin{abstract}
% We consider a transient symmetric branching random walk (BRW) on a free group $\mathbb{F}$ indexed by a Galton-Watson tree $\mathcal{T}$ without leaves. 
% The limit set $\Lambda$ is defined as 
% the random subset of $\partial \mathbb{F}$ (the boundary  of $\mathbb{F}$) consisting of all ends in $\partial \mathbb{F}$ to which particle trajectories converge. 
% Hueter--Lalley [Probab. Theory Relat. Fields, 2000] determined the Hausdorff dimension of the limit set, and 
% found that $ \dim_{\mathrm{H}} \Lambda \leq  \frac{1}{2} \dim_{\mathrm{H}} \partial \mathbb{F}$. 
% In this paper, we conduct a multifractal analysis for the limit set $\Lambda$. 
% We compute almost surely and simultaneously, the Hausdorff dimensions of the sets $\Lambda(\alpha) \subset \Lambda$ consisting of all ends  in  $\partial \mathbb{F}$ to which  particle trajectories converge at rate of escape $\alpha$.
% Moreover, for isotropic BRWs, we  obtain the dimensions of the sets $\Lambda(\alpha,\beta) \subset \Lambda$ which consist of all ends  in  $\partial \mathbb{F}$  to which  particle trajectories, with the average rates of escape having limit points $[\alpha,\beta]$, converge.  
% Finally, analogous to   results of Attia--Barral (2014), we obtain the Hausdorff dimensions of the level sets $E(\alpha,\beta)$ of infinite branches in $\partial \mathcal{T}$ along which the averages of the BRW have $[\alpha,\beta]$ as the set of accumulation points.
%   \end{abstract}
 
\maketitle
 
%\tableofcontents

\section{Introduction and Main Results}
A branching random walk (BRW) on a group $\Gamma$ with a finite symmetric generating set $S$ and the identity $e$ is constructed   as follows:
\begin{enumerate}[(i)]
  \item  Sample a Galton-Watson tree $\mathcal{T}$ rooted at $\emptyset$ with   offspring distribution $p=(p_k)_{k\geq 0}$ and mean offspring $r:=\sum_{k \geq 1} k p_{k}$.
  \item Independently assign each  non-root vertex $u \in \mathcal{T}$ a random element $Y_{u} \in S$ according to a symmetric\footnote{Symmetric here means that $\mu(a)=\mu(a^{-1})$ for all $a \in S$.} probability measure $\mu$ on $S$. For every non-root vertex $u \in \mathcal{T}$, define $ V(u)=Y_{u_{1}} Y_{u_{2}} \cdots Y_{u_{n}}$  where $(\emptyset,u_{1},\cdots,u_{n}=u)$ is the geodesic in $\mathcal{T}$ from $\emptyset$ to $u$. Moreover, set   $V(\emptyset)=e$.
\end{enumerate}
  We refer to $(V(u): u \in \mathcal{T})$ as a $\mu$-branching random walk on $\Gamma$.
  Let $G$ be the Cayley graph of $(\Gamma,S)$. Then, the process $(V(u) : u \in \mathcal{T})$ can also be regarded as a random walk on $G$ indexed by the Galton–Watson tree $\mathcal{T}$ in the sense of Benjamini and Peres \cite{BP94}.

Random walks (RWs) and Branching random walks (BRWs) have been extensively studied in Euclidean spaces. However, it has been discovered that BRWs on hyperbolic spaces (and more generally, nonamenable graphs) exhibit a phase that is not present in the corresponding processes in Euclidean spaces.
For instance,  let $(\Gamma,S)$ denote a nonelementary hyperbolic group. The random walk on $\Gamma$ with a symmetric step distribution $\mu$ supported on $S$ has spectral radius $R^{-1}$ strictly less than $1$, i.e.,
\begin{equation}
 R^{-1} := \limsup_{n \to \infty} p_{n}(x,  y) ^{1/n} < 1,
\end{equation} where  $p_{n}(x,y)$ is the $n$-step transition probability  for the $\mu$-random walk on $\Gamma$. Due to this property, BRWs on $\Gamma$ are particularly interesting because of the following \textit{double phase transition} (see Benjamini and Peres \cite{BP94} for the case $r \neq R$ and Gantert and Müller \cite{GM06} for the critical case $r =R$):
\begin{itemize}
  \item If the mean offspring $r \leq 1$ and   $p_1<1$, the process becomes extinct almost surely. 
   \item If $1<r \leq R$, with positive probability the  process survives forever, but eventually vacates every compact subset of the state space with probability one. This is known as the \emph{transient} phase or \emph{weak survival} phase.
  \item If $r > R$, the process survives forever with positive probability  and when it survives, every vertex of $G$ is visited by infinitely many particles of the BRW. This is known as the \emph{recurrent} phase or \emph{strong survival} phase.
\end{itemize}
In contrast, for a nearest-neighbor symmetric random walk on $\mathbb{Z}^{d}$, the spectral radius is exactly $1$, so the transient phase cannot be observed in BRWs on $\Z^{d}$. See \cite{Lal06} for additional processes that exhibit this weak/strong survival transition.

Interesting questions arise in the transient regime $r \in (1, R]$, where the \textit{limit set} $\Lambda_{r}$, defined as the random subset of the Gromov boundary
$\partial \Gamma$ (endowed with the visual metric) consisting of those points to which BRW's particle trajectories converge, is a proper random subset  of $\partial \Gamma$.
Perhaps the most basic and important characteristic of a random fractal is its Hausdorff dimension. Under the previous setting,  the Hausdorff dimension of $\Lambda_{r}$ is related to the growth rate of the trace of the BRW defined by
\begin{equation}
  \mathsf{Gr}(r) := \limsup_{n \to \infty}   \# \{ x \in \Gamma: |x|=n, x \text{ is visited by the BRW}  \} ^{1/n}  .
\end{equation}
where $|\cdot|$ denotes the word length of $x$ in group $ (\Gamma, S)$. Precisely,  it has been shown that
for any $r \in (1,R]$,
\begin{equation}\label{eq-hdim-lim-set-hyper}
  \dim_{\mathrm{H}} (\Lambda_{r}) \propto \ln \mathsf{Gr}(r) \,\text{ and } \, \dim_{\mathrm{H}} (\Lambda_{r}) \leq \frac{1}{2}  \dim_{\mathrm{H}} ( \partial \Gamma) \quad \text{ a.s. }
\end{equation}
In particular, the dimension of $\Lambda_{r}$ has a phase transition at the critical value $r=R$.
We refer to Sidoravicius, Wang, and Xiang \cite{SWX23} for the case  $r \in (1,R)$ and  Dussaule, Wang, and Yang
\cite{DWY22} for the critical case  $r = R$. The latter paper
\cite{DWY22} extends the corresponding results  to BRWs on relative hyperbolic groups.
Notably,   problems of  type  of  \eqref{eq-hdim-lim-set-hyper}
were initially investigated by  Liggett \cite{Liggett96} for isotropic BRWs on $d$-regular trees (Cayley graph of the free product $(\Z_{2})^{*d}$); by Lalley and Sellke \cite{LS97} for  branching Brownian motion on hyperbolic plane $\mathbb{H}^2$;  by Hueter and Lalley  \cite{HL00} for  anisotropic BRWs on $d$-regular trees;  and by Candellero, Gilch and M\"{u}ller \cite{CGM12}  for BRWs on free products of groups.  Here we say the BRW is isotropic if random walk is simple possibly with laziness, i.e., $\mu(x)=\mu(y)$ for all $x,y$ in $S$; otherwise it is anisotropic. Remarkably, Hueter and Lalley \cite{HL00} showed that for BRWs on $d$-regular trees,
\begin{quotation}
  \normalsize\emph{$\dim_{\mathrm{H}} (\Lambda_{r}) = \frac{1}{2}  \dim_{\mathrm{H}} ( \partial \Gamma)$  in \eqref{eq-hdim-lim-set-hyper}    holds if and only if  $r= R$ and the underlying random walk is isotropic.}
\end{quotation}

 From a fractal geometry standpoint, a single exponent---the fractal dimension given in \eqref{eq-hdim-lim-set-hyper}---is insufficient to fully capture the characteristics of a random fractal like $\Lambda_r$. Instead, a continuous spectrum of exponents, known as the multifractal spectrum, is required. 
Roughly speaking, the multifractal spectrum  reflects the spatial heterogeneity of fractal patterns.
(We refer  \cite{Mande99} and  \cite[Chapter 17]{Falconer03} for an introduction to multifractal analysis.)
A natural and intriguing question is to establish the multifractal analysis of the limit set $\Lambda_{r}$ in the transient regime $r \in (1,R]$. In this case, the branching random walk can be extended to a continuous mapping from $\partial \mathcal{T}$ to $\partial \Gamma$ (see \cite[Section 5]{SWX23}):
\begin{equation}\label{eq-BRW-mapping}
  V : \quad \partial \mathcal{T}  \ni t=(t_{n})_{n \geq 0} \longmapsto   \lim_{n\to\infty}V(t_n)\in\partial\Gamma.
\end{equation}
We remark that in a forthcoming paper by the first and third authors, it will be shown that the map $V$ is injective on $\partial \mathcal{T}$. In the setting current paper where $\Gamma$ is a free group, this  injectivity follows directly from   Hutchcroft’s result \cite{Hut20}, see \S\ref{sec-BRWonF}. Consequently, for each $\omega \in \Lambda_{r}$, there is a  unique ray $t \in \partial\mathcal{T}$ such that $V(t)=\omega$.  Inspired by the work of Attia and Barral \cite{AB14}, we propose using  the  escape rate  of the walk $(V(t_n))_{n \geq 1}$ to  describe the  degree of singularity  around the point $\omega=V(t)$ in the fractal $ \Lambda_{r}$.   For each $\alpha \in [0,1]$, define
\begin{equation}
  \Lambda_r(\alpha) := \Bigl\{    \omega \in \partial \Gamma :\:  \exists \, t \in \partial \mathcal{T}, V(t)= \omega \text{ s.t. } \lim_{n \to \infty} \frac{| V(t_n) |}{n}= \alpha   \Bigr\}.
\end{equation}
Let $\Lambda_{r}^{\mathrm{nl}} := \{   \omega \in \partial \Gamma :   t \in \partial \mathcal{T}, V(t)= \omega \text{ and } \lim\limits_{n \to \infty}  { | V(t_{n})  |}/{n} \text{ does not exist}   \}$. Since $V|_{\partial \mathcal{T}}$ is injective,  with probability one, $\Lambda_{r}$ can be decomposed into the following disjoint subsets:
\begin{equation}
  \Lambda_{r} =   \Lambda_{r}^{\mathrm{nl}}  \cup \bigcup_{\alpha \in [0,1]} \Lambda_{r}(\alpha) .
\end{equation}
The problem of multifractal analysis can then be formalized as follows:

\begin{question}\label{ques1}For a BRW  on a nonelementary hyperbolic group, find a nonnegative function $f_{r}$ with domain $J_{r} \subset [0,1]$ such that almost surely  for every $\alpha \in [0,1]$, $\Lambda_{r}(\alpha)$ is non-empty if and only if $\alpha \in J_{r}$ and in this case
\begin{equation}
 \dim_{\mathrm{H}}  \Lambda_r(\alpha)=   f_{r}(\alpha) .
\end{equation}
\end{question}

Our primary goal is to determine the multifractal spectrum of the limit
set $\Lambda_r$ in the transient regime $r\in(1,R]$. We address this
problem for symmetric BRWs on \textbf{free groups}.
Free groups provide a natural model setting: their Cayley graphs are
trees, making the boundary geometry explicit, while BRWs on them
already exhibit the competition between population growth and the
rarity of atypical escape speeds that determines the multifractal
spectrum.

Specifically, let $\mathbb{F}=\mathbb{F}^{d}$ be the free group of rank
$d\geq2$ with symmetric generating set
$\mathcal{A}=\{a_i,a_i^{-1}:1\leq i\leq d\}$,
and equip its boundary $\partial\mathbb{F}$ with the standard
ultrametric(see \S\ref{sec-free-group}). Let $\mu$ be a symmetric,
possibly lazy,   probability measure on
$\mathcal{A}\cup\{e\}$ such that
\[
\mu(a)=\mu(a^{-1})>0
\quad\text{for every }\ a\in\mathcal{A}.
\]
We assume that the offspring distribution
$p=(p_k)_{k\geq0}$ satisfies $p_0=0$ and has a finite exponential
moment, namely,
$
\sum_{k\geq0}e^{sk}p_k<\infty$ for some $s>0$.
Under these assumptions, we solve Question~\ref{ques1} in the
free-group setting.

To state our first main result, let $(Z_n)_{n\geq0}$ be the
$\mu$-random walk on $\mathbb{F}$ started at $e$, and let $L^*$ denote
the large-deviation rate function for its normalized word length
$|Z_n|/n$, established in  \cite[Theorem 7.2]{L93}.
The function $L^*$ is convex, is finite precisely on
$[0,1]$, and satisfies $L^*(0)=\ln R$ (see
Propositions~\ref{prop-word-length-ldp} and~\ref{prop-LDP-|Zn|}).
For each $r>1$, define
\begin{equation}\label{def-I-r}
 I(r):=\bigl\{q\in[0,1]:L^*(q)\leq\ln r\bigr\}
      =[I_-(r),I_+(r)].
\end{equation}
 Since $L^*(0)=\ln R$, we have
$I_-(r)>0$ for $r<R$ and $I_-(r)=0$ for $r\geq R$.  Therefore, for $\alpha \in I(r)$, $\alpha = 0$ is permissible only when $r \geq R$, and in this case we interpret $\frac{\ln R - L^{*}(\alpha)}{\alpha}$ as $\lim_{q \downarrow 0} \frac{L^{*}(0) - L^{*}(q)}{q} =-(L^{*})'(0)$ throughout this paper.

\begin{theorem}\label{thm-Hdim-Lambda-alpha}
Let $r \in (1,R]$. Almost surely  simultaneously for every $\alpha \in [0,1]$,  $\Lambda_r(\alpha)$ is  nonempty   if and only if $\alpha \in I(r)$. In this case, the Hausdorff dimension of $\Lambda_r(\alpha)$ is given by
\begin{equation}
  \dim_{\mathrm{H}} \Lambda_r(\alpha) =  \frac{\ln r - L^{*}(\alpha)}{\alpha} \ .
\end{equation} 
Furthermore, we have
\begin{equation}
  \dim_{\mathrm{H}} \Lambda_{r} = \max_{ \alpha \in I(r) } \frac{\ln r - L^{*}(\alpha)}{\alpha}   =\max_{ \alpha \in I(r) } \, \dim_{\mathrm{H}} \Lambda_{r}(\alpha) .
\end{equation}
\end{theorem}

\begin{figure}[htbp]
  \begin{subfigure}[t]{0.45\textwidth}
  \begin{tikzpicture}[decoration={brace,amplitude=5}, declare function={ phi(\a, \x) = \a*(\x - 2)^2; },]
    \pgfmathsetmacro{\a}{0.6}
    \pgfmathsetmacro{\r}{1.4}
    \pgfmathsetmacro{\b}{sqrt(\r/\a)}
    \pgfmathsetmacro{\q}{sqrt(4-\r/\a)}
    % \draw[style=help lines,step=1cm] (0,0) grid (1,1);
    \draw[->,> = latex', thin] (0,0) -- (4.5,0) node[right] {$\scriptstyle \alpha$};
    \draw[->,> = latex', thin] (0,0) -- (0,3);
    \draw[domain=0:3.7,smooth,variable=\x] plot ({\x},{phi(\a, \x)}) node[right] {$\scriptstyle L^{*}$};
    \draw[densely dashed, very thin] (0, {\r}) -- (3.8, {\r});
    \draw[densely dashed, very thin] ({2-\b}, {\r}) -- ({2-\b}, 0) node[below]{$\scriptstyle I_-(r)$};
    \draw[densely dashed, very thin] ({2+\b}, {\r}) -- ({2+\b}, 0) node[below] {$\scriptstyle I_+(r)$};
    \draw[thin,red] (0, {\r}) -- ({\q}, {\a*(\q - 2)^2}) node[circle,fill,inner sep=0.8pt]{};
    \draw[densely dashed, very thin] (0, {\r}) -- (3, {\a}) node[circle,fill,inner sep=0.8pt]{};
    \draw[densely dashed, very thin, densely dashed] (3, {\r}) -- (3, {\a}); % node[below] {$\scriptstyle \alpha$};
    \draw[densely dashed, very thin, densely dashed] ({\q}, {\a*(\q - 2)^2}) -- (\q, 0) node[circle,fill,inner sep=0.8pt]{};
    \node[below] at ({\q}, 0) {$\textcolor{Red}{\scriptstyle \alpha(r)}$};
    \node[left] at (0, 0) {$\scriptstyle 0$};
    \node[left] at (0, {\r}) {$\scriptstyle \ln r$};
    \node[left] at (0, {4*\a}) {$\scriptstyle \ln R$};
    \node[below] at (2.1, 0) {$\scriptstyle \mathtt{C}_{\mathrm{RW}}$};
    \node[circle, fill, inner sep=0.8pt] at (2, 0) {};
    \node[circle, fill, inner sep=0.8pt] at (0, {\r}) {};
    \node[circle, fill, inner sep=0.8pt] at (0, {\a*4}) {};
  \draw[decorate, very thin] (3, {\r}) -- (3, {\a});
  \draw[->, very thin] (3, {0.5*(\r + \a)}) ++(0.18cm, 0) to[out=0,in=180] (4, {0.45*(\r + \a)}) node[right] {$\scriptstyle\ln r - L^{*}(\alpha)$};
  \end{tikzpicture}
  % \caption{Subcritical case $1 < r < R$.}
  \end{subfigure}
  ~
  \begin{subfigure}[t]{0.45\textwidth}
  \begin{tikzpicture}[decoration={brace,amplitude=5}, declare function={ phi(\a, \x) = \a*(\x - 2)^2; },]
    \pgfmathsetmacro{\a}{0.6}
    \pgfmathsetmacro{\r}{1.4}
    \pgfmathsetmacro{\b}{sqrt(\r/\a)}
    \pgfmathsetmacro{\q}{sqrt(4-\r/\a)}
    % \draw[style=help lines,step=1cm] (0,0) grid (1,1);
    \draw[->,> = latex', thin] (0,0) -- (4.5,0) node[right] {$\scriptstyle \alpha$};
    \draw[->,> = latex', thin] (0,0) -- (0,3);
    \draw[domain=0:3.7,smooth,variable=\x] plot ({\x},{phi(\a, \x)}) node[right] {$\scriptstyle L^{*}$};
    \draw[densely dashed, thin] (0, {4*\a}) -- (3, {\a});
    \draw[red] (-0.2, {4.8*\a}) -- (0.7, {1.2*\a});
    \node[circle, fill, inner sep=0.8pt] at (3, {\a}) {};
    \node[below] at (2.1, 0) {$\scriptstyle \mathtt{C}_{\mathrm{RW}}$};
    % \node[left] at (0, 0) {$\scriptstyle 0$};
    \node[left] at (0, {4*\a}) {$\scriptstyle \ln R$};
    \node[below] at (0, 0) {$\textcolor{Red}{\scriptstyle \alpha(R) = 0}$};
    \node[circle, fill, inner sep=0.8pt] at (0, 0) {};
    \node[circle, fill, inner sep=0.8pt] at (2, 0) {};
    \node[circle, fill, inner sep=0.8pt] at (0, {\a*4}) {};
    \draw[densely dashed, very thin] (0, {4*\a}) -- (4, {4*\a});
    \draw[densely dashed, very thin] (3, {4*\a}) -- (3, {\a});
    \draw[decorate, very thin] (3, {4*\a}) -- (3, {\a});
    \draw[->, very thin] (3, {0.5*(4*\a + \a)}) ++(0.18cm, 0) to[out=0,in=180] (4, {0.45*(4*\a + \a)}) node[right] {$\scriptstyle\ln R - L^{*}(\alpha)$};
  \end{tikzpicture}
  \end{subfigure}
  \caption{Illustration for $\alpha(r)$ in subcritical case $1 < r < R$ (left) and critical case $r = R$ (right). }
  \label{fig:alphar}
  \end{figure}
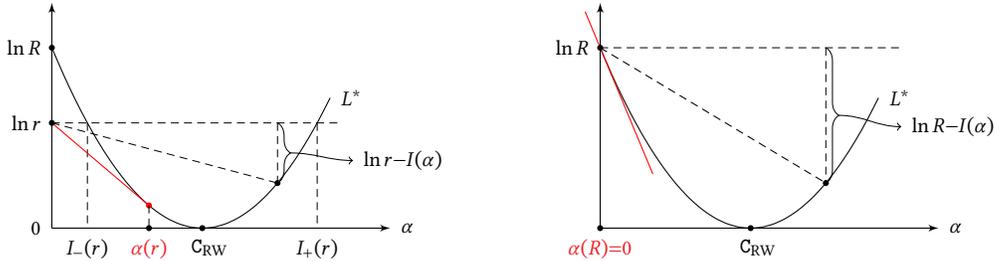

Since the rate function $L^{*}$ is convex, vanishes at  $\mathtt{C}_{\mathrm{RW}}$ (the escape rate of the RW) and is strictly decreasing on $[0, \mathtt{C}_{\mathrm{RW}} ]$ (see \S\ref{sec-RW}),
for any $r \in (1,R]$, there exists a unique $\alpha(r) \in [0, \mathtt{C}_{\mathrm{RW}})\cap I(r)$ such that
 $\max_{ \alpha \in I(r) } [\ln r - L^{*}(\alpha)]/\alpha $ is attained at $\alpha(r)$. See Figure \ref{fig:alphar} for an illustration of the location of $\alpha(r)$.
As a consequence of Theorem  \ref{thm-Hdim-Lambda-alpha}, there holds
  \begin{equation}
   \dim_{\mathrm{H}} \Lambda_{r}(\alpha(r)) =\dim_{\mathrm{H}} \Lambda_{r} \ \text{ and } \ \dim_{\mathrm{H}} \Lambda_{r}(\alpha) < \dim_{\mathrm{H}} \Lambda_{r}, \forall \, \alpha \neq \alpha(r).
  \end{equation}
That is, the subfractal $  \Lambda_{r}(\alpha(r))$, consisting of all points in $\partial \mathbb{F}$ to which particle trajectories with an escape rate $\alpha(r)$ converge, contributes the  dimension of the limit set $\Lambda_{r}$.
All other subfractals $  \Lambda_{r}(\alpha)$ are some lower dimensional structures filled in the seams of $  \Lambda_{r}(\alpha(r))$.

It may be unexpected that $\alpha(r) \neq \mathtt{C}_{\mathrm{RW}}$ which is the rate of escape of the $\mu$-random walk on $\mathbb{F}$. More surprisingly, when $r= R$, it holds $\alpha(R) = 0$, while  Theorem \ref{thm-dim-E-alpha-beta} yields that the set of genealogical rays in $\partial \mathcal{T}$ along which particle trajectories have escape rate zero, has Hausdorff dimension zero!
We interpret this as an \textit{energy-entropy competition}: although particle trajectories with lower escape rates are significantly fewer in number, the corresponding shadows in
$\partial\mathbb{F}$ shrink at a slower exponential rate; this slower
spatial contraction can compensate for their lower abundance and yield
a larger Hausdorff dimension. In the critical case $r=R$,  the energy dominates; whereas in the subcritical case,   $I_{-}(r)<\alpha(r)<\mathtt{C}_{\mathrm{RW}}$   represents a compromise between competing energy and entropy.

 Although  the multifractal analysis of the limit set $\Lambda_{r}$, in the specific context where $\Gamma= \mathbb{F}$, is established,  many intriguing questions beyond multifractal analysis remain to be explored. Inspired by Attia and Barral \cite{AB14}, we are particularly interested  in further understanding  $\Lambda^{\mathrm{nl}}_{r}$.
  For $0 \leq \alpha \leq \beta \leq 1$,  define
\begin{equation}
  \Lambda_r(\alpha, \beta) :=  \Bigl\{    \omega \in \partial \Gamma :  \exists t \in \partial \mathcal{T}, V(t)= \omega \text{ s.t. }  \varliminf_{n \to \infty} \frac{\left| V(t_n) \right|}{n} = \alpha, \, \varlimsup_{n \to \infty} \frac{\left| V(t_n) \right|}{n} = \beta \Bigr\} .
\end{equation}
In particular, $\Lambda_{r}(\alpha)= \Lambda_{r}(\alpha, \alpha)$ and $  \Lambda_{r}^{\mathrm{nl}} = \cup_{0 \leq \alpha < \beta \leq 1}  \Lambda_r(\alpha, \beta)$. Since  $V|_{\partial\mathcal{T}}  $ is injective,  the limit set $\Lambda_{r}$ can be written as a disjoint union:
 \begin{equation}
  \Lambda_{r} = \bigcup_{0 \leq \alpha \leq \beta \leq 1}  \Lambda_r(\alpha, \beta)  .
\end{equation}
Our second result concerns the Hausdorff dimension of each sub-fractal $\Lambda_{r}(\alpha,\beta)$.

\begin{theorem}\label{thm-Hdim-Lambda-alpha-beta}
Let $ r \in (1, R]$.
Almost surely, simultaneously for every  $[\alpha,\beta] \subset [0,1]$, $ \Lambda_r(\alpha, \beta) $ is nonempty if and only if $[\alpha,\beta] \subset I(r)$, and in this case,
\begin{equation}\label{eq-hdim-Lambda-alpha-beta}
    \min_{q \in \{\alpha,\beta\}} \frac{\ln r - L^{*}(q)}{q} \leq   \dim_{\mathrm{H}} \Lambda_r(\alpha, \beta) \leq   \frac{\ln r -L^{*}(\alpha)}{\alpha}. 
\end{equation}
 Additionally, if    $\mu$ is  isotropic (i.e., $\mu(a_{i})=\mu(a_{i}^{-1})= \frac{1- \mu(e)}{2d}$ for $i \leq d$), then almost surely,
\begin{equation}\label{eq-hdim-Lambda-alph}
    \dim_{\mathrm{H}} \Lambda_r(\alpha, \beta) = \min_{q \in \{\alpha,\beta\}} \frac{\ln r - L^{*}(q)}{q} \quad \text{   for every } \ [\alpha,\beta] \subset I(r).
\end{equation}
  \end{theorem}

    Taking $\beta=\alpha$ in
Theorem~\ref{thm-Hdim-Lambda-alpha-beta} recovers
Theorem~\ref{thm-Hdim-Lambda-alpha}. For anisotropic BRWs, however, the
bounds in \eqref{eq-hdim-Lambda-alpha-beta} need not coincide. In fact, we prove in 
Lemma~\ref{lem-assume-Hypothesis}  that, for a symmetric step
distribution $\mu$, the identity \eqref{eq-hdim-Lambda-alph} holds
whenever the pressure comparison condition
\eqref{eq-pressure-comparison} is satisfied. This condition is
automatic in the isotropic case and may also hold for some anisotropic
step distributions. Preliminary numerical experiments, however, do not
support its validity for arbitrary symmetric step distributions.
Accordingly, we do not expect \eqref{eq-hdim-Lambda-alph} to hold in
full generality, although we currently have no counterexample.
We suspect, although this is far from clear, that the pressure
comparison condition is not only sufficient but also necessary for \eqref{eq-hdim-Lambda-alph} to hold throughout the transient regime. 
In other words, condition \eqref{eq-pressure-comparison} may characterize
the symmetric step distributions for which the dimension formula
\eqref{eq-hdim-Lambda-alph} remains valid. Determining whether such a
characterization holds, and identifying the correct dimension formula
when the condition fails, remain open problems.

More broadly, it is natural to ask to what extent the dimensions of the
irregular sets $\Lambda_r(\alpha,\beta)$ are determined by the
multifractal spectrum in more general settings. This motivates the
following question. 

\begin{question}\label{ques2}
Suppose that the admissible set $J_r$ and the multifractal spectrum
$f_r$ in Question~\ref{ques1} have been determined. Characterize the
intervals $[\alpha,\beta]\subset[0,1]$ for which
$\Lambda_r(\alpha,\beta)$ is nonempty, and determine its Hausdorff
dimension. In particular, under what conditions  do we have,
\begin{equation}
  \dim_{\mathrm H}\Lambda_r(\alpha,\beta)
  =\min_{q\in[\alpha,\beta]}f_r(q)?
\end{equation}
Can this formula fail in general? If so, is
$\dim_{\mathrm H}\Lambda_r(\alpha,\beta)$ nevertheless determined
solely by the multifractal spectrum $f_r$? 
\end{question}

In Theorem \ref{thm-Hdim-Lambda-alpha} and \ref{thm-Hdim-Lambda-alpha-beta}, we study a family of fractals living on  $  \partial \mathbb{F}$. It is  worth noting that the BRW also induces a family of fractals living on $\partial \mathcal{T}$, the boundary of the underlying Galton-Watson tree endowed with the standard ultrametric distance.
Specifically, for each $0 \leq \alpha \leq \beta \leq 1$, define
\begin{equation}\label{def-E-alpha-beta}
E_r(\alpha,   \beta) =  \Bigl\{   t \in \partial \mathcal{T} \colon \liminf_{n \to \infty} \frac{| V(t_n)|}{n} = \alpha, \, \limsup_{n \to \infty} \frac{| V(t_n)|}{n} = \beta  \Bigr\}
\end{equation}
and write $E_r(\alpha) := E_r(\alpha, \alpha)$.  In the transient regime $r \in (1, R]$, the set $E_r(\alpha, \beta)$ coincides with the preimage of $\Lambda_r(\alpha, \beta)$ under the mapping $V|_{\partial \mathcal{T}}  $ defined in~\eqref{eq-BRW-mapping}.
Note that $E_r(\alpha,\beta)$ is well-defined not only in the transient
regime but also in the recurrent regime and, more generally, whenever
the state space of the BRW is equipped with a metric and a distinguished
origin. By contrast, the definition of $\Lambda_r(\alpha,\beta)$
requires a natural boundary of the state space. This broader
applicability makes $E_r(\alpha,\beta)$ the natural object for extending
the preceding dimension questions.
Motivated by   the work of Attia and
Barral~\cite{AB14}, we ask:

\begin{question}\label{ques3}
Consider a BRW on a rooted graph $(G,o)$, for instance in the setting
of~\cite{BP94}. Suppose that $d_G(o,Z_n)/n$ satisfies a
large deviation principle with rate function $L^*$. Does the following
statement hold?
Almost surely, simultaneously for every
$0\leq\alpha\leq\beta<\infty$,  $E_r(\alpha,\beta)$ is nonempty
if and only if
$\max_{\alpha \le q \le \beta } L^{*}(q) \le \ln r$
and, whenever it is nonempty,
\[
  \dim_{\mathrm{H}}  E_r(\alpha,\beta) =   \ln r-\max_{q \in [\alpha, \beta]}  L^{*}(q) .
\]
\end{question}

For BRWs in Euclidean space $\mathbb{R}^d$, one can refine the
decomposition in \eqref{def-E-alpha-beta} by recording the full set of
accumulation points of
$(V(t_n)/n)_{n\geq1}$. Related multifractal spectra in dimension one
have been extensively studied; see
\cite{Barral00,BHJ11,Fal94,HW92,Mol96}. The vector-valued case was
treated by Attia~\cite{Attia14} and Attia and Barral~\cite{AB14}. In
particular, the latter determined the dimensions of level sets with a
prescribed closed connected set of accumulation points, thereby going
beyond classical multifractal analysis, which corresponds to singleton
sets. Theorem~\ref{thm-dim-E-alpha-beta} provides a free-group analogue
for radial escape speeds and answers Question~\ref{ques3} in the
present setting.

\begin{theorem}\label{thm-dim-E-alpha-beta}
 Let $r \in (1,\infty)$. Almost surely, simultaneously for every  $[\alpha,\beta] \subset [0,1]$,  $E_r(\alpha,  \beta) $ is non-empty if and only if $[\alpha,\beta] \subset I(r)$ and  in this case,  $E_r(\alpha,  \beta)$ is a dense subset of  $\partial \mathcal{T}$ with Hausdorff dimension
\begin{equation}
  \dim_{\mathrm{H}}   E_r(\alpha,  \beta)   = \ln r - \max\{   L^{*}(\alpha)  , L^{*}(\beta)\} .
\end{equation}
 In particular, for any $\alpha \in I(r)$,  $ \dim_{\mathrm{H}} E_r(\alpha) =  \ln r -   L^{*}(\alpha) $; and $\dim_{\mathrm{H}} \partial \mathcal{T} = \dim_{\mathrm{H}} E_r(\mathtt{C}_{\mathrm{RW}})$.
\end{theorem}

\subsection{Related work}
The trace $\mathsf{Tr}$ of a transient simple BRW on a Cayley graph $G$ is a proper  (random) subgraph induced by the vertices that are ever visited. (The limit set is actually the boundary of this trace.)
 Benjamini and M\"{u}ller \cite{BM12} showed that the $\mathsf{Tr}$ is a.s. transient for the simple RW but is a.s. strongly recurrent for any non-trivial simple BRW.
 A large number of questions concerning the geometric properties of the trace were posed in  \cite{BM12}, including  Question 4.1  which  asked whether   $\mathsf{Tr}$  is infinitely-ended with no isolated ends.
This question was partially answered in \cite{CR15, GM17} and resolved  in \cite{Hut20} with full generality.
Recent  work by \cite{Hut22} explored
how the recurrence and transience of space-time sets for a BRW on a graph depends on the offspring distribution. For a different perspective, \cite{CH23,KW23,Woess24}  investigates  the random limit boundary measures arising from the empirical distributions of sample populations. Finally, we refer to \cite{DKLT22} for a study of scaling limits of BRWs on trees in a completely different context.

\subsection{Overview of the proofs}
Although we formulate our results for free groups, many of our arguments
rely primarily on the tree structure of the Cayley graph and should extend,
with suitable modifications, to symmetric BRWs on
$(\mathbb{Z}_2)^{*d}$ for $d\geq3$  (whose Cayley graph is the
$d$-regular tree.)
 In what follows, for $x\in\mathbb{F}$, let $\mathbb{F}(x)$ denote the
forward cone consisting of all words having $x$ as a prefix, and let
$[x]_{\mathbb{F}}$ denote the shadow consisting of all ends in
$\partial\mathbb{F}$ having $x$ as a prefix.  Similarly, for
$u\in\mathcal{T}$, let $\mathcal{T}(u)$ denote the descendant subtree
rooted at $u$, and let $[u]_{\mathcal{T}}$ denote the set of
genealogical rays passing through $u$.

\medskip
\noindent
\underline{Input 1: Uniform random-walk estimates.}
Our first input is the uniform local large-deviation estimate
\begin{equation}
  \label{eq:intro-RW}
 \P\bigl(|Z_n|=m\bigr)
 = \exp\bigl\{-nL^*(m/n) + o (n)\bigr\} \
\text{ for every  admissible }\  m
\end{equation}
see Proposition~\ref{prop-LDP-|Zn|}.
Moreover, our method yields a new representation of the rate function $L^{*}$, which differs from that in \cite{L93}; see \eqref{eq-forluma-Lstar}.
The uniformity in $m$ is essential
at the critical point $r=R$, where speeds close to zero must be retained.
Indeed, the many-to-one formula shows that the expected number of
generation-$n$ particles $u$ with  distance exactly $m$ from $e \in \mathbb{F}$ is
$r^n\P(|Z_n|=m)$.  Thus its exponential rate is
$\ln r-L^*(m/n)$, and the condition $L^*(q)\leq\ln r$ defining $I(r)$
is exactly the condition that this rate be nonnegative.

\medskip
\noindent
\underline{Upper Bounds.} To provide an upper bound for the Hausdorff dimension of a fractal, a common approach is to construct an appropriate covering of it. Since each
ray in $E_r(\alpha,\beta)$ has subsequences
with   speed arbitrarily close to every $q\in[\alpha,\beta]$, we have
\begin{equation}\label{eq-overview-proof-1}
  E_r(\alpha,\beta)\subset
  \bigcap_{\epsilon>0}
  \bigcup\left\{[u]_{\mathcal{T}}:u\in\mathcal{T},
  \ |V(u)|/|u|\in[q-\epsilon,q+\epsilon]\right\}.
\end{equation}
For  generation $n$ particle $u$, the diameter of $[u]_{\mathcal{T}}$  is $e^{-n}$. While by  \eqref{eq:intro-RW},
the
number of generation $n$ particles  with average speed  $  q$ has exponential order
$\exp\{n[\ln r-L^*(q)]\}$. Thus the $s$-dim Hausdorff out-measure   converges for all $s> \ln r-L^*(q)$.  Optimizing over $q$ and the convexity of $L^*$ give
\begin{equation}
 \dim_{\mathrm{H}}E_r(\alpha,\beta)
 \leq \min_{q\in[\alpha,\beta]}\bigl\{\ln r-L^*(q)\bigr\}
 =\ln r-\max\bigl\{L^*(\alpha),L^*(\beta)\bigr\}.
\end{equation}
If $q\notin I(r)$, the same first-moment argument   and the
Borel--Cantelli lemma then rules out any ray having $q$ as a subsequential
escape rate.

 The difficulty is to construct a valid cover in the spatial
boundary. A shadow $[x]_{\mathbb{F}}$   has diameter
$e^{-m}$ if $|x|=m$. For $ m\approx qn$, the expected number of words of length
$m$ occupied by generation-$n$ particles has exponential order at most
$\exp\{n[\ln r-L^*(q)]\}$. Comparing this growth rate with the shadow
diameter $e^{-m}$ suggests the exponent
\begin{equation}
 h_r(q)=\frac{\ln r-L^*(q)}{q}.
\end{equation}
However, particle trajectories descending from $u$ may backtrack past
$V(u)$ before converging to the boundary. Hence,  we have
$V([u]_{\mathcal{T}})\nsubseteq[V(u)]_{\mathbb{F}}$ in general. This is
the obstruction to transferring the covering argument in
\eqref{eq-overview-proof-1} directly to
$\Lambda_r(\alpha,\beta)$.
Lemma~\ref{lem-dim-hat-L-1} provides a coarser valid cover. It shows that
every $\omega\in\Lambda_r(\alpha,\beta)$ has infinitely many prefixes
$\omega_m$ visited by its associated particle trajectory at times $n$
for which $m/n$ comes arbitrarily close to $\alpha$. Covering these ends
by the corresponding shadows $[\omega_m]_{\mathbb{F}}$ yields
\begin{equation}
 \dim_{\mathrm{H}}\Lambda_r(\alpha,\beta)\leq h_r(\alpha).
\end{equation}

A different idea in Lemma \ref{lem-assume-Hypothesis} is to replace $V(u)$ by a ``protected prefix".
Namely,  let $z(u)$
be the prefix of $V(u)$ with length $|z(u)|=\min_{v\in\mathcal{T}(u)}|V(v)|$.
Then we have $[z(u)]_{\mathbb{F}} \supset V([u]_{\mathcal{T}})$.
The price is that $z(u)$ may be much shorter than $V(u)$.
To control this loss, we decompose according to the suffix erased by the
descendant BRW. Hitting probability estimates and the pressure comparison assumption in
\eqref{eq-pressure-comparison} show that, replacing
$[V(u)]_{\mathbb{F}}$ by the larger  shadow
$[z(u)]_{\mathbb{F}}$ has no positive exponential cost.  Therefore we can
perform the covering argument with any speed $q\in[\alpha,\beta]$ and obtain
\begin{equation}
 \dim_{\mathrm{H}}\Lambda_r(\alpha,\beta)
 \leq\min_{q\in[\alpha,\beta]}h_r(q).
\end{equation}
For isotropic walks, \eqref{eq-pressure-comparison}  follows directly from convexity of the pressure function.

\smallskip
\noindent
\underline{Input 2: Level-set growth.}
Let $\mathcal{N}_{n,m}$ be the number of generation-$n$ particles at
distance $m$ from $e$, and let $\mathcal{N}_{n,m}^{\mathbb{F}}$ be the
number of distinct sites occupied by them.  From \eqref{eq:intro-RW}, we have
$ \E\bigl[\mathcal{N}_{n,m}\bigr]
 =e^{o(n)}\,a_{n,m}$ where $
 a_{n,m}:=r^n\exp \{-nL^*(m/n) \}$.
Lemmas~\ref{lem-level-set-size} and~\ref{lem-level-set-size-F} show  show that this exponential scale is
typically attained by both $\mathcal{N}_{n,m}$ and $\mathcal{N}_{n,m}^{\mathbb{F}}$.
Precisely,
for every  $\epsilon>0$ and all  large $n$,
\begin{equation}
 \P\bigl(\mathcal{N}_{n,m}^{\mathbb{F}}
       \leq a_{n,m}^{1-\epsilon}\bigr)
 \leq e^{-\sqrt n} \ ,
 \quad
 \inf_{g\in\mathbb{F}}
 \P_g\bigl(\mathcal{N}_{n,m}^{\mathbb{F}}(g)
       \geq a_{n,m}^{1-\epsilon}\bigr)
 \geq   {1}/{3},
 \label{eq-overview-level-set-input}
\end{equation}
uniformly over admissible $(n,m)$.  Here
$\mathcal{N}_{n,m}^{\mathbb{F}}(g)$ is the corresponding count in the forward cone $\mathbb{F}(g)$ rooted at $g$.  The same estimates hold for
$\mathcal{N}_{n,m}$.

To prove these estimates, Proposition~\ref{prop-LDP-sample-paths}
allows us to restrict attention to ancestral paths lying within a narrow
band around the linear profile, making a truncated second-moment argument
controllable.  A bootstrap then gives the lower-tail estimate for
$\mathcal{N}_{n,m}$.  Finally, Lemma \ref{lem-Nnx-is-small} yields the maximal multiplicity at any site among
the truncated particles is subexponential, so
$\mathcal{N}_{n,m}^{\mathbb{F}}$ has the same exponential growth rate.

\medskip
\noindent
\underline{Lower bounds.}
We explain the  lower bound for $\dim_{\mathrm{H}} \Lambda_{r}(\alpha,\beta)$,   the
proof for $E_{r}(\alpha,\beta)$ case is analogous and simpler.
Input~2 yields many distinct occupied sites at one scale, but does not
produce infinite rays along which this abundance persists.
We therefore
iterate it, retaining only particles that remain fecund at the
next stage.

\smallskip
\noindent
\emph{Fecund skeleton}.
We fix a sequence $(m_n)_{n\geq1}$ of spatial increments. For each
positive speed $q\in I(r)$, let $\ell_n(q)\approx m_n/q$ be the
corresponding genealogical increment, and choose the target branching
number $A_n(q)  $ slightly below the natural   growth scale, say $A_n(q) \approx  \exp \{(1-\epsilon)\ell_n(q)[\ln r-L^*(q)] \} $.
For  $v \in \mathcal{T}$,
let $\mathbb X_n(v,q)$ be the set of distinct sites lying $m_n$ steps
forward from $V(v)$ that are occupied by descendants of $v$ after
$\ell_n(q)$ generations. Call $v$ \emph{$q$-fecund for stage $n$} if
\[
   \#\mathbb X_n(v,q)\geq 6A_n(q).
\]
where the factor $6$ leaves room for the subsequent thinning.

Fix a sequence of speed labels
$\mathbf q=(q_n)_{n\geq1}$ and set
$\mathcal S_0=\{\emptyset\}$. Given $\mathcal S_{n-1}$, for each
$v\in\mathcal S_{n-1}$ and site $x\in\mathbb X_n(v,q_n)$, we choose one
representative particle $\operatorname{rep}_n(v,q_n,x)$ of this site and retain it
only if it is $q_{n+1}$-fecund for stage $n+1$.  The retained particles
form the next selected cohort:
\[
\mathcal S_n
:=
\bigcup_{v\in\mathcal S_{n-1}}
\bigl\{
\operatorname{rep}_n(v,q_n,x):
x\in\mathbb X_n(v,q_n),\
\operatorname{rep}_n(v,q_n,x)
\text{ is $q_{n+1}$-fecund for stage $n+1$}
\bigr\}.
\]
The genealogically linked sequence $(\mathcal S_n)_{n\geq0}$ is the
\emph{fecund skeleton associated with $\mathbf q$}.
By construction, distinct particles in the same selected cohort occupy distinct sites.
Consequently, once two skeleton branches separate, they determine
distinct spatial sites.

\smallskip
\noindent
\emph{The energy method.}
Let $B(\mathbf q)$ be the event that the root is $q_1$-fecund and, for
every $n\geq1$, each $v\in\mathcal S_{n-1}$ produces at least
$A_n(q_n)$ particles in the next selected cohort.
For a suitable choice of the scales, Input~2 and the branching property
give $\mathbb{P}(B(\mathbf q))>0$.

On $B(\mathbf q)$, we assign unit mass to the root and recursively
distribute the mass of each $v\in\mathcal S_{n-1}$ uniformly among the
particles in $\mathcal S_n$ descending from $v$. This yields a
probability measure $\nu^{\mathcal T}$ supported on the infinite rays
of the skeleton, with
\begin{equation} \label{eq:i-cylinder-mass}
     \nu^{\mathcal T}([u]_{\mathcal T})
     \leq
     \prod_{j=1}^n A_j(q_j)^{-1},
     \qquad u\in\mathcal S_n.
\end{equation}
Let
$  \nu^{\mathbb F}
   :=
   \nu^{\mathcal T} \circ (V|_{\partial\mathcal T})  ^{-1}$
be the pushforward of $\nu^{\mathcal{T}}$ to $\partial\mathbb F$.

For a given $[\alpha,\beta]\subset I(r)$, we choose $\mathbf q$ with
$\operatorname{dist}(q_n,\{\alpha,\beta\})\to0$ so that the sequence $ \frac{\sum_{j\leq n}m_j}{ \sum_{j\leq n}\ell_j(q_j) }  $  has liminf $\alpha$ and limsup $\beta$.
Then every infinite skeleton ray
has the same lower and upper speeds, which gives
\[
   \operatorname{supp}\nu^{\mathcal T}\subset E_r(\alpha,\beta),
   \qquad
   \operatorname{supp}\nu^{\mathbb F}\subset\Lambda_r(\alpha,\beta).
\]
The spatial separation of
the skeleton and the cylinder-mass bound \eqref{eq:i-cylinder-mass} give
\[
  \mathbf{I} (
 \theta,\nu^{\mathbb{F}}
 )
 := \int  \frac{\nu^{\mathbb{F}} (\dif \omega ) \nu^{\mathbb{F}} (\dif \omega' ) }{\mathrm{d}_{\partial \mathbb{F}}   ( \omega,  \omega'   )^{\theta}}   \,
   \lesssim
   \sum_{n\geq0}
   e^{\theta  \sum_{j=1}^{n+1} m _j}
   \prod_{j=1}^n A_j(q_j)^{-1}.
\]
Since
$\log A_n(q_n)\approx(1-\epsilon)m_nh_r(q_n)$, a sufficiently regular
choice of the scales, with
$  m_{n+1}=o (\sum_{j=1}^n m_j )$,
makes the preceding series convergent for every
$\theta<\min_{q\in[\alpha,\beta]}h_r(q)$, after choosing $\epsilon$
sufficiently small. The energy criterion therefore gives
\[
   \dim_{\mathrm H}\Lambda_r(\alpha,\beta)
   \geq\min_{q\in[\alpha,\beta]}h_r(q) \quad \text{ on } \ B(\mathbf{q}).
\]

To make the construction uniform, we discretize the speed labels at
stage $n$ by a finite grid $\mathbb D_n$, with
$\bigcup_n\mathbb D_n$ dense in $I(r)$, and construct the cohorts for
all label sequences simultaneously. A union bound yields a single event
of positive probability on which every interval admits a suitable
sequence. The branching property then provides arbitrarily many
independent attempts in descendant subtrees, upgrading this event to
probability one. Since re-rooting preserves asymptotic speeds and
Hausdorff dimensions, the lower bounds hold almost surely,
simultaneously for all $[\alpha,\beta]\subset I(r)$.

  \medskip
 \noindent\textbf{Notation convention.}
 Let $\mathbb{N}_0=\{0,1,\cdots \}$ denote the set of nonnegative integers. For a continuous function $f$ defined on a compact set, we denote by $\omega_{f}(\delta)$ its modulus of continuity. That is $\omega_{f}(\delta):= \sup_{|x-y| \leq \delta} |f(x)-f(y)|$. By the uniform continuity of $f$, we have $\lim_{\delta \to 0} \omega_{f}(\delta)=0 $.
We use notation  $C$ and $c$ for  positive constants whose actual values may vary from line to line. We write $C_{\delta}$ or $c_{\delta}$ if the constant depends on parameter $\delta$.

\section{Preliminaries}
\label{sec:ldp}

This section is organized as follows. First we clarify the definitions and notation.
In \S\ref{sec-free-group} we introduce  fundamental concepts related to free groups and the boundaries of trees.   In \S\ref{sec-RWonF} we recall the definition of random walks on free groups. In \S\ref{sec-BRWonF}, we precisely define the branching random walk on free groups and demonstrate that in the transient regime, the mapping \eqref{eq-BRW-mapping} is injective.
 In \S\ref{sec-Hdim} we recall the definition of the Hausdorff dimension.

 Next we introduce the essential results that will be frequently used later.
 \S\ref{sec-lltRW} presents the local limit theorems for random walks on free groups established by Lalley \cite{Lal91,L93}.
 \S\ref{sec-RW} discusses the large deviations for the word length of random walks on free groups and provides a new representation of the rate function \eqref{eq-forluma-Lstar}.
 \S\ref{sec-asymp-F} studies the asymptotic behavior of the free energy of  some Gibbs-Boltzmann distributions on free group.

\subsection{Free groups and trees}
\label{sec-free-group}
Let $\mathbb{F} \equiv \mathbb{F}^{d}$ be a free group of rank $d$ with symmetric generating set $\mathcal{A} = \left\{ a_1, \, \ldots, \, a_d, \, a_1^{-1}, \, \ldots, \, a_d^{-1} \right\}$, $d \geq 2$. That is, $\mathbb{F}$ consists of all finite reduced words from the alphabet $\mathcal{A}$ (a word is  reduced if no letter is adjacent to its inverse). Multiplication in $\mathbb{F}$ consists of concatenation followed by reduction, and the group identity $e$ is the empty word. The Cayley graph of $(\mathbb{F},\mathcal{A})$  is a $2d$-regular tree. (That is, we view each element $x \in \mathbb{F}$ as a vertex, and there is an edge between two vertices $x$, $y$ if and only if
$x^{-1}y \in \mathcal{A}$.) In this paper we shall not distinguish between group elements in $\mathbb{F}$ and vertices in the $2d$-regular tree to which they correspond.

Given an infinite tree $T$ with root $o$ (which may represent either a free group $(\mathbb{F}, e)$ or a Galton-Watson tree $(\mathcal{T}, \emptyset)$), we denote by $\dist(\cdot,\cdot)$ the standard graph distance on $T$ and write $|x|:= \dist (o,x) $ for $x \in T$.
If $T=\mathbb{F}$, we usually call $|x|$ the word length of   $x \in \mathbb{F}$; and if $T=\mathcal{T}$ we call $|u|$ the generation of $u \in \mathcal{T}$. For each $n \geq 0$, let
\begin{equation}
  \mathbb{F}_{n}:= \{ x \in \mathbb{F} : |x|=n\} \ \text{ and } \ \mathcal{T}_{n}:= \{ u \in \mathcal{T} : |u|=n\}.
\end{equation}

For   two vertices $x,y$ on the tree $T$, let $[x,y]$ denote the unique  geodesic from $x$ to $y$. There is
a natural partial order on $T$ called the genealogical order: $x \prec_{T} y$ if $x$ belongs to $[o,y]$,  in which case $x$ is called an ancestor of $y$. We say that $y $ is a child of $x$ if $x \prec_{T} y$  and $\dist(x,y)=1$.
 For $x$ and  $y$ in $T$, let   $x \wedge y $ denote the most recent common ancestor of $x, y$, which is   the unique   vertex $z \in T$ satisfying $[o,x] \cap [o,y]=[o,z]$. Define $T(x)$ the subtree of $T$ rooted at $x$ and $T_n(x)$ its $n$-th level:
\begin{equation}\label{eq:forward-cone}
     T(x):= \{ y \in T: x \prec_{T} y \} \quad , \quad  T_n(x):= \{ y \in T: x \prec_{T} y , \dist(x,y)=n\} .
\end{equation}
 We call $\mathbb{F}(x)$ the forward cone rooted at $x$.

A ray $\omega = (\omega_0 = o, \omega_1, \ldots)$ in $T$ is a semi-infinite self-avoiding path starting from the root.  The boundary   $\partial T$ of $T$ is defined to be the collection of all rays, endowed with the standard ultrametric distance:
\begin{equation}
  \mathrm{d}_{\partial T}(\omega, \omega') = e^{-|\omega \wedge \omega'|} \quad \text{for } \omega \neq \omega' \in \partial T.
  \end{equation}
 % Here as before,   $\omega \wedge \omega'$ represents the unique  vertex $z \in T$ such that $[o, z] = \omega \cap \omega'$.
 The corresponding topology on $\partial T$ is the topology of coordinatewise convergence.
It is known that for $\partial \mathbb{F}$ and $\partial \mathcal{T}$ equipped with the ultrametric distance,
\begin{equation}
\dim_{\mathrm{H}} \partial \mathbb{F} = \ln(2d - 1) \quad \text{and} \quad
\dim_{\mathrm{H}} \partial \mathcal{T} = \ln r \quad \text{a.s.}
\end{equation}
 For $x\in T$, define the shadow (or boundary cylinder) of $x$ by
\begin{equation}
  \label{eq:shadow}
  [x]_T:=\{\omega\in\partial T:\omega_{|x|}=x\}.
\end{equation}
It is a closed and open subset of $\partial T$ and satisfies
$\operatorname{diam}_{\partial T}([x]_T)\leq e^{-|x|}$.

The boundary $\partial T$ can also be  identified with   equivalent classes of infinite paths on $T$.  Precisely, let $\mathsf{Path}_{\infty}$ denote the set of all semi-infinite paths $\gamma$ in $T$ that go to infinity, i.e., $|\gamma_n| \to \infty$ as $n \to \infty$.  We say $\gamma,\gamma' \in  \mathsf{Path}_{\infty}$ are equivalent, denoted by  $\gamma \sim \gamma'$, if their  intersect with each other infinitely many times.
Notice that for any $\gamma \in \mathsf{Path}_{\infty}$, there exists a  unique ray $\omega$ in $T$ such that $\gamma \sim \omega$. That is, each ray $\omega$ can be regarded as  a representative of the equivalence class $\{ \gamma \in \mathsf{Path}_{\infty} :\gamma \sim \omega \}$. As a result,  $\partial T$ can be identified with the quotient space $\mathsf{Path}_{\infty} / \sim$.

% For $x\in T$, define the boundary cylinder
% \[
%   [x]_T:=\{\omega\in\partial T:\omega_{|x|}=x\}.
% \]
% It is a clopen subset of $\partial T$ satisfying
% \[
%   \operatorname{diam}_{\partial T}[x]_T\leq e^{-|x|}.
% \]

\subsection{Random walks  on free groups}
\label{sec-RWonF}
Let $\mu$ be a symmetric probability measure on $\mathcal{A} \cup \{ e\}$ such that $\mu(a) = \mu(a^{-1}) > 0$ for $a \in \mathcal{A}$.
We call $\{(Z_{n})_{n \geq 0},\mathbf{P}\}$ a random walk on  $\mathbb{F}$ with   step distribution $\mu$ (or $\mu$-random walk for short), if $(Z_{n})_{n \geq 0}$ is a Markov chain taking values in $\mathbb{F} $  with transition probabilities
\begin{equation}
\mathbf{P}\left(Z_{n+1}=g a  \mid Z_n=g\right)=  \mu(a) \,,   \text{ for }  a \in \mathcal{A} \cup \{e\}, n \geq 0.
\end{equation}
We denote $\mathbf{P}_{x}$ the law of $(Z_{n})_{n \geq 0}$ with initial value $Z_{0}= x$. For simplicity let $\mathbf{P}:= \mathbf{P}_{e}$.

The  transition kernel of $(Z_n)$ is denoted by
$ p_n(x,y) : = \mathbf{P}_{x}(Z_n= y)=\mathbf{P}(Z_n=x^{-1}y)$. Define
\begin{equation}
  \mathtt{h}_\mu=
  \begin{cases}
  1,&\mu(e)>0,\\
  2,&\mu(e)=0,
  \end{cases}
  \quad \text{ and } \quad
  \mathtt{L}_n^\mu
  :=
  \{0 \le m \le n: m \equiv  n   \pmod{\mathtt{h}_\mu}
\},  \label{eq:solve-periodicity-1}
\end{equation}
Then $p_n(x,y)>0 $ if and only if $|x^{-1}y|\in\mathtt{L}_n^\mu$.
 For $s \geq0$ and $n \ge 1$, define
\begin{equation}
  \label{eq:solve-periodicity-2}
  \langle s\rangle_{n,\mu}
  :=  \min\{ m \in \mathbb{N}_0 : |m-s| = \dist(s, n + \mathtt{h}_{\mu} \mathbb{Z} ) \}  .
\end{equation}
Then
$|\langle s\rangle_{n,\mu}-s|\leq1$,
$\langle s\rangle_{n,\mu}\in \mathsf{L}_n^\mu$ when $0\leq s\leq n$, whereas
$n\in\mathtt{L}_{\langle s\rangle_{n,\mu}}^\mu$ when $s\geq n$.

Given $x,y $ in $\mathbb{F}$, we define the Green function
\begin{equation}
  G(x,  y \mid r) :=\sum_{n=0}^{\infty} r^n p_n(x,y) \quad \text{for}\quad   r > 0.
\end{equation}
 Since the random walk is irreducible,  the radius of convergence $R$ of the Green function  $G(x,  y \mid r)$, given by
   \begin{equation}
   R^{-1} =  \limsup_{n \to \infty}  p_{n}(x,y)^{1/n},
  \end{equation}
 does not depend on $x,  y \in \mathbb{F}$ (see \cite[Lemma 1.7]{W00}).  And $R^{-1}$ is called the \textit{spectral radius} of the $\mu$-random walk. It is known that $G(x,y\mid R)< \infty$ for any $x,y \in \mathbb{F}$.

\subsection{Branching random walks on free groups}
\label{sec-BRWonF}

Let $\mathcal{T}$ be a Galton-Watson tree  rooted at $\emptyset$ with  offspring law $p=(p_{k})_{k \geq 0}$ and  mean offspring $r:=\sum_{k \geq 1} k p_{k} <\infty$.  We always assume that $p_0=0$, so that $\mathcal{T}$ is an infinite tree without leaves.

A branching random walk $\left( V(u), \, u \in \mathcal{T} \right)$ on free group $\mathbb{F}$ is a random map from the Galton-Watson tree $\mathcal{T}$ into $\mathbb{F}$ constructed as follows.
Conditionally on $\mathcal{T}$, let $(Y_{u}:u \in \mathcal{T})$ be a family of independent random element in $S$ with the common distribution $\mu$.   For every non-root vertex $u \in \mathcal{T}$, define
\begin{equation}
  V(u)= V(\emptyset) Y_{u_{1}} Y_{u_{2}} \cdots Y_{u_{n}}
\end{equation}
where $(\emptyset,u_{1},\cdots,u_{n}=u)=[\emptyset,u]$ is the geodesic in $\mathcal{T}$ from $\emptyset$ to $u$. We denote by $\mathbf{P}_{x}$ the law of the BRW starting from $x \in \mathbb{F}$, i.e., $\mathbf{P}_{x}(V(\emptyset)=x)=1$. Write  $\mathbf{P}:=\mathbf{P}_{e}$ for simplicity.

Let
$(\mathcal{F}^{V}_k := \sigma  \{V(w):  w \in \mathcal{T}, |w| \le k  \}, k \ge 0)$ be the natural filtration of the BRW.
For $w\in\mathcal{T}$, write $V^{w}$ for the descendant BRW rooted at $w$ and started from $V(w)$, that is:
  \begin{equation}
    V^{w}=( V(u), u \in \mathcal{T}(w))   \label{eq:sub-BRW}
  \end{equation}
The branching property yields that   conditionally on  $\mathcal{F}_k$,    $\{ V^{w} : w \in \mathcal{T}_{k}  \}$ are independent BRWs starting from their respective locations.

In the transient regime $r \in (1,R]$, since the BRW vacates every compact subset of the state space, it holds almost surely that
for every $t=(t_{n})_{n \geq 0} \in \partial \mathcal{T}$, $|V(t_{n})| \to \infty$ as $n \to \infty$; i.e., the path $(V(t_{n}) )_{n \geq 0}$ belongs to $\mathsf{Path}_{\infty}$.
For simplicity we shall not distinguish the $(V(t_{n}) )_{n \geq 0}$ and its equivalent class, so
 we will simply denote $(V(t_{n}))_{n \geq 0} \in \partial \mathbb{F}$.
 Then the branching random walk $V$: $\mathcal{T} \to \mathbb{F}$ can be extended to a continuous map $V$: $\mathcal{T} \cup \partial \mathcal{T} \to \mathbb{F} \cup \partial \mathbb{F}$,  by defining
\begin{equation}
    V: \partial \mathcal{T} \to  \partial \mathbb{F} ;  \quad t=(t_{n})_{n \geq 0} \mapsto  (V(t_{n}):n \geq 0) .
\end{equation}
Then the limit set $\Lambda_{r}$ is defined to be the image of $\partial \mathcal{T}$  under this map  $V$:
\begin{equation}
  \Lambda_{r} := \text{Image}(V|_{\partial \mathcal{T}}).
\end{equation}

Actually, $V: \partial \mathcal{T} \to \partial \mathbb{F}$ is an injection with probability one. To see this, it suffices to show that for any two rays $t, t' \in \partial \mathcal{T}$, the paths $(V(t_{n}))_{n \geq 0}$ and $(V(t'_{n}))_{n \geq 0}$ intersect at most finitely often. A special case of a general result by Hutchcroft (see \cite[Theorem 1.2 and Remark 1.3]{Hut20}) states that for two independent $\mu$-BRWs starting from $x$ and $y$ respectively, there are almost surely at most finitely many vertices in $\mathbb{F}$ that are visited by both BRWs. By the branching property, given  the particles of generation $1$ and their positions, the subsequent processes are independent BRWs with possibly different starting vertices; and hence they will intersect at most finitely many times. This argument establishes the injectivity of $V$.

\subsection{Hausdorff dimension}
\label{sec-Hdim} In this section we briefly  state the definition of the
Hausdorff measure and Hausdorff dimension   introduced by Felix Hausdorff in 1919. Let $X$ be a metric space. A $\delta$-cover of a set $F \subset X$ is a countable  collection of sets ${A_i}$ with diameters $\diam (A_{i})< \delta$ and $F \subset \cup_{i} A_{i}$.
Fix $s \geq 0$. For each $\delta>0$, we define
\begin{equation}
  \mathcal{H}_{\delta}^s (F)=\inf  \Bigl\{ \ \sum_{i} \diam( A_i)^{s}  : (A_i) \text { is a $\delta$-cover  of } F  \,  \Bigr\}
\end{equation}
As $\delta$ decreases, the class of $\delta$-overs of $F$ is reduced. Thus the infimum $ \mathcal{H}_{\delta}^s (F)$ increases as $\delta \downarrow 0$.   We define
\begin{equation}
  \mathcal{H}^s (F) := \lim_{\delta \downarrow 0}  \mathcal{H}_{\delta}^s (F) \in [0,\infty].
\end{equation}
 $\mathcal{H}^{s}(F)$ the s-dimensional Hausdorff measure of $F$.
Notice that if  $0 \leq s < s'$, and $\mathcal{H}^{s} (F) =0$, then   $\mathcal{H}^{s'} (F) =0$; and if $ \mathcal{H}^{s'}(F) =\infty$ then $\mathcal{H}^{s} (F) =\infty$. Thus we   define the Hausdorff dimension of the set $F$ by
\begin{equation}
  \dim_{\mathrm{H}} F := \inf \left\{ s \geq  0: \mathcal{H}^s (F) =0\right\}    =\inf \left\{ s \geq  0: \mathcal{H}^s (F) < \infty \right\} .
\end{equation}

In this paper, we show the upper bound of the Hausdorff dimension
by finding an appropriate cover. To demonstrate the lower bound, we apply  Frostman's well-known criterion.

\begin{lemma} \label{energy-method}
Let $\theta \geq 0$ and $\nu$ be a mass distribution on a metric space $(X,  d)$.  Define the $\theta$-potential of a point $x \in X$ with respect to $\nu$ as
 \begin{equation}
   \mathbf{I} (  \theta, \nu  ) := \int_{X \times X}  \frac{1 }{d(x , y)^{\theta}} \nu (\dif x)  \nu (\dif y) .
 \end{equation}
If  $ \mathbf{I} (  \theta, \nu  )<\infty$,  then we have  $\dim_{\mathrm{H}} X \geq \theta$.
\end{lemma}

% \begin{proof}
%   Suppose for $a<\infty$  that $v\left(B_a\right)>0$, where $B_a=\left\{x:\mathbf{I}(\theta, v, x)\leq a\right\}$, and define a probability measure $v_*$ by $v_*(A)=v\left(A \cap B_a\right) / v\left(B_a\right)$. Notice that for every $x \in B_{a}$, $ \int_{B_{a}} \frac{1}{d(x, y)^{\theta}}  \nu(\dif y) \leq \int_{X}  \frac{1}{d(x, y)^{\theta}}  \nu(\dif y) \leq a$. Hence
% \begin{align}
%   \mathbf{I} \left(\theta, v_*\right)& = \int_{X \times X} \frac{1 }{d(x , y)^{\theta}} \nu_{*}(\dif y) \nu_{*}(\dif x)  \\
%   &  = \frac{ 1}{\nu(B_a)^2} \int_{B_{a}}   v(\dif  x) \int_{B_{a}} \frac{1}{d(x, y)^{\theta}}  \nu(\dif y) \leq \frac{a}{\nu(B_{a}) } < \infty .
% \end{align}
% The usual Frostman Lemma (see \cite[Theorem 4.13]{Fal94}) states that the Hausdorff dimension of $E$ is at least $\theta$.
% \end{proof}

\subsection{Local limit theorems for RWs}
\label{sec-lltRW}
 Recall that $\{(Z_{n})_{n \geq 0}, \mathbf{P}\}$ represents the $\mu$-random walk on $\mathbb{F}$. Using the method of saddle-point approximations,
 Lalley \cite{Lal91} derived the precise asymptotic estimates of the $n$-step transition probabilities $p_n(e, x)$ of the random walk $(Z_{n})$.
 To state his results, we first introduce some necessary notation.

For each fixed $x \in \mathbb{F}$, let $T(x):=\inf \{ n \geq 1: Z_{n}=x \}$. Define the  generating function of $T(x)$ by
\begin{equation}\label{eq:generating-fun-F}
  F_{x}(r) :=  \mathbf{E}  \bigl[    r^{ T(x)} 1_{\{ T(x) < \infty \} } \bigr]  = \sum_{n=0}^{\infty}  r^{n} \mathbf{P}( T(x)= n) \ , \quad  r \geq 0 .
\end{equation}
Let $\psi(s)$ denote the vector $(\psi_{a}(s))_{a \in \mathcal{A}}$, where $\psi_{a}(\cdot)$ is the logarithmic moment generating function given by
\begin{equation}\label{eq-def-psi-a-s}
   \psi_{a}(s)   := \ln F_{a}(e^{s}) = \ln \mathbf{E}   \bigl[  e^{ sT(a)} 1_{\{ T(a) < \infty \} }  \bigr]  \ , \quad   s \in \mathbb{R}.
\end{equation}
Lalley \cite[Proposition 1,2]{Lal91} showed that for any $a \in \mathcal{A}$, $\psi_{a}(s)$ is a strictly increasing and strictly convex function  of $s \in(-\infty, \ln R]$ with the following properties:
\begin{equation}\label{eq-properties-psi}
\psi_{a}^{\prime \prime}(s)>0  \ \forall s < \ln R , \
\lim _{s \uparrow \ln R} \psi_{a}^{\prime}(s)=\infty    \ ,\
\lim _{s \downarrow-\infty} \psi_{a}^{\prime}(s)=1   \ , \
 \text{ and } \, \psi_{a}(\ln R) = \ln F_{a}(R)
 < 0.
\end{equation}

 Set $\Omega := \{ \xi=\left( \xi_{a} \right)_{ a \in \mathcal{A} } : \xi_{a} \geq 0,  \, \sum_a \xi_a \leq 1\}$.
 For each $x \in \mathbb{F}$,  let $\Xi(x) =\left(\Xi_{a}(x) \right)_{a \in \mathcal{A}}$ denote the vector of nonnegative integers where   $\Xi_{a}(x)  $ is the number of times  the letter $a$ occurs   in the reduced word representation of $x$. Let $\xi(n,x) := \frac{1}{n} \Xi(x) \in  \Omega$ for all $n \geq |x|$ and $\xi(x):= \frac{1}{|x|} \Xi(x)$.  We introduce a function $\hyperref[eq-def-Psi*]{\Psi^{*}}$ defined on $\Omega$ by
\begin{equation}\label{eq-def-Psi*}
\Psi^{*}(\xi) :=\inf _{-\infty<s \leq \ln R}  \Bigl( \  \sum_{a \in \mathcal{A}} \xi_{a} \psi_{a}(s)-s \, \Bigr) \quad  \text{ for } \xi \in \Omega.
\end{equation}
By \eqref{eq-properties-psi}, we have  $\Psi^{*}(\xi) \leq \sum_{a} \xi_{a} \psi_{a}(\ln R)-\ln R \leq -\ln R$ for all $\xi \in \Omega$.
Furthermore, if $\sum \xi_{a}=0$ (i.e., $\xi= 0$), the infimum is attained at $s=\ln R$,  yielding $\Psi^{*}(0)=-\ln R$. If $0<\sum \xi_{a}<1$, the infimum is attained at the unique $s=s  ( \xi ) \in(-\infty, \ln R)$ satisfying
\begin{equation}\label{eq-s-xi}
\sum_{a \in \mathcal{A}} \xi_{a} \psi_{a}^{\prime}( s  ( \xi ) )=1 .
\end{equation}
 If $\sum \xi_{a}=1$, the infimum is attained as $s=-\infty$,    in which case $
\Psi^{*}(\xi)=\sum_{a \in \mathcal{A}} \xi_{a} \ln \mu(a)$ because $(\psi_{a}(s) - s)=\ln ( F_{a}(e^{s})/e^{s} ) \to \ln \mu(a)$ as $s \to -\infty$.

\medskip
The main results of Lalley \cite{Lal91} read as follows.  
\footnote{Lalley assumes $\mu(e)>0$ to ensure aperiodicity and explicitly
states that his results readily extend to the case $\mu(e)=0$. In the latter case, the estimates are understood on the admissible parity class $|x|\in\mathtt L_n^\mu$, equivalently, $n-|x|$ is even.}
\begin{proposition}[{\cite[Theorem 1, Proposition 4-5]{Lal91}}] \label{T:LLT}
For $n \geq 1$,  there exists function $\beta_{n}$ on $\Omega$ such that
for any $x \in \mathbb{F}$   with $|x| \in \mathtt{L}^{\mu}_n$
\begin{equation}\label{eq-asymptotic-pn(e,x)}
   \mathbf{P}(Z_{n}=x )  = \beta_n \left( \xi(n,x) \right) \exp  \bigl\{  n \Psi^{*}  (  \xi(n,x)  )  \bigr\}   ,
\end{equation}
and the following assertions hold:
\begin{enumerate}
  \item \label{def-widetilde-beta-n} $\hyperref[def-widetilde-beta-n]{\widetilde{\beta}}_{n}:= \max_{1 \leq k \leq n} \sup_{\xi \in \Omega}  |\ln \beta_k( \xi )| = o(n) $ as $n \to \infty$.
  \item There exists a constant $C$ such that  for all $n \geq 1$  and $ \xi \in \Omega$, $ \beta_n( \xi ) \leq C$.
  \item There exists a constant $C$  for all $n \geq 1$,  $ \xi \in \Omega $ with $  \|\xi \|_{1} \leq n^{-1/4}$,  $ \beta_n( \xi ) \leq  \frac{C}{n^{3/2}} (1+ n \|\xi\|_{1})  $.
\end{enumerate}

\end{proposition}

\begin{remark}
  Lalley  \cite[Propositions 4-5]{Lal91} also provided  the explicit asymptotic estimates of $\beta_{n}(\xi)$.  As he observed,  ``These saddle-point approximations are not entirely routine, because $\beta_n(\xi)$ makes a transition from $\mathrm{C}n^{-3 / 2}$ to $C^{\prime} n^{-1 / 2}$ as $\sum \xi_a$ varies from 0 to $\epsilon$, then another transition from $C^{\prime \prime} n^{-1 / 2}$ to $C^{\prime \prime \prime}$ as $\sum \xi_a$ varies from $1-\epsilon$ to 1. ''
\end{remark}

By  Proposition \ref{T:LLT},
$\Psi^{*}$ serves as the rate function  for the large deviation probabilities  $\mathbf{P}(Z_{n}=x)$.
The following lemma shows that $\Psi^{*}$ is  concave on the interior points of  its domain. The proof is included in Appendix~\ref{app-pf-concave}.

\begin{lemma} \label{l:concave}  The rate function $ \Psi^{*}(\xi) $ is concave on  $\Omega$ and the following assertions hold.
  \begin{enumerate}
   \item For each $\xi \in \Omega$ with $\|\xi\|_{1}<1$, $\nabla \Psi^{*}(\xi) = \psi(s(\xi))$. In particular
   $\nabla \Psi^{*}( 0) =\psi(\ln R) $.
   \item For all $\xi \in \Omega \setminus \{\mathbf{0}\} $, $\lambda \mapsto \Psi^{*}(\lambda\xi) $ is strictly decreasing in $\lambda \in [0,1]$.
   \item Let $\xi \in \Omega$ be an inner point and let $h=h(\xi)$ be orthogonal to $(\psi_{a}'(s(\xi)))_{a}$. Then $\Psi^{*}$ is linear on the segment $\Omega  \cap \{ \xi + t h(\xi) : t \in \mathbb{R}\}$.
  \end{enumerate}

  \end{lemma}

\subsection{Large deviations for word length of RWs}
\label{sec-RW}
Guivarc'h \cite{G80} proved that the random walk  $(Z_{n})_{n \geq 1}$ on the free group $\mathbb{F}$ has a rate of escape  $\mathtt{C}_{\mathrm{RW}}$. Specifically, the following  strong law of large numbers holds:
  \begin{equation}\label{eq-LLN}
    \lim_{n \to \infty} \frac{\left| Z_n \right|}{n} = \lim_{n \to \infty} \frac{\mathbf{E} \left[ \left| Z_n \right| \right]}{n} = \mathtt{C}_{\mathrm{RW}} >0 \  \text{ a.s. } \footnote{This result is now a consequence of  Kingman's subadditive ergodic theorem; and the positivity of  $\mathtt{C}_{\mathrm{RW}} $ can be shown  by comparing $|Z_{n}|$ with a biased RW on half-line $\mathbb{Z}_{+}$. Or see  \cite[Chaper 8]{W00}.}
  \end{equation}
  Sawyer and Steger \cite{SS87}
  derived a central limit theorem (and a law of the iterated logarithm)  for the sequence $(|Z_{n}|)$. They showed that there exists $\sigma^2 \geq 0$
   such that
  \begin{equation}\label{eq-CLT-word-length}
    \frac{ |Z_{n}|- \mathtt{C}_{\mathrm{RW}} n}{\sqrt{n}} \overset{law}{\longrightarrow}  N(0,\sigma^{2} ) \text{ as } n  \to \infty
  \end{equation}
  Once the SLLN and CLT are established, the next step is to investigate the large deviation behaviors.
  Lalley \cite{L93} obtained a  precise result on large deviations for  $(|Z_{n}|/n)$ where the step distribution $\mu$ is allowed to have a finite-range support.
  Specifically, the rate function $\hyperref[eq-L*-via-P]{L^{*}}$ is given by the the following Legendre transform
\begin{equation}\label{eq-L*-via-P}
  L^{*}(q) :=\sup _{-\infty<s \leq \ln R}\{ s- q P(s)   \} .
\end{equation}
  of the pressure function (or the logarithm of the Perron-Frobenius eigenvalues) $P(s)$, $s \in (-\infty,
  \ln R]$ of certain Ruelle operators (see \cite[Section 7]{L93}). Furthermore $P(0)=0, P''(s) > 0$ and $ P'(s) \geq 1$ for $ s < \ln R$ and $\lim_{s \to \ln R} P'(s)= \infty$.

The supremum in \eqref{eq-L*-via-P} is attained uniquely at some $s=s(q) \in$ $[-\infty, \ln R]$. If $q=0$, then $s(0)=\ln R$; if $q \in\left(0, 1\right)$, then $s(q) \in$ $(-\infty, \ln R)$; and if $q \geq 1$, then $s(q)=-\infty$. Note that $s(q)=s \in(-\infty, \ln R)$ if and only if $P'(s)=1 / q$. Clearly, $s(q)$ is a decreasing function. By the chain rule,
\begin{equation}\label{eq-L*prime}
  (L^{*})'(q)= - P(s(q)) \quad \text { and } \quad   (L^{*})''(q)=\frac{P'(s(q))^2}{q P''(s(q))}>0.
\end{equation}
So $L^{*}(q)$ is strictly convex on $\left(0, 1\right)$, $(L^{*})'(0)=-P(\ln R)$  is finite but $(L^{*})'(1)=\infty$.
Clearly $L^{*}(q) \geq 0$, and $L^{*}$ attains it minimum value  $0$ uniquely at $q= 1/P'(0)$.

\begin{proposition}[{\cite[Theorem 7.2]{L93}}]\label{prop-word-length-ldp}
There are positive constants $C(q)$, $D(q)$, $q \in (0, 1)$, such that as $n\to \infty$ and  uniformly for $m \in \mathtt{L}^{\mu}_n$, $m/n$ in any compact subset of $(0, \, 1)$,
  \begin{equation}
    \mathbf{P}( |Z_n|=m ) \sim \frac{  C\left( {m}/{n} \right) }{ \sqrt{2\pi m D(m/n)} }e^{ -n L^{*}( \frac{m}{n}) }.
  \end{equation}
\end{proposition}

Proposition \ref{prop-word-length-ldp} implies the law of large numbers \eqref{eq-LLN}, giving us that $\mathtt{C}_{\mathrm{RW}} $  equals $1/P'(0)$, and  the central limit theorem \eqref{eq-CLT-word-length}, confirming that  $\sigma^2 = (L^{*})''(\mathtt{C}_{\mathrm{RW}}) >0$. It also establishes the following large deviation principle:
For any Borel subset $B \subset \mathbb{R}$,
  \begin{equation}\label{eq-ldp-word-length}
    - \inf_{q \in B^o} L^{*}(q) \leq \liminf_{n \to \infty} \frac{1}{n} \ln \mathbf{P}  \Bigl(\frac{|Z_{n}|}{n} \in B  \Bigr) \leq \limsup_{n \to \infty} \frac{1}{n} \ln \mathbf{P} \Bigl( \frac{|Z_{n}|}{n} \in B  \Bigr) \leq - \inf_{q \in \overline{B}} L^{*}(q),
  \end{equation}
 where $B^o$ and $\overline{B}$ denote the interior and the closure of $B$, respectively. See also Corso \cite{C21} for a generalization of \eqref{eq-ldp-word-length} for random walks on free products of finitely generated groups.

 The goal of this section is to prove that uniformly in $m \in \mathtt{L}^{\mu}_{n}$, as $n \to \infty$,
 \begin{equation}\label{eq-Z-n-is-m-0}
  \mathbf{P}( |Z_{n}|= m ) =  \exp \left(- n L^{*}(m/n) +  o(n)  \right).
 \end{equation}
 We emphasize that the uniformity of  \eqref{eq-Z-n-is-m-0} for all $m \in \mathtt{L}^{\mu}_{n}$ is crucial in the proofs of Theorems \ref{thm-Hdim-Lambda-alpha-beta} and \ref{thm-dim-E-alpha-beta} in the critical case $r=R$ and $\alpha=0$. In contrast,   Proposition \ref{prop-word-length-ldp} applies only when $m/n$ in an arbitrarily fixed compact subset of $(0,1)$. We have to admit that our method can not give the explicit order of $O(\log n)$ as in Proposition \ref{prop-word-length-ldp}, from which one can deduce the central limit theorem.

 We also obtained additional  expressions for the rate function $L^{*}$   via  the function $\Psi^{*}$ defined in \eqref{eq-def-Psi*}, see \eqref{eq-forluma-Lstar} and \eqref{eq-def-f(qxi)}. It is not surprising that one can estimate  $\mathbf{P}(|Z_{n}|=m)$ if the asymptotic of $\mathbf{P}(Z_{n}=x)$ was known. However, to our knowledge, the expression  \eqref{eq-forluma-Lstar} %and \eqref{eq-forluma-P}
 has not previously appeared in the literature.

 Before the main result  of this section, we introduce notation and a supporting lemma. Let  $\Omega^{1} := \{ \xi =(\xi_{a})_{a \in \mathcal{A}}: \xi_{a} \geq 0,    \sum_a \xi_{a}=1   \}$ be the set of  probability measures on $\mathcal{A}$.
 For $n \geq 1$, set $\Omega^{1}_{n}: = \{ \xi \in \Omega^{1}:   n\xi_a \in \mathbb{N}_0 \text{ and } [\xi_a+\xi_{a^{-1}}]\ind{\xi_a\xi_{a^{-1}}>0}<1, \forall a\in\mathcal{A} \} $.
 For each vector $\lambda = (\lambda_{a})_{a \in \mathcal{A}} \in \mathbb{R}^{2d}$, define
 \begin{equation}
  \varrho(\lambda)   := \text{ the largest positive solution  to }\sum_{a\in \mathcal{A}} e^{\lambda(a)}\frac{  e^{\varrho} - e^{ \lambda(a^{-1})} }{ e^{2\varrho} -e^{\lambda(a) }e^{ \lambda(a^{-1}) } } =1  .
   \label{eq-def-of-varrho}
 \end{equation}
(Indeed  $ \varrho(\lambda)$ is the logarithm of the leading eigenvalue of the matrix $ \left( e^{\lambda_a}  \ind{b \neq a^{-1}} \right)_{(a, b) \in \mathcal{A}} $  by Lemma \ref{lem-PF-eigen-equa}.)
Let  $\varrho^{*}$ denote the Legendre transform of $\varrho$,  defined by
 \begin{align}
  \varrho^{*} (\xi)   :=   \sup_{ \lambda  \in \mathbb{R}^{2d}} \left\{ \langle \xi, \lambda \rangle-  \varrho(\lambda)  \right\} + \ln (2d-1)  \ , \quad    \forall \, \xi \in \mathbb{R}^{2d}.
 \end{align}
 According to \cite[Theorem 3.1.2, 3.1.6]{DZ09}, $\varrho^{*}$ is a convex, good rate function and  $\varrho^{*}(\xi) < \infty$ if and only if $\xi \in \Omega^{1}$.
 Since  Legendre transform is an involution (by the Fenchel--Moreau theorem),
we have for any $\lambda \in \mathbb{R}^{2d}$,
\begin{equation}\label{eq-rho-via-rho*}
  \varrho(\lambda) =  \max_{\xi \in \Omega^{1}} \big\{ \langle  \xi,\lambda \rangle - \varrho^{*}(\xi) + \ln (2d-1)  \big\} .
\end{equation}

The following lemma gives a specific value of $\varrho$ at $\ln R$ and follows directly from the crucial identity in \cite{HL00} that  $ \sum_{a \in \mathcal{A}}\frac{F_{a}(R)^{2} }{1 + F_{a}(R)^{2}} = 1$.

\begin{lemma}[{\cite[Proposition 3]{HL00}}]\label{lem-rho2psi}
  $\varrho(2 \psi(\ln R))=0$.
\end{lemma}

 The  next Lemma \ref{lem-ldp-empirical-measure} is based on the classical large deviation theory of empirical measures for Markov chains. The estimate
$ \mathbf{P} \left( \xi(X_n)  = \xi  \right)=n^{O(1)} e^{- n \varrho^{*}(\xi)}$ which holds uniformly in  $\xi$,  appears to be less readily available in the literature,  
 so we provide a proof in Appendix \ref{app-pf-empirical-measure} for completeness.

 \begin{lemma}\label{lem-ldp-empirical-measure}
 Let $X_{n}$ be a random variable uniformly distributed on $\mathbb{F}_{n}$.
 There exists a function $B_{n}$ on $\Omega^{1}_{n}$ such that  $  \sup_{ \xi \in  \Omega^{1}_n } |\ln B_{n}( \xi)| =O(\ln n ) $ and
 \begin{equation}
   \mathbf{P} \left( \xi(X_n)  = \xi  \right)= B_{n}( \xi ) e^{-n \varrho^{*}(\xi) } \ \text{ for all } \xi \in \Omega^{1}_{n}.
 \end{equation}
 \end{lemma}

 Define the function $f : [0,1] \times \Omega^1 \to \mathbb{R}$ by
 \begin{equation}\label{eq-def-f(qxi)}
   f (q,  \xi) := q [ \varrho^{*} (\xi) - \ln(2d-1) ]  - \Psi^{*}(q \xi)   \quad \text{ for }q \in [0,1], \xi \in \Omega^{1}.
 \end{equation}
 Then $f$ is  continuous since $ \varrho^{*}$ is a convex good rate function on  $\Omega^{1}$.
 The next Proposition shows the assertion \eqref{eq-Z-n-is-m-0} and gives an expression for $L^{*}$ via the function $f$.

 \begin{proposition}\label{prop-LDP-|Zn|}
 For $n \geq 1$ and $m \in \mathtt{L}^{\mu}_{n}$, there exists   $B_{n}(m)>0$ such that
    \begin{equation}
     \mathbf{P}( |Z_n|=m ) =    B_{n}(m) \exp \bigl\{   - n L^{*}(m/n)   \bigr\}
    \end{equation}
    and the following assertions hold.
    \begin{enumerate}[(i)]
  \item \label{def-widetilde-B-n}  As $n \to \infty$, we have $\hyperref[def-widetilde-B-n]{\widetilde{B}}_{n}:= \max_{1 \leq k \leq n} \max_{m \in \mathtt{L}^{\mu}_{k}} | \ln B_{k}(m) | = o(n)$.
  \item  There is a function $B (m)$  such that $| \ln B(m) | =O(\ln m) $ as $m \to \infty$ and  $ B_{n}(m) \leq  B(m) \sup_{ \xi \in \Omega^1} \beta_{n}(\frac{m}{n}\xi) \leq C B(m)$ where  $C$ is the constant given in Proposition \ref{T:LLT}.
  \item For each $q \in [0, 1]$, with the function $f(q, \xi)$   defined in \eqref{eq-def-f(qxi)}, the following holds:
  \begin{equation} \label{eq-forluma-Lstar}
    L^{*}(q) = \min_{\xi \in \Omega^{1}} f(q, \xi).
  \end{equation}
  As a result, $L^{*}(0)=-\Psi^{*}(0)= \ln R$, and
  \begin{equation}\label{eq-forluma-P}
    P(s)=\varrho (\psi(s)) \ \forall s \in (-\infty,\ln R] \quad \text{ and }  P(\ln r)= \dim_{\mathrm{H}} \Lambda_{r} \ \forall r \in (1,R].
  \end{equation}
    \end{enumerate}
  \end{proposition}

\begin{proof}[Proof of Proposition \ref{prop-LDP-|Zn|}]
 We first consider the case $m=0$. By \eqref{eq-def-Psi*} and \eqref{eq-L*-via-P}, we have $-\Psi^{*}(0)=\ln R=L^{*}(0)$. If $0\in\mathtt{L}^{\mu}_{n}$, Proposition \ref{T:LLT}, applied with $x=e$, yields
\begin{equation}
  \label{eq:m=0case}
 \mathbf{P}(|Z_n|=0)=\beta_n(0)e^{n\Psi^{*}(0)}=\beta_n(\mathbf{0})e^{-nL^{*}(0)}
\end{equation}
with $ \beta_n(\mathbf{0}) = o_n(1)$. Setting $B_n(\mathbf{0})= \beta_n(\mathbf{0})$, it thus suffices to prove assertions (i) and (ii) for   $m \in \mathtt{L}^{\mu}_{n}$ with $m\ge 1$.
Proposition \ref{T:LLT} yields
  \begin{equation}\label{eq-basic-sum-P-Zn=m}
    \mathbf{P}( |Z_n|=m ) = \sum_{x \in \mathbb{F}_{m}} \mathbf{P}( Z_n=x ) = \sum_{x \in \mathbb{F}_{m}} \beta_n  \bigl(  \frac{m}{n} \xi(x)  \bigr)  \exp \Bigl\{   n \Psi^{*}  \bigl( \frac{m}{n} \xi(x)  \bigr)   \Bigr\} .
  \end{equation}

  \smallskip
\noindent
\underline{\textit{Step 1}}.
Our main objective is to analyze  the sum  $  \sum_{x \in \mathbb{F}_{m}}  e^{ n \Psi^{*} ( \frac{m}{n} \xi(x)  )}$. To this end, we introduce   $(X_{m},\mathbf{P})$, a random variable with uniform  distribution on $\mathbb{F}_{m}$. Then,  we have
\begin{align}
	\qquad \Sigma_{\refeq{eq-sum-exp-psi-00}}(m,n) &:= \sum_{x \in \mathbb{F}_{m}}   \exp  \Bigl\{   n \Psi^{*} \bigl(  \frac{m}{n} \xi(x)  \bigr) \Bigr\}   \label{eq-sum-exp-psi-00} =  \# \mathbb{F}_{m}   \cdot    \mathbf{E}  \Bigl[     \exp \Bigl\{   n \Psi^{*} \bigl(   \frac{m}{n}  \xi(X_m)  \bigr)    \Bigr\}\Bigr]   \\
	&=  \frac{2d}{2d-1} (2d-1)^{m}  \sum_{ \xi \in \Omega_{m}^{1}}   \exp \Bigl\{ n \Psi^{*} \bigl(   \frac{m}{n} \xi  \bigr)   \Bigr\}   \mathbf{P} (   \xi(X_m) 	 =\xi ).
\end{align}
 Lemma \ref{lem-ldp-empirical-measure}   yields $
  \mathbf{P}  (   \xi(X_m) =\xi )= B_{m}( \xi ) \, e^{-m  \varrho^{*}( \xi ) }
  $.  Plugging this into the expression above, since $m \ln(2d-1) + n \Psi^{*}  (  \frac{m}{n} \xi   )   - m  \varrho^{*}( \xi)  = - n f   (  \frac{m}{n},  \xi  )$ by
    \eqref{eq-def-f(qxi)},    we obtain
\begin{equation}
  \Sigma_{\refeq{eq-sum-exp-psi-00}}(m,n)
   =  \frac{2d}{2d-1} \sum_{ \xi \in \Omega_{m}^{1}} B_{m}( \xi )  \exp \Bigl\{   - n f  \bigl(  \frac{m}{n},  \xi   \bigr) \Bigr\}    .
\end{equation}
Define $B_{+}(m) :=2 \frac{2d}{2d-1}  (m+1)^{2d}    \sup_{ \xi \in \Omega^{1}_m } B_{m}( \xi)$ and $B_{-}(m) :=  \frac{2d}{2d-1} \inf_{ \xi \in \Omega^{1}_m }B_{m}( \xi )$.   Lemma  \ref{lem-ldp-empirical-measure} shows
 $  | \ln B_{\pm} (m) | =O(\ln m) $ as $m \to \infty$. Since   $\# \Omega_{m}^{1} \leq (m+1)^{2d}$, we get
\begin{equation}
B_{-}(m) \exp \Bigl\{   - n\min_{\xi\in \Omega^{1}_{m}} f\bigl( \frac{m}{n} ,  \xi  \bigr) \Bigr\}   \leq  \Sigma_{\refeq{eq-sum-exp-psi-00}}(m,n) \leq   B_{+}(m) \exp \Bigl\{   - n\min_{\xi\in \Omega^{1}_{m}} f\bigl( \frac{m}{n} ,  \xi  \bigr) \Bigr\}.
\end{equation}

  \smallskip
\noindent
\underline{\textit{Step 2}}.
We show  the error $|\min_{\xi\in \Omega^{1}_{m}} f(\frac{m}{n},  \xi )-\min_{\xi\in \Omega^{1}} f(\frac{m}{n},  \xi )|$ is negligible.
 Specifically,
\begin{equation}
 \Delta_{f}(m):= \sup_{q \in [0,1]} \Bigl|\, \min_{\xi \in \Omega^{1}_{m}} f(q,  \xi) -  \min_{\xi\in \Omega^{1}} f( q,  \xi )  \,\Bigr|   \to 0 \ \text{ as } m \to \infty.
\end{equation}
Since  $f$ is continuous, there must exist $\xi(q)\in \Omega^{1} $ such that
 $\min_{\xi\in \Omega^{1}} f(q,  \xi) = f(q,    \xi(q))$. Let $\xi^m(q)$ be the closest point  to $\xi(q)$  in $\Omega^{1}_{m}$. Denote by $\omega_{f}$  the modlous of continuity of $f$.  $  \Delta_{f}(m)  $ is then bounded above by
\[     \sup_{q \in [0,1]}  | f(q, \xi^m(q))-  f(q, \xi(q))|
     \leq  \sup_{q \in [0,1]} \omega_{f}( \|\xi^m(q)-\xi(q)\|_{1} ) \leq \omega_{f}(2d/m) = o_m(1) \]
 Plugging this into the last display in Step 1 gives
\begin{equation}\label{eq-goal-bound}
 B_{-}(m) e^{-n \Delta_{f}(m)}    \leq   \Sigma_{\refeq{eq-sum-exp-psi-00}}(m,n) \exp \Bigl\{  n \min_{\xi\in \Omega^{1}} f(\frac{m}{n} , \xi ) \Bigr\} \leq   B_{+}(m)    .
\end{equation}

 Since $\Delta_f(m)$ may not vanishes for $m=O(1)$, we provide a straightforward lower bound for $\Sigma_{\refeq{eq-sum-exp-psi-00}}(m,n)$
for small  values of $m$, say $1 \leq m \leq \sqrt{n}$.
From \eqref{eq-def-f(qxi)} and the fact  $\varrho^{*}(\xi) \ge 0$,  we have  $f(q,\xi) \geq -\max_{\xi \in \Omega^{1}}\Psi^{*}(q\xi)-q\ln(2d-1)$. It thus holds, with  $\omega_{\Psi^{*}}$ denoting the modulus of continuity of $\Psi^{*}$, that
\begin{align}
  \Sigma_{\refeq{eq-sum-exp-psi-00}}(m,n)
  &\geq  \exp \Bigl\{  n    \max_{\xi \in \Omega^{1}_{m} } \Psi^{*} ( \frac{m}{n} \xi  )  \Bigr\}   \geq  \exp \Bigl\{  n    \max_{\xi \in \Omega^{1}  } \Psi^{*} ( \frac{m}{n} \xi  ) - n \, \omega_{\Psi^{*}}(n^{-\frac{1}{2}})   \Bigr\} \notag \\
  & \geq  \exp \Bigl\{  -n \min_{\xi\in \Omega^{1}} f(\frac{m}{n},  \xi )  - n [ \omega_{\Psi^{*}}(n^{-\frac{1}{2}}) + n^{-\frac{1}{2}} \ln(2d-1) ]  \Bigr\} .
    \label{eq-goal-bound-2}
\end{align}

  \smallskip
\noindent
\underline{\textit{Step 3}}.
Define    $\Delta(n)$ as follows:
\begin{equation}
 \Delta(n):= \sup_{\sqrt{n} \leq m \leq n} \Bigl\{   \Delta_{f}(m)+\frac{1}{n}|\ln B_{-}(m)| \Bigr\} +\omega_{\Psi^{*}}(n^{-\frac{1}{2}})+n^{-\frac{1}{2}}\ln(2d-1).
\end{equation}
It clear that  $\Delta(n) \to 0$ as $n \to \infty$ since  $  | \ln B_{\pm} (m) | =O(\ln m) $ and $\Delta_{f}(m) \to 0$ as $m \to \infty. $

Combining \eqref{eq:m=0case}  \eqref{eq-basic-sum-P-Zn=m}
  \eqref{eq-goal-bound} and \eqref{eq-goal-bound-2}, and using that
$f(0,\xi)=-\Psi^{*}(\mathbf{0})=\ln R$, we conclude that, for every
$m\in\mathtt{L}^{\mu}_{n}$,
\begin{align}
	\mathbf{P}( |Z_n|=m )  =    B_{n}(m) \exp  \Bigl\{    - n \min_{\xi\in \Omega^{1}} f \bigl( \frac{m}{n},   \xi   \bigr)   \Bigr\}   \label{eq-ldp-Zn=m-pre-00}
\end{align}
where  prefactors $B_{n}(m)$  satisfy $ e^{-n \Delta(n)} \inf_{ \xi \in \Omega} \beta_{n}( \xi)   \leq  B_{n}(m) \leq B_{+}(m) \sup_{ \xi \in \Omega^1} \beta_{n} \bigl( \frac{m}{n}\xi \bigr)   $.
Since Proposition \ref{T:LLT} gives $  \sup_{\xi \in \Omega}  |\ln \beta_n ( \xi )|  =o(n) $    while $\Delta(n)\to0$ and  $  | \ln B_{\pm} (m) | =O(\ln m) $,
it follows that $ \max_{  m \in \mathtt{L}^{\mu}_{n}  }|\ln B_{n}(m)| = o (n)$. This implies assertion (i).

Set  $B(m) :=  B_{+}(m)$.  By Proposition \ref{T:LLT} there exists a constant  $C>0$ such that $\beta_{n}(\xi) \leq C$ for every $n \ge 1$ and $ \xi \in \Omega$. This shows $B_{n}(m) \leq C B(m)$ and hence assertion (ii) follows.

 \smallskip
\noindent
\underline{\textit{Step 4}}.
Comparing formula \eqref{eq-ldp-Zn=m-pre-00} with \eqref{eq-ldp-word-length}, it follows  that  \eqref{eq-forluma-Lstar} holds for  all $q \in [0,1]$.
It remains to prove  $P(s)= \varrho(\psi(s))$. According to the definition of $\Psi^{*}$, we have
\begin{align}
  - L^{*}(q) &= \max_{\xi \in \Omega^{1}} \left\{  \Psi^{*}(q \xi) + q[ \ln (2d-1) - \varrho^{*}(\xi) ]  \right\} \\
  &= \max_{\xi \in \Omega^{1}} \inf_{-\infty<s \leq \ln R} q \langle   \xi, \psi(s) \rangle - s +   q[ \ln (2d-1) - \varrho^{*}(\xi) ] \label{eq-max-inf-L*}
\end{align}

For fixed $q\in[0,1]$, define
\[
  \Phi_q(s,\xi)
  :=
  q\langle\xi,\psi(s)\rangle-s
  +q\bigl[\ln(2d-1)-\varrho^{*}(\xi)\bigr].
\]
For every $\xi\in\Omega^{1}$, the function
$s\mapsto\Phi_q(s,\xi)$ is continuous and convex, whereas, for every
$s\leq\ln R$, the function $\xi\mapsto\Phi_q(s,\xi)$ is continuous and
concave. Since $\Omega^{1}$ is compact and convex, Sion's minimax
theorem \cite[Corollary 3.3]{Sion1958} allows us to interchange the maximum and the infimum in
\eqref{eq-max-inf-L*}. Consequently,
\begin{align}
  - L^{*}(q)  &= \inf_{-\infty<s \leq \ln R} \ \Bigl(  \max_{\xi \in \Omega^{1}}  \left\{   q \langle \xi, \psi(s) \rangle +   q[ \ln (2d-1) - \varrho^{*}(\xi) ] \right\} - s \Bigr) \\
  &=  \inf_{-\infty<s \leq \ln R} q \, \varrho(\psi(s)) - s  \label{eq-inf-max-L*}
\end{align}
where the second equality follows from \eqref{eq-rho-via-rho*}.
Comparing \eqref{eq-inf-max-L*} with \eqref{eq-L*-via-P}, we see that
$P$ and $s\mapsto\varrho(\psi(s))$ have the same Legendre--Fenchel
dual. The Fenchel--Moreau theorem therefore yields $P(s)= \varrho(\psi(s))= \sup_{\alpha>0} \alpha s - \alpha L^{*}(1/\alpha)$, $s\leq\ln R$.

Finally, fix $r\in(1,R]$. By \cite[Theorem~1]{HL00},
$\dim_{\mathrm{H}}\Lambda_r$ is the unique solution $h$ of
$
  \sum_{a\in\mathcal{A}}
  \frac{F_a(r)}{e^h+F_a(r)}
  =1.$
Since $e^{\psi_a(\ln r)}=F_a(r)$ and
$F_a(r)=F_{a^{-1}}(r)$, this is precisely the defining equation
\eqref{eq-def-of-varrho} for
$\varrho(\psi(\ln r))$. Hence
$
  P(\ln r)
  =
  \dim_{\mathrm{H}}\Lambda_r$. This completes the proof.
\end{proof}

 \subsection{Asymptotics  of a free energy}
 \label{sec-asymp-F}
  In this section, we examine the asymptotic behavior of the free energy for a specific system defined as follows. Let $\lambda=(\lambda_{a})_{a \in \mathcal{A}} \in \mathbb{R}^{2d}$ be a given vector.
 For each element $x \in \mathbb{F}$,   we define its energy by $
 H_{\lambda}(x) := - \sum_{a \in \mathcal{A}} \lambda_{a} \Xi_a(x)  $.
 The corresponding Gibbs-Boltzmann distribution with parameter $\beta>0$ is given by $ \mathcal{G}^{(\lambda)}_{n,\beta} (x) = \frac{ 1}{\mathsf{Z}_{n,\beta}(\lambda) } e^{- \beta H_{\lambda}(x)} \ind{ x \in \mathbb{F}_{n}}$, where  $\mathsf{Z}_{n,\beta}(\lambda)$ is the  partion function   defined as
 \begin{equation}
   \mathsf{Z}_{n,\beta}(\lambda) : = \sum_{x \in \mathbb{F}_{n}} e^{-\beta H_{\lambda}(x)}  = \sum_{x \in \mathbb{F}_{n}} \exp   \Bigl(  \beta \sum_{a}  \lambda_{a} \Xi_a(x)    \Bigr) .
 \end{equation}
The following lemma yields that the normalized free energy $\frac{1}{n}\ln  \mathsf{Z}_{n,\beta}(\lambda) $ converges as $n \to \infty$.

\begin{lemma}\label{lem-eigenvalue} Let the function $\varrho$ be defined as in~\eqref{eq-def-of-varrho}.
Fix $\beta > 0$. Given any vector $\lambda =(\lambda_a)_{a \in \mathcal{A}} \in \mathbb{R}^{2d}$, we have
\begin{equation}
 \lim_{n\to\infty}\frac{1}{n} \ln \mathsf{Z}_{n,\beta}(\lambda)  =   \varrho(\beta \lambda) .
\end{equation}
\end{lemma}

\begin{proof}
 Without loss of generality we may assume  $\beta=1$.
As in the proof of Proposition \ref{prop-LDP-|Zn|},  let  $(X_{n}, \mathbf{P})$ be a random variable uniformly distributed on $ \mathbb{F}_{n} $.  Then, we have
\begin{equation}
    \mathsf{Z}_{n,1}(\lambda)   =  \# \mathbb{F}_{n} \cdot \mathbf{E} \left[ e^{ \sum_{a} \lambda_{a} \Xi_a(X_n) } \right].
\end{equation}
We represent  $X_{n}$  as a concatenation of the first $n$ elements of a Markov chain as follows. Let $\{ (W_k)_{k\geq 1}, \mathbf{P}\}$ be the  Markov chain on $\mathcal{A}$ with transition probabilities $
	 p(a,  b) = \frac{1}{2d-1} \mathbf{1}_{ \{ b\neq a^{-1} \} } $ for   $a,  b \in\mathcal{A}$
and initial distribution $\mathbf{P} ( W_1=a)=\frac{1}{2d}, \forall a\in\mathcal{A}$. The $ (W_k)_{k\geq 1}$, is stationary and hence
 $X_n := W_1 W_2 \cdots W_n  \, \text { for } n \geq 1$ is uniformly distributed on $\mathbb{F}_{n}$.

We shall employ the first step analysis to derive an iterative equation for   $  \mathsf{Z}_{n,\beta}(\lambda) $. Let
\begin{equation}
	g_a(n) := \mathbf{E}  \bigl[   e^{  \sum_{b} \lambda_{b} \Xi_b(X_n)  } \mid W_1=a \bigr]   \quad \text{ for } \  a \in \mathcal{A} \ , \  n \ge 1 .
\end{equation}
In particular, $g_a(1)= e^{\lambda_{a}}$. Applying the  Markov property gives
\begin{equation}
\begin{aligned}
	g_a(n) &=   \frac{1}{2d-1}\sum_{c \neq a^{-1}}  \mathbf{E} \bigl[  e^{  \sum_{b} \lambda_{b} \Xi_b(X_n) } \,\bigm|\, W_1=a , W_{2}=c \bigr]   \\
	& =\frac{e^{\lambda_{a}}}{2d-1}\sum_{c \neq a^{-1}}  \mathbf{E} \bigl[   e^{  \sum_{b} \lambda_{b}  \Xi_b(X_{n-1})  } \,\bigm|\,  W_{1}=c \bigr] =\frac{e^{\lambda_{a}}}{2d-1} \sum_{c \neq a^{-1}} g_c(n-1).
\end{aligned}\end{equation}
Let $M=(M_{a,b})$ be the  $2d \times 2d$ matrix with entries given by  $
  M_{a,   b} := e^{\lambda_{a}} 1_{ \{ b \neq a^{-1} \} } $ for $ a,b \in \mathcal{A}$.
%  (This form is analogous to the backscattering matrix defined in  \cite[equation (16)]{HL00}.)
Let $g(n)$ and $e^{\lambda}$
be the column vector  $( g_{a}(n))_{a \in \mathcal{A}}$   and   $(  e^{\lambda_{a}})_{a \in \mathcal{A}}$ respectively.
 We then get
\begin{equation}
g(n) =\frac{1}{2d-1} M \cdot g(n-1) = \cdots = \frac{1}{ (2d-1)^{n-1} } \, M^{n-1} \cdot  e^{\lambda}
\end{equation}
Since $ \#\mathbb{F}_{n}=2d(2d-1)^{n-1}$, it follows that
\begin{equation}
	\mathsf{Z}_{n,1}(\lambda) = 2d(2d-1)^{n-1} \cdot \frac{1}{2d} \sum_{a\in\mathcal{A}} g_a(n) = \sum_{a , b \in\mathcal{A}}  M^{n-1}_{a,   b} \, e^{\lambda_{b}} ,
\end{equation}
where $ M^{n-1}_{a,   b}$ is the $(a,  b)$-th element of $M^{n-1}$.
Since $M$  is aperiodic and  irreducible,  the Perron-Frobenius Theorem (see e.g. \cite[Theorem 3.1.1]{DZ09}) yields
\begin{equation}
 \lim_{n\to\infty} \frac{1}{n}  \ln 	\mathsf{Z}_{n,1}(\lambda) = \ln \rho^{\mathrm{PF}}(\lambda)  \ , \
\end{equation}
  where $   \rho^{\mathrm{PF}}(\lambda)$ is   the Perron--Frobenius eigenvalue of $M$.
Applying Lemma \ref{lem-PF-eigen-equa} to the matrix $M$, together with the definition  \eqref{eq-def-of-varrho},
  we obtain  $ \ln \rho^{\mathrm{PF}}(\lambda) = \varrho(\lambda) $ as desired.
 \end{proof}

\section{Proof of the Upper Bound}
In this section, we establish upper bounds for the Hausdorff dimensions of the random fractals $\Lambda_r(\alpha,\beta)$ and $E_r(\alpha,\beta)$, as stated in Theorem \ref{thm-Hdim-Lambda-alpha-beta} and Theorem \ref{thm-dim-E-alpha-beta}, respectively. We begin with introducing some notation. Let
\[
 \mathfrak{I}  := \bigl\{  [\alpha,\beta]:
      0\leq\alpha\leq\beta\leq1
 \bigr\}.
\]
For each $J \in \mathfrak{I}$ and $\epsilon>0$ define $J^\epsilon=\{x \in [0,1]: \dist(x,J) \le \epsilon \}=[\alpha-\epsilon,\beta+\epsilon]\cap[0,1]$.

Recall \eqref{eq:forward-cone}.
For $n \ge m \ge 1$ and $x, g \in \mathbb{F}$  define
\begin{equation}
  \mathcal{N}_{n,x}
   :=\sum_{u\in\mathcal{T}_n}\ind{V(u)=x}, \
  \mathcal{N}_{n,m}(g)
   :=\sum_{x\in\mathbb{F}_m(g)}\mathcal{N}_{n,x}, \
  \mathcal{N}^{\mathbb{F}}_{n,m}(g)
   :=\sum_{x\in\mathbb{F}_m(g)}
     \ind{\mathcal{N}_{n,x}\geq1}.
 \label{eq:def-level-counts}
\end{equation}
Thus $\mathcal{N}_{n,m}(g)$ counts particles with multiplicity,
whereas $\mathcal{N}^{\mathbb{F}}_{n,m}(g)$ counts occupied sites.
When $g=e$, we omit $g$ from the notation.

\subsection{Upper bound for $\dim_{\mathrm{H}} E_{r}(\alpha, \beta)$}
\label{sec-up-E}

Recall that a ray $t \in \partial \mathcal{T}$ is a semi-infinite self-avoiding path
$t=(t_n)_{n \geq 0}$ in $\mathcal{T}$ starting from the root; see \S\ref{sec-free-group}. For  $J \in \mathfrak{I}$, define
\begin{equation}
  \hat{E}_r( J ) : =   \bigl\{ t \in \partial \mathcal{T} :   \text{the sequence }
\bigl( |V(t_n)|/n \bigr)_{n \geq 1}
\text{ has a subsequential limit in }  J     \bigr\}  .
\end{equation}
Equivalently, we have
\begin{equation}\label{eq-E-hat-alpha-beta-0}
\hat{E}_r( J ) =  \bigcap_{\epsilon > 0} \bigcap_{k \geq 1} \bigcup_{n \geq k}  \, \bigl\{t \in \partial \mathcal{T} \colon    | V(t_n)  |/n \in J^{\epsilon}    \bigr\}.
\end{equation}

\begin{lemma}
  \label{P:upboundEhat}
Almost surely,  simultaneously for every interval
$J \in \mathfrak{I}$, if
$J\cap\hyperref[def-I-r]{I(r)}\neq\emptyset$, then
\begin{equation}
  \dim_{\mathrm{H}}\hat{E}_r(J)
  \leq \ln r-\min_{q\in J }
  \hyperref[eq-L*-via-P]{L^*}(q).
\end{equation}
If $J\cap I(r)=\emptyset$, then
$\hat{E}_r(J)=\emptyset$.
\end{lemma}

\begin{proof}
Fix  an arbitrary interval $J \in \mathfrak{I}$.
Applying \eqref{eq-E-hat-alpha-beta-0}, we have for every $k\geq1$ and
$\epsilon>0$,
\[
\hat{E}_r(J)
  \subset
  \bigcup_{n\geq k} \
  \bigcup_{ u\in\mathcal{T}_n: |V(u)|/n\in J^\epsilon} \
  [u]_{\mathcal{T}} \ .
\]
where $[u]_{\mathcal{T}}$ is defined in \eqref{eq:shadow}.
Since $\operatorname{diam}_{\partial\mathcal{T}}([u]_{\mathcal{T}})
\leq e^{-|u|}=e^{-n}$, for every $s\geq0$,
\begin{equation}
  \mathcal{H}^s_{e^{-k}}\bigl(\hat{E}_r(J)\bigr)
  \leq \sum_{n\geq k}e^{-sn} \sum_{m/n \in  J^{\epsilon}} \mathcal{N}_{n,m}
  \label{eq-cover-Ehat-count}
\end{equation}
The many-to-one formula, and the large deviation principle
\eqref{eq-ldp-word-length} (or Proposition \ref{prop-LDP-|Zn|}) give
\begin{align}
  \limsup_{n\to\infty}\frac1n\ln  \Bigl( \sum_{m/n \in  J^{\epsilon}} \E [ \mathcal{N}_{n,m}] \Bigr)
  & =\ln r+\limsup_{n\to\infty}\frac1n
    \ln\P  \Bigl(  \frac{|Z_n|}{n}\in J^\epsilon  \Bigr)  \\
    & \leq \ln r-\inf_{q\in J^\epsilon}L^*(q).
  \label{eq-growth-Ehat-count}
\end{align}
We now distinguish between two cases.
\begin{itemize}
  \item Suppose  $J\cap I(r)\neq\emptyset$. For any $s>\ln r- \min_{q\in J}L^*(q)$, the continuity of
$L^*$ allows us to choose $\epsilon>0$ such that
$s>\ln r-\inf_{q\in J^\epsilon}L^*(q)$. It follows from
\eqref{eq-growth-Ehat-count} that
$\sum_{n\geq1}e^{-sn}\sum_{m/n \in  J^{\epsilon}}  \E[\mathcal{N}_{n,m}]<\infty$. Hence, by
\eqref{eq-cover-Ehat-count} and monotone convergence,
\[
  \E\bigl[  \mathcal{H}^s\bigl(\hat{E}_r(J)\bigr) \bigr]
  =\lim_{k\to\infty}\E \bigl[
  \mathcal{H}^s_{e^{-k}}\bigl(\hat{E}_r(J)\bigr)\bigr]
  \leq\lim_{k\to\infty}
  \sum_{n\geq k}e^{-sn}\sum_{m/n \in  J^{\epsilon}}  \E[\mathcal{N}_{n,m}]=0.
\]
This implies $\mathcal{H}^s(\hat{E}_r(J))=0$ almost surely. Taking a
countable sequence of $s \downarrow \ln r- \min_{q\in J}L^*(q) $
proves the desired upper
bound for this fixed interval.

\item Suppose  $J\cap I(r)=\emptyset$.  The continuity of gives an $\epsilon>0$ such that
$\inf_{q\in J^\epsilon}L^*(q)>\ln r$. By
\eqref{eq-growth-Ehat-count}, we have
 \[
  \sum_{n\geq1}\P  ( \sum_{m/n \in  J^{\epsilon}}   \mathcal{N}_{n,m}  \ge 1  )
  \leq\sum_{n\geq1} \sum_{m/n \in  J^{\epsilon}}  \E[\mathcal{N}_{n,m}]<\infty. \]
The Borel--Cantelli lemma and \eqref{eq-E-hat-alpha-beta-0} therefore
imply that $\hat{E}_r(J)=\emptyset$ almost surely.
\end{itemize}

The preceding conclusions hold simultaneously for all intervals with
rational endpoints.
For an arbitrary $J$, choose rational intervals
  $J_m \downarrow J$. Then $\hat{E}_r(J)\subset\hat{E}_r(J_m)$, and
continuity of $L^*$ gives
$\min_{q\in J_m}L^*(q)\to\min_{q\in J}L^*(q)$, proving the upper bound.
If $J\cap I(r)=\emptyset$, one may instead choose a rational interval
$ J' \supset J$ disjoint from $I(r)$; then
$\hat{E}_r(J)\subset\hat{E}_r(J')=\emptyset$.
\end{proof}

Observe that  $E_r(\alpha,  \beta) \subset  \hat{E}_r(\{\gamma\} )$ for every $\gamma \in [\alpha,  \beta]$. The following corollary follows immediately from Lemma \ref{P:upboundEhat}.

\begin{corollary}
Almost surely,
simultaneously for every  $[\alpha , \beta ] \subset  I(r)$ we have
\begin{equation}
  \dim_{\mathrm{H}} E_r(\alpha,  \beta) \leq \ln r - \max\{   L^{*}(\alpha)  , L^{*}(\beta)\}.
\end{equation}
\end{corollary}

\subsection{Upper bound for $\Lambda_{r}(\alpha, \beta)$}
\label{sec-up-Lambda}

Recall that a ray $\omega   \in \partial\mathbb{F}$ is a semi-infinitely self-avoiding path $(\omega_{n})_{n \geq 0}$ starting from the root in $\mathbb{F}$ (see \S\ref{sec-free-group}); and we write $V(t)= \omega$ for $t \in \partial\mathcal{T}$ if the traces of $(V(t_{n}))_{n \geq 0}$ and $(\omega_{n})_{n \geq 0}$  intersect infinitely many times.
Let $r \in (1,R]$. Define
\[
 \mathfrak{I}_r :=
 \bigl\{
       [\alpha,\beta]:
      0\leq\alpha\leq\beta\leq1,\
     [\alpha,\beta]\subset I(r)
 \bigr\}.
\]
For each $J\in\mathfrak{I}$, define
\begin{equation}
  \underline{E}_r(J)
  :=
 \Bigl\{
    t\in\partial\mathcal{T}:
    \liminf_{n\to\infty}\frac{|V(t_n)|}{n}\in J
\Bigr\} ,
  \qquad
  \underline{\Lambda}_r(J)
  :=
  V \bigl(\underline E_r(J)\bigr).
\end{equation}
Equivalently, $ \underline\Lambda_r(J)
  =
 \{ \omega\in\partial\mathbb{F}:
    V(t)=\omega
    \text{ for some }t\in\underline E_r(J)\} $. Finally we define
    \[  h_r(q) = \frac{\ln r -L^{*}(q)}{q}  \ \text{ for } \  q \in (0,1] \quad \text{ and } \quad   h_R(0)
  =-(L^*)'(0) \, . \]
Consequently, for any $r \in (1,R]$,  $J \in \mathfrak{I}_r$ there exists $\epsilon>0$ such that  $h_{r}$ is continuous on $J^{\epsilon}$.

\begin{lemma}\label{lem-dim-hat-L-1}
 Let $r \in (1,R]$. Almost surely, simultaneously for every $J \in \mathfrak{I}_r$,
  \begin{equation}
    \label{eq:udim-J}
  \dim_{\mathrm{H}} \underline{\Lambda}_{r}(J)  \leq \max_{ q \in J} \, h_r(q)  \, .
  \end{equation}
  \end{lemma}

\begin{proof}
As in th proof of Lemma \ref{P:upboundEhat},  it suffices to show \eqref{eq:udim-J} for an arbitrary fixed  $J \in \mathfrak{I}_r$.

We assert that for any $ \omega = V(t)$ with $t \in \underline{E}_{r}(J)$,  there exists a sequence $\ell(m)$ such that $V(t_{\ell(m)})= \omega_{m}$, and   $\liminf\limits_{m \to \infty} \frac{\left| V(t_{\ell(m)}) \right|}{\ell(m)} $  lies in $ J$.
Consequently, for every $\epsilon>0$ and $k \ge 1$
\begin{align}
  \underline\Lambda_r(J)
 \subset
  \bigcup_{m\geq k} \
  \bigcup_{n :\: m/n \in J^{\epsilon}} \
  \bigcup_{ x\in\mathbb{F}_m:\:
  \exists u\in\mathcal{T}_n,  V(u)=x} \ [x]_{\mathbb{F}} \ .
  \label{eq-cover-L-underline}
\end{align}
To prove this, fix an end $\omega=(\omega_m)_{m \geq 1} \in \partial \mathbb{F}$ and a path $\gamma \in \mathsf{Path}_{\infty}$ in the equivalent class of $\omega$. Let $\ell(m)= \sup\{ \ell \geq 1: \gamma_{\ell}= \omega_{m}\}$ for $m \geq 1$.
 Since $|\gamma_{n}| \to \infty$ we have $\ell(m) < \infty$,  and   $ \gamma_{\ell(m)}= \omega_{m}$. By the geometry of the tree, $\ell(m)$ is increasing in $m$, and $|\gamma_n| \ge m$ whenever $\ell(m)\leq n\leq\ell(m+1)$.
 This  gives  $\liminf\limits_{m \to \infty} |\gamma_{\ell(m)}| / \ell(m) =  \liminf\limits_{n \to \infty} |\gamma_{n}|/n$. Taking $ \gamma = t \in \underline{E}_{r}(J) $ proves the preceding assertion.

 Each shadow $[x]_{\mathbb{F}}$ in \eqref{eq-cover-L-underline}   has diameter less than
$e^{-m} $,
and the number of occupied sites $x\in\mathbb{F}_m$ is at most the
number of particles in $\mathbb{F}_m$.
Hence, for any $s\geq0$, $k\geq1$, and
$\epsilon>0$,
\[\mathcal{H}^s_{e^{-k}} \bigl(   \underline{\Lambda}_r (J)  \bigr)
\leq   \sum_{ m\geq k} e^{-sm}    \sum_{m/n \in J^{\epsilon}}  \mathcal{N}_{n,m}. \]

 Fix $s>\max_{q\in J}h_r(q)$.  By continuity of $h_{r}$, one can choose
$\epsilon,\eta>0$ such that
\[
 \max_{q\in J^\epsilon } h_r(q)\leq s-\eta \, .
\]
Applying the many-to-one formula and Proposition \ref{prop-LDP-|Zn|} gives
\begin{align}
  \E \bigl[  \mathcal{H}^s_{e^{-k}}
  \bigl(\underline\Lambda_r(J)\bigr)  \bigr]
  &\leq
  \sum_{m\geq k} \
  \sum_{n:\: m\in\mathtt{L}_n^\mu, m/n\in J^\epsilon} \
  e^{-sm} \, r^n\, \P(|Z_n|=m) \\
  &=
   \sum_{m\geq k} \
  \sum_{n:\: m\in\mathtt{L}_n^\mu, m/n\in J^\epsilon} \
  B_n(m)\exp \bigl\{  m[h_r(m/n)-s]   \bigr\} \\
  &\leq
  \sum_{m\geq k}e^{-\eta m}S_m \ , \qquad   \text{ where } \quad S_m:= \!\!\!
  \sum_{n:\: m\in\mathtt{L}_n^\mu, m/n\in J^\epsilon} B_n(m).
  \label{eq-hat-L-a-b-0}
\end{align}

 We next show that
$S_m=e^{o(m)}$.  Split its defining
sum at $e^{\sqrt m}$.  Proposition~\ref{prop-LDP-|Zn|} gives
$B_n(m)\leq B(m)\sup_{\xi\in\Omega^1}\beta_n((m/n)\xi)$ and
$\ln B(m)=O(\ln m)$.  For $n>e^{\sqrt m}$, we have
$m\leq(\ln n)^2$.  Using, respectively, the global and small-argument
bounds for $\beta_n(\cdot)$ in Proposition~\ref{T:LLT}, we obtain
\begin{align*}
  S_m
  &\lesssim B(m)e^{\sqrt m}
  +B(m)\sum_{n>e^{\sqrt m}}
  \frac{1+(\ln n)^2}{n^{3/2}}
  \leq e^{o(m)}.
\end{align*}
It then follows from \eqref{eq-hat-L-a-b-0} that $\mathbf{E} [  \mathcal{H}^s_{e^{-k}} ( \underline{\Lambda}_r (J)   ) ]$ tends to zero as $k\to\infty$.  Monotone convergence therefore
gives $\mathcal{H}^s(\underline\Lambda_r(J))=0$ almost surely.  Letting
$s\downarrow\max_{q\in J}h_r(q)$ proves the required result.
\end{proof}

Note that  $\Lambda_r([\alpha,  \beta]) \subset  \underline{\Lambda}_r([\alpha, \gamma] )$ for every $\gamma \geq \alpha$; and $\Lambda_{r} =\underline{\Lambda}_r([0,1]) = \Lambda_{r}( I(r) ) $ almost surely.
As an immediate consequence of Lemma \ref{lem-dim-hat-L-1}, we obtain the following corollary.

\begin{corollary}
	Let $r \in (1,R]$. Almost surely, simultaneously for every $[\alpha,\beta] \subset I(r)$
  \begin{equation}
    \dim_{\mathrm{H}}  \Lambda_{r} (\alpha,  \beta) \leq   h_r(\alpha)  \ ; \  \text{ and  moreover } \ \dim_{\mathrm{H}}  \Lambda_{r}  \leq   \max_{q \in I(r)} h_r(q) \ .
  \end{equation}
\end{corollary}

\begin{remark} \label{rmk-difficulty-1}
The gap in the preceding corollary  arises from a geometric obstruction to
proving a spatial analogue of Lemma~\ref{P:upboundEhat}.  For
$J\in\mathfrak{I}$, set
$ \hat\Lambda_r(J)
  := V(\hat E_r(J) )$. Such an analogue would assert that
\[
  \dim_{\mathrm{H}}\hat\Lambda_r(J)
  \leq \max_{q\in J}h_r(q)\ ,
  \quad J \subset I(r).
\]
Since
$\Lambda_r(\alpha,\beta)\subset\hat\Lambda_r(\{\gamma\})$
for every $\gamma\in[\alpha,\beta]$, applying this   to
$J=\{\gamma\}$ and minimizing over $\gamma$ would yield the  
upper bound $\dim_{\mathrm{H}}\Lambda_r(\alpha,\beta)
  \leq\min_{q\in[\alpha,\beta]}
 h_r(q).$

The obstruction lies in the covering step.  If
$t\in\hat E_r(J)$ and $|V(t_{n_k})|/n_k\to q\in J$, the vertices
$V(t_{n_k})$ need not lie on the limiting ray $\omega=V(t)$; equivalently,
$\omega$ need not belong to $[V(t_{n_k})]_{\mathbb{F}}$, since a finite
excursion may subsequently be erased by backtracking.  This explains
precisely why the shadows used in
\eqref{eq-cover-L-underline} need not cover $\hat\Lambda_r(J)$.

A valid cover should instead use a prefix $z(u)$ of $V(u)$ such that every
descendant $v\in\mathcal{T}(u)$ remains in the forward cone
$\mathbb{F}(z(u))$; then every limiting end generated by a genealogical ray
through $u$ belongs to $[z(u)]_{\mathbb{F}}$.  The vertex $z(u)$ introduced
below has precisely this property, but its shadow has diameter
$e^{-|z(u)|}$ rather than
$e^{-|V(u)|}$.  Controlling the resulting larger shadows is exactly where
the additional estimate in \eqref{eq-pressure-comparison} is needed.
\end{remark}

 Recall that  $\psi(s):= (\psi_{a}(s))_{a \in \mathcal{A}}$,  $\varrho(\cdot)$   and  $ P(s) :=  \varrho( \psi(s) )$ are defined in \eqref{eq-def-psi-a-s}, \eqref{eq-def-of-varrho} and \eqref{eq-forluma-P} respectively.  Let
\begin{equation}
      \hat{P}(s)  :=  \varrho( \psi(s) + \psi(\ln R)) \ \text{ for }\ s  \in (-\infty,\ln R].
\end{equation}
Note that  $\hat{P}'(s) = \langle \nabla \varrho(\psi(s)+\psi(\ln R)), \psi'(s) \rangle$. Since $ \nabla \varrho(\lambda) \in \Omega^{1}$ for all $\lambda \in \mathbb{R}^{2d}$ satisfying $\lambda_{a}=\lambda_{a^{-1}} , \forall \, a \in \mathcal{A}$ (in fact for such $\lambda$, $\frac{\partial}{\partial \lambda_{b}} \varrho(\lambda) =  \frac{ e^{\lambda_{b}}}{ (e^{\varrho(\lambda) } +   e^{\lambda_{b}})^2 } /  \sum_{a \in \mathcal{A}} \frac{ e^{\lambda_{a}}}{ (e^{\varrho(\lambda) } +   e^{\lambda_{a}} )^2} $), by applying \eqref{eq-properties-psi}, we have $\hat{P}'(s) \geq 1 $ for $s \in (-\infty,\ln R)$,  $ \lim _{s \uparrow \ln R} \hat{P}'(s) =\infty $ and  $ \lim _{s \downarrow -\infty} \hat{P}'(s) =1$. Additionally, $ \hat{P}''(s)>0$ for all $s \in (-\infty, \ln R)$. To verify this, let $H\varrho$ denote the Hessian matrix of $\varrho$. We then compute:
\begin{align}
  \hat{P}''(s) &=  \langle (H\varrho)(\psi(s)+\psi(\ln R))\psi'(s), \psi'(s) \rangle +  \langle \nabla \varrho(\psi(s)+\psi(\ln R)), \psi''(s) \rangle \\
  &  \geq \langle \nabla \varrho(\psi(s)+\psi(\ln R)), \psi''(s) \rangle \geq \min_{a \in \mathcal{A}} \psi''_{a}(s) > 0.
\end{align}
Here we used the facts that  $H\varrho$ is non-negative definite as $\varrho$ is convex,  $ \nabla \varrho(\lambda) \in \Omega^{1}$  and the property \eqref{eq-properties-psi} of $\psi$.  Moreover from Lemma \ref{lem-rho2psi} it follows
 \begin{equation}\label{eq-tilde-P-lnR-is-0}
  \hat{P}(\ln R)= \varrho(2 \psi(\ln R)) =  0 .
\end{equation}

In the following lemma,  we require the following  pressure comparison condition on the step distribution $\mu$:
 \begin{equation}\label{eq-pressure-comparison}
  \frac{\hat{P}(s)}{P'(s)} - s + \ln R \leq  0 \quad \text{ for all } \quad s \in (-\infty,\ln R). \tag{PC}
\end{equation}
This condition is automatically satisfied when the random walk is isotropic,   in which case all coordinates of $\psi(\cdot)$ coincide. Together with \eqref{eq-forluma-P} and \eqref{eq-def-of-varrho}, this implies that the difference between $\hat{P}(s)$ and $P(s)$ is  independent of $s$, and therefore $\hat{P}'(s) = P'(s)$.
Since $\hat{P}(s)$ is convex and $\hat{P}(\ln R)=0$,  \eqref{eq-pressure-comparison} follows immediately.

\begin{lemma}\label{lem-assume-Hypothesis}
  Let $r \in (1,R]$. Assume  \eqref{eq-pressure-comparison}. Almost surely, simultaneously for every $J \in \mathfrak{I}_r$,
\begin{equation}
  \dim_{\mathrm{H}} \hat{\Lambda}_{r}(J)  \leq \max_{ q \in J} h_r(q).
\end{equation}
Consequently, almost surely for every $[\alpha,\beta] \subset I(r)$, $
  \dim_{\mathrm{H}}  \Lambda_{r}(\alpha, \beta)  \leq \min_{ q \in [\alpha,  \beta]} h_r(q)$.
\end{lemma}

\begin{proof} %% We can consider those paths with positive liminf speed.
%----------------------------------
We first dispose of the case $J=[0,\beta]$.  Then necessarily $r=R$, and
Lemma~\ref{lem-dim-hat-L-1}, along with the  the fact $ \hat{\Lambda}_{R}([0, \beta])  \subset  \Lambda_{R}([0, 1])$  gives   the desired upper bound:
 \begin{equation}
  \dim_{\mathrm{H}} \hat{\Lambda}_{r}([0, \beta])  \leq  \max_{ q \in [0, 1]}  h_R(q) = h_R(0) =\max_{ q \in [0, \beta]}    h_R(q) .
 \end{equation}
 Henceforth we assume $J \in \mathfrak{I}_r$ with $\min J>0$.
Fix
$s>\max_{q\in J}h_r(q)$, and we may choose $\epsilon,\eta>0$ so that $J^{\epsilon} \subset (0,1]$ and
$ s\geq\max_{q\in J^{\epsilon}}h_r(q)+2\eta$.

\smallskip \noindent
 \underline{\textit{Step 1}}.
 For each $u \in \mathcal{T}$,  let $z(u)$ be the vertex on $[e,V(u)]$ with
$  |z(u)|=\min_{v\in\mathcal{T}(u)}|V(v)|$.  Then
every descendant of $u$ remains in the forward cone  $\mathbb{F} (z(u))$.
For each $t\in\hat{E}_r(\alpha,\beta)$ with $V(t)=\omega$,  for every $l\geq1$, there must exist $m\geq l$ and $n$
such that
$|V(t_n)|=m$ and $  m/n \in J^{\epsilon}$.   Since
$V(t_j)\in\mathbb{F}(z(t_n))$ for all $j\geq n$, we have
$\omega_{|z(t_n)|}=z(t_n)$. Consequently,
\begin{equation}
  \hat\Lambda_r(J)
  \subset \bigcap_{l\geq1} \ \bigcup_{m\geq l}
 \  \bigcup_{  n:\: m/n\in J^{\epsilon}} \
 \bigcup_{x \in \mathbb{F}_m }
  \ \bigcup_{ u\in\mathcal{T}_n :  V(u)=x }[z(u)]_{\mathbb{F}}.
\end{equation}

The shadow $[z(u)]_{\mathbb{F}}$ has diameter $e^{-|z(u)|}$.
Set $\delta_l^V:=\max\{e^{-|z(u)|}:u\in\mathcal{T}_n,
n\geq l/(\beta+\epsilon) \}$.
Transience of the BRW implies that $\delta_l^V\downarrow0$ almost surely. The branching and many-to-one properties yield
\begin{equation}
\E \bigl[  \mathcal{H}^s_{\delta^{V}_{l} }   \bigl(  \hat{\Lambda}_r(J)\bigr) \bigr]
 \leq   \sum_{m\geq l} \
 \sum_{    n:\: m/n\in J^{\epsilon}} \ r^{n} \sum_{x \in \mathbb{F}_m }  \E \bigl[  e^{ -s |z(u)|} \mathbf{1}_{ \{  V(u)= x \} }  \bigr]
\end{equation}

\smallskip \noindent
 \underline{\textit{Step 2}}.
Define
$ \kappa_s(x):=\E_{x} [\exp\{-s\min_{v\in\mathcal{T}}|V(v)|\}]$ for each $x \in \mathbb{F} $.
The branching property yields
$ \E  [  e^{ -s |z(u)|} \mathbf{1}_{ \{  V(u)= x \} }  ]= \P(Z_n= x) \kappa_{s}(x) $.

Fix $x\in\mathbb{F}_m$. To bound $\kappa_s(x)$, note that for $1\leq k\leq m$,
if the minimum in the definition of $\kappa_s(x)$
equals $m-k$, then the descendant BRW must hit $z^{-1}$; where we  write uniquely $x=yz$
with $|y|=m-k$ and $|z|=k$.
By symmetry, the
union bound, and the many-to-one formula,
 \begin{align}
 \P \bigl(  \min_{v\in\mathcal{T}}|V(v)| = m-k \bigr)
 &\leq \P(\text{the BRW hits }z^{-1})
  =\P(\text{the BRW hits }z) \\
 &\leq F_z(r)
  =\prod_{a\in\mathcal{A}}F_a(r)^{\Xi_a(z)}
  =\exp\{\langle\psi(\ln r),\Xi(z)\rangle\}.
\end{align}
Here recall $F_z$ is defined in \eqref{eq:generating-fun-F}  and we used  the branching property.

We thus set,
for $0\leq k\leq m$,
\begin{equation}
\Sigma_{\refeq{eq-H-s-delta-4}}(s;m,n,k):=r^n e^{-s(m-k)}
\sum_{\substack{y\in\mathbb{F}_{m-k},\ z\in\mathbb{F}_k\\yz\in\mathbb{F}_m}}
\P (Z_n=yz)
e^{\langle\psi(\ln r),\Xi(z)\rangle}.
\label{eq-H-s-delta-4}
\end{equation}
  Combining the previous estimates, we obtain that, for any   $ l \geq 0$,
\begin{equation}
\E \bigl[  \mathcal{H}^s_{\delta^{V}_{l} }   \bigl(  \hat{\Lambda}_r(J)\bigr) \bigr]
  \le  \sum_{m\geq l} \
 \sum_{    n:\: m/n\in J^{\epsilon}}  \   \sum_{k=0}^{m}   \, \Sigma_{\refeq{eq-H-s-delta-4}}(s;m,n,k).
\end{equation}

\smallskip \noindent
\underline{\textit{Step 3}}.  It suffices to show the existence of a constant $c  >0$ (depending on  $\epsilon, \eta, J $) such that
\begin{equation}\label{eq-claim-upp-Lambda-10}
   \sum_{k=0}^{m}  \Sigma_{\refeq{eq-H-s-delta-4}}(s;m,n,k)
   \leq e^{-cn}
   \quad\text{whenever }m/n\in J^{\epsilon}\text{ and }n\text{ is sufficiently large}.
\end{equation}
Together with the fact that $\lim_{l \to \infty} \delta^{V}_{l} =0$ almost surely, and the dominated convergence theorem,  we then get  $\E  [  \mathcal{H}^s_{\delta^{V}_{l} }   (  \hat{\Lambda}_r(J) )  ]	  \leq \sum_{m \ge l} O(m) e^{- c' m } \xrightarrow{l \to \infty} 0$. Thus
$\dim_{\mathrm{H}}\hat\Lambda_r(\alpha,\beta)\leq s$ almost surely, and
letting $s\downarrow\max_{q\in[\alpha,\beta]}h_r(q)$ proves the desired result.

To prove \eqref{eq-claim-upp-Lambda-10},
 applying Proposition \ref{T:LLT} gives
 \begin{align}
  & \Sigma_{\refeq{eq-H-s-delta-4}}(s;m,n,k)  \lesssim   r^{n}  e^{- s(m-k)} \sum_{\substack{  y \in \mathbb{F}_{m-k}, z \in \mathbb{F}_{k} \\ yz \in \mathbb{F}_{m}}}   e^{  n \hyperref[eq-def-Psi*]{\Psi^{*}}\left( \frac{(m-k)\xi(y) +k\xi(z)}{n} \right)  + k \langle \psi(\ln r) ,\xi(z) \rangle } \\
  & \lesssim r^{n} e^{- s(m-k)} \sum_{\xi \in \Omega^{1}_{m-k},    \eta \in \Omega^{1}_{k}}   e^{  n \Psi^{*}\left( \frac{(m-k)\xi  +k\eta }{n} \right) + k \langle \psi(\ln r) ,\eta\rangle } \sum_{y \in \mathbb{F}_{m-k}}   \ind{ \xi(y)= \xi }  \sum_{z \in \mathbb{F}_{k} } \ind{ \xi(z)= \eta} .
 \end{align}   For simplicity, set $q:= m/n$ and $\theta := k/m$.
From Lemma \ref{lem-ldp-empirical-measure} we have $\sum_{x \in \mathbb{F}_{n}}   \ind{ \xi(x)= \xi }  = \frac{2d}{2d-1} B_{n}(\xi)  (2d-1)^{n} e^{-n\varrho^{*}(\xi) } $ with $  \sup_{ \xi \in  \Omega^{1}_n } |\ln B_{n}( \xi)| =O(\ln n ) $. Consequently, uniform in $0\leq k\leq m$,
\begin{equation}\label{eq:Sigma-F-bnd-1}
  \Sigma_{\refeq{eq-H-s-delta-4}}(s;m,n,k) \lesssim   \exp \bigl\{
    n \bigl(   \mathsf{F}_{r}(q,\theta) +  (1-\theta)q [  h_r(q)- s] \bigr) + O(\ln m)  \bigr\} ,
\end{equation}
where
\begin{align}
\mathsf{F}_{r}(q,\theta)  : =\max_{\xi,\eta \in \Omega^{1}} &  \big\{ \Psi^{*} ( q [(1-\theta) \xi + \theta \eta ] )   +  \theta q  [\langle \psi(\ln r), \eta \rangle +  \log (2d-1) - \varrho^{*}(\eta)]
 \\
  & \qquad  +  \theta \ln r  + (1-\theta) L^{*}(q)   +  (1-\theta)q [\log (2d-1) - \varrho^{*}(\xi)] \big\}   \label{eq-bound-sigma-smnk-2}
\end{align}
Here, the endpoints  $k \in \{0,m\}$ indeed should be  handled separately: the case $k=0$ follows from  Proposition
\ref{prop-LDP-|Zn|}, \eqref{eq-forluma-Lstar}, and $\mathsf{F}_r(q,0) \equiv 0$ while the case $k=m$ follows from Proposition~\ref{T:LLT}.

 We need to obtain a uniform negative exponent in \eqref{eq:Sigma-F-bnd-1}.
We claim there is $c>0$ such that
\begin{equation}\label{eq:negative-exp}
\mathsf{D}(q,\theta):=\mathsf{F}_R(q,\theta)-2\eta q(1-\theta) \le - c \quad  \text{ on }  J^\epsilon\times[0,1]
\end{equation}
Since each
$\psi_a$ is increasing, $\mathsf{F}_r(q,\theta)\leq
\mathsf{F}_R(q,\theta)$.
Thus as  $s\geq\max_{q\in J^\epsilon}h_r(q)+2\eta$,
\eqref{eq:Sigma-F-bnd-1} now gives
$
\Sigma_{\refeq{eq-H-s-delta-4}}(s;m,n,k)
\leq e^{-c_0n+O(\ln m)}$ uniformly in $0\leq k\leq m$,.
Because $m\leq n$, summing over $0\leq k\leq m$ and decreasing
$c$ if necessary proves \eqref{eq-claim-upp-Lambda-10}.

\smallskip \noindent
\underline{\textit{Step 4}}.
It remains to show \eqref{eq:negative-exp}. Since $\mathsf{D}$ is continuous on the compact set $J^\epsilon\times[0,1]$, it suffices to show
 \begin{equation}
   \mathsf{F}_R(q, \theta) \le 0 \quad \text{ and } \quad  \mathsf{F}_R(q, 1) < 0 \quad  \text{ for every } \quad q \in J^{\epsilon}, \theta \in [0,1].  \label{eq:negative-F-value}
 \end{equation}
Recalling the definition of $\Psi^{*}$ in \eqref{eq-def-Psi*}, we rewrite the function $\mathsf{F}_{R}(q,\theta)$ as
\begin{align}
   \max_{\xi,\eta \in \Omega^{1}} \inf_{s \leq \ln R} & \big\{
    q \theta [ \langle \psi(s)+\psi(\ln R), \eta  \rangle  +\ln (2d-1) - \varrho^{*}(\eta) ] + \theta \ln R \\
   & \quad   + q(1-\theta) [ \langle \psi(s), \xi \rangle + \ln (2d-1) - \varrho^{*}(\xi) ] + (1-\theta) L^{*}(q) - s \big\}.
\end{align}
Notice that the function in parentheses is convex in $s$ and concave in $(\xi,\eta)$.
Since $\Omega^{1}\times \Omega^{1}$ is compact, using Sion's minimax
theorem \cite[Corollary 3.3]{Sion1958} again,
we get
\begin{align}
 \mathsf{F}_{R}(q,\theta) &=  \inf_{s \leq \ln R} \Big( q \theta \max_{\eta \in \Omega^{1}} \big\{  \langle \psi(s)+\psi(\ln R), \eta  \rangle + \ln(2d-1) - \varrho^{*} (\eta)  \big\} + \theta \ln R \\
  & \qquad + q (1-\theta) \max_{\xi \in \Omega^{1}} \big\{ \langle \psi(s), \xi \rangle + \ln(2d-1) - \varrho^{*}(\xi)  \big\} + (1-\theta) L^{*}(q) - s \Big) \\
  &= \inf_{s \leq \ln R} \Big( \theta [ q \hat{P}(s) -s + \ln R ] + (1-\theta) [ q P(s) - s +L^{*}(q) ] \Big) .  \label{eq-bound-sigma-smnk-3}
\end{align}

For fixed $q\in(0,1)$, the function $\theta \mapsto F_R(q,\theta)$ is
concave in $\theta$. Indeed, this follows from the inequality
$\inf (f_{1} + f_{2}) \geq \inf f_{1} + \inf f_{2}$ for arbitrary
functions $f_{1},f_{2}$. Moreover, by
\eqref{eq-inf-max-L*} and \eqref{eq-forluma-P}, we have
$\mathsf{F}_R(q,0) \equiv 0$.
For $q\in(0,1)$, let $s(q)$ be the unique solution of
$qP'(s(q))=1$. Such a solution exists since $P'(s)$ is strictly
increasing on $(-\infty,\ln R)$, with
$\lim _{s \uparrow \ln R} P'(s) =\infty$ and
$\lim _{s \downarrow -\infty} P'(s) =1$.
We claim that
\begin{equation}\label{eq:d-theta-F-neg}
    \frac{\partial}{\partial\theta} \Big|_{\theta=0+}
  \mathsf{F}_R(q,\theta)
   =\frac{\hat{P}(s(q))}{P'(s(q))}-s(q)+\ln R\leq0.
\end{equation}
Using $\mathsf{F}_R(q,0)=0$ and  concavity,
\eqref{eq:d-theta-F-neg} implies
$\mathsf{F}_R(q,\theta)\leq0$ for $q\in(0,1)$ and
$\theta\in[0,1]$. The same conclusion at $q=1$ follows from the
continuity of the maximization formula
\eqref{eq-bound-sigma-smnk-2}.

To show the equality in \eqref{eq:d-theta-F-neg}, note that
$\hat{P}'$ has the same monotonicity and endpoint behavior as $P'$.
Hence, for each $q \in (0,1)$, there is a unique solution
$\hat{s}= \hat{s}(q,\theta)$ of
$ \theta q \hat{P}'(s) + (1- \theta) q P'(s) = 1 $.
Consequently,
$
   \mathsf{F}_{R}(q,\theta) = \theta ( q \hat{P}(\hat{s}) -\hat{s} + \ln R ) + (1-\theta)  [ q P(\hat{s}) - \hat{s} +L^{*}(q) ]. $
Differentiating with respect to $\theta$ and using the defining
equation for $\hat{s}$ gives
\[  \frac{\partial }{\partial \theta}  \mathsf{F}_R(q,\theta) = q \hat{P}(\hat{s}) - \hat{s} + \ln R - [ q P(\hat{s})-\hat{s}+L^{*}(q) ].   \]
Letting $\theta \to 0+$ and using
$\hat{s}(q,\theta) \to s(q)$ and   $L^{*}(q)=qP(s(q))-s(q)$ proves the claimed equality.    The
inequality in \eqref{eq:d-theta-F-neg} then follows from
\eqref{eq-pressure-comparison}.

We finally show $\mathsf{F}_R(q, 1)<0$ for all $q \in {J}^{\epsilon}$.
Set $q_-:=\min J^\epsilon>0$. Since
$\hat{P}'(s)\to\infty$ as $s\uparrow\ln R$, one can choose
$s_0<\ln R$ such that
$q_-\hat{P}'(s)>1$ for $s\in[s_0,\ln R)$. As
$\hat{P}(\ln R)=0$, the function
$s\mapsto q_-\hat{P}(s)-s+\ln R$ is strictly increasing on
$[s_0,\ln R]$ and vanishes at $\ln R$. Hence there is a constant
$c_1>0$ such that, for every $q\in J^\epsilon$,
\begin{equation}
\mathsf{F}_R(q,1)
\leq q\hat{P}(s_0)-s_0+\ln R
\leq q_-\hat{P}(s_0)-s_0+\ln R=-c_1<0.
\label{eq-bound-sigma-smnk-strict-endpoint}
\end{equation}
This completes the proof.
\end{proof}

\section{Level Sets of BRWs on Free Groups}
 \label{sec-level-set}In this section, we establish key deviation estimates for the sizes of level sets of the BRW on $\mathbb{F}$. We continue to use the level-set counts introduced in
\eqref{eq:def-level-counts}.

% Fix $x \in \mathbb{F}$, $n \geq 1$ and $0 \leq m \leq n$, and let
% $\mathcal{N}_{n,x}$ denote the number of generation-$n$ particles sitting at $x$,
% i.e.
% \begin{equation}
%    \mathcal{N}_{n,x} := \sum_{u \in \mathcal{T}_{n}} 1_{\{V(u)=x\}}.
% \end{equation}

%  Recall that
% $\mathbf{P}_g$ denotes the law of the BRW started from $g$, and
% write $\prec$ for $\prec_{ \mathbb{F}}$ throughout this section.
% Let $\mathcal{N}_{n,m}(g)$ denote the number of generation-$n$
% particles whose locations lie in the forward cone rooted at $g$ and
% are at distance $m$ from $g$; and
% $\mathcal{N}^{ \mathbb{F}}_{n,m}(g)$ denote the number of distinct
% locations occupied by these particles:
% \begin{equation}
%   \mathcal{N}_{n,m}(g) := \sum_{x\in\mathbb{F}_m(g)} \mathcal{N}_{n,x} \quad ,
%   \quad
%   \mathcal{N}_{n,m}^{\mathbb{F}}(g) := \sum_{x\in\mathbb{F}_m(g)}
%      1_{\{\mathcal{N}_{n,x} \geq 1\}}.
% \end{equation}
% We write $\mathcal{N}_{n,m}$ and $\mathcal{N}_{n,m}^{\mathbb{F}}$ for $\mathcal{N}_{n,m}(e)$ and $\mathcal{N}_{n,m}^{\mathbb{F}}(e)$.

\begin{lemma}\label{lem-level-set-size}
 Let $r \in (1,\infty)$.  There is a decreasing sequence $(\hyperref[lem-level-set-size]{\delta}_{n})_{n \geq 1}$ defined in \eqref{def-of-deltan}, satisfying $\hyperref[lem-level-set-size]{\delta}_{n} \ll 1 \ll n\hyperref[lem-level-set-size]{\delta}_{n}$, such that, for any fixed $\epsilon >0$, for all sufficiently large $n$ and all $m\in\mathtt{L}^{\mu}_{n}$ satisfying $m \leq n(1-\hyperref[lem-level-set-size]{\delta}_{n})$ and $\ln r-\hyperref[eq-L*-via-P]{L^{*}}(m/n) \geq \hyperref[lem-level-set-size]{\delta}_{n}$, one has
  \begin{align}
     & \mathbf{P} \bigl(  \mathcal{N}_{n,m}  \leq [ r^{n} e^{-n L^{*}(\frac{m}{n})}  ]^{1- \epsilon} \bigr)
   \leq   e^{-\sqrt{n}} \ ,  \quad \text{ and }   \label{eq:level-set-size}\\
   & \mathbf{P}_{g} \bigl(  \mathcal{N}_{n,m}(g) \geq [ r^{n} e^{-n L^{*}(\frac{m}{n})}  ]^{1- \epsilon}  \bigr) \ge \frac{1}{3}  \quad \text{for all} \quad g \in \mathbb{F}. \label{eq:forward-cone-size}
\end{align}
   \end{lemma}

 \begin{lemma}\label{lem-level-set-size-F}
Let  $r \in(1,R]$. There exists
 a decreasing sequence $(\hyperref[lem-level-set-size-F]{\tilde{\delta}}_{n} )_{n \geq 1}$ defined in \eqref{def-of-tilde-deltan}, satisfying $\hyperref[lem-level-set-size-F]{\tilde{\delta}}_{n} \ll 1 \ll n \hyperref[lem-level-set-size-F]{\tilde{\delta}}_{n}$, such that for all sufficiently large $n$ and all $m\in\mathtt{L}^{\mu}_{n}$ satisfying $m\leq n(1-\hyperref[lem-level-set-size-F]{\tilde{\delta}}_{n})$ and $\ln r-L^{*}(m/n)\geq\hyperref[lem-level-set-size-F]{\tilde{\delta}}_{n}$, one has
  \begin{align}
   & \mathbf{P}   \bigl(  \mathcal{N}_{n,m}^{\mathbb{F}}  \leq [ r^{n} e^{-n L^{*}(\frac{m}{n})}  ]^{1- \epsilon}  \bigr)
   \leq \, e^{-\sqrt{n}} \ ,  \quad \text{ and }      \label{eq:level-set-size-F}\\
  & \mathbf{P}_g \bigl(  \mathcal{N}_{n,m}^{\mathbb{F}}(g) \geq [ r^{n} e^{-n L^{*}(\frac{m}{n})}  ]^{1- \epsilon} \bigr)   \ge \frac{1}{3}    \quad \text{for all} \quad g \in \mathbb{F}.\label{eq:forward-cone-size-F}
  \end{align}
  \end{lemma}

    \begin{remark}
  Lemmas~\ref{lem-level-set-size} and \ref{lem-level-set-size-F} imply Law of Large Numbers theorems for the size of level sets. By applying Markov's inequality to control the upper deviation probabilities and using the Borel-Cantelli Lemma, we conclude that  for $r \in (1,R]$,  almost surely
  \begin{equation}
\lim_{n \to \infty} \frac{1}{n} \ln \mathcal{N}_{n, \langle qn\rangle_{n,\mu}}=\lim_{n \to \infty} \frac{1}{n} \ln \mathcal{N}_{n, \langle qn\rangle_{n,\mu}}^{\mathbb{F}}=\ln r-L^{*}(q), \quad \forall \, q \text{ such that } L^{*}(q)<\ln r .
  \end{equation}
  Here we used the notation \eqref{eq:solve-periodicity-2} to make sure $\mathcal{N}_{n, \langle qn\rangle_{n,\mu}}$ and $\mathcal{N}_{n, \langle qn\rangle_{n,\mu}}^{\mathbb{F}}$ are nonempty.
  \end{remark}

  \begin{remark}
 In general, the probability bounds \eqref{eq:forward-cone-size} and \eqref{eq:forward-cone-size-F} in Lemmas~\ref{lem-level-set-size} and \ref{lem-level-set-size-F} cannot be tightened to hold with probability approaching 1. Indeed, suppose that $1<r\leq R$, fix a generator
$a\in\mathcal{A}$, and start the BRW from $g=a$.  Choose $k\geq1$
with $p_k>0$.  With probability at least  \[  p_k\, \bigl[ \mu(a^{-1})  (1-F_a(r) ) \bigr] ^{k} >0, \]
the root has $k$ children, all of them move from $a$ to $e$, and
none of the $k$ descendant BRWs ever come back to $a$.  Here the last
factor is positive because the union bound and the many-to-one formula
give  $  \P_e(\text{the BRW ever visits }a)
   \leq F_a(r)\leq F_a(R)<1$.
 \end{remark}

  \subsection{Sample paths large deviations}
  \label{sec-pathLD}
In this subsection, we study a specific sample path large deviation property for the word length of the random walk $(Z_{n})_{n \geq 0}$ on $\mathbb{F}$.

To illustrate,   consider a random walk $(S_n)$ on $\mathbb{R}^{d}$ with i.i.d. Gaussian increments $\mathcal{N}(0,\Sigma)$. Since  $S_{k}- kS_{n}/n$ is independent to $S_{n}$ for each $k \leq n$,
conditioning on $S_n=x \in \mathbb{R}^{d}$  yields $\mathbf{E}[ S_{k} \mid S_n=x] =kx / n$. Thus standard Gaussian concentration implies that
for any $\delta>0$   there is a constant $c_\delta>0$ satisfying
\begin{equation}\label{eq-smaple-path-RW}
\sup_{x \in \mathbb{R}^{d}} \max_{1 \leq k \leq n}\mathbf{P} \Bigl(  \Big\|S_{k}-\frac{k x}{n} \Big\| >\delta n  \,\Bigm|\, S_n=x    \Bigr)
 \lesssim  e^{-c_\delta n} .
\end{equation}

The following proposition establishes     an analogous result for $Z_{n}$: given   $|Z_{n}|$, the trajectory $|Z_{k}|$ at time $k$   concentrates closely around $\frac{k}{n}|Z_{n}|$.

  \begin{proposition}\label{prop-LDP-sample-paths}
 For every $\delta \in (0,1/4)$,  there exist  a constant $C_{\refeq{e:ZkZnl}}(\delta) > 0$  and a
 sequence $\hyperref[eq-def-C_1-varepsilon-n]{\varepsilon}_n =o(1)$ (independent of $\delta$), both defined in \eqref{eq-def-C_1-varepsilon-n}, such that for all $n \geq 1$,
  \begin{equation}   \label{e:ZkZnl}
 \max_{l\in\mathtt{L}^{\mu}_{n}} \mathbf{P} \Bigl( \Bigl| |Z_k|-\frac{k}{n}|Z_n| \Bigr|>\delta n \text{ for some } k\leq n \,\Bigm|\, |Z_n|=l \Bigr) \leq \exp\{ -n[C_{\refeq{e:ZkZnl}}(\delta)-\hyperref[eq-def-C_1-varepsilon-n]{\varepsilon}_n] \} .
  \end{equation}
  Moreover $C_{\refeq{e:ZkZnl}}(\delta)$ is continuous and strictly increasing, with  $ C_{\refeq{e:ZkZnl}}(\delta)\to 0 $ as $\delta \to 0$.
  \end{proposition}

  We emphasize that the definition of $\hyperref[eq-def-C_1-varepsilon-n]{\varepsilon}_n$ is independent of $\delta$. This will be used in the proof of Lemma~\ref{lem-level-set-size} (precisely, in Subsection \ref{sec:bootstrap}).

 \begin{remark}\label{rmk-rw-LDP}
Combining Propositions~\ref{prop-LDP-sample-paths} and \ref{prop-LDP-|Zn|} yields
\begin{equation}
  \lim_{\delta \to 0} \lim_{n \to \infty} \frac{1}{n} \log \mathbf{P} \left( |  |Z_{n}|- q n | \leq \delta n , \forall 0 \leq k \leq n  \right) = - L^{*}(q) \ \text{ for any } q \in [0,1].
\end{equation}
This result is of interest in its own right, representing a special case of large deviations for the sample paths of $(|Z_{n}|)$.
For classical results on random walks in $\mathbb{R}^d$ and abstract spaces, we refer to \cite[Theorem 5.1.2]{DZ09} and \cite{DZ95}.
In analogy with \cite[Theorem 5.1.2]{DZ09}, we expect that for i.i.d. elements  $Y_{k}$   on $\mathbb{F}$ with distribution $\mu$,   the    process $ (  |\prod_{k=1}^{\lfloor nt \rfloor} Y_{k}| : t \in [0,1]  )$ in $D[0,1]$ satisfies an LDP with the good rate function
\begin{equation}
I(\phi)=\left\{\begin{array}{cl}
  \int_0^1 L^*( \phi'(q) ) \dif q, & \text { if } \phi \text{ is absolutely continuous }, \phi(0)=0 \\
  \infty & \text { otherwise.}
  \end{array}\right.
\end{equation}
\textit{Update:} Recently, such a sample path LDP was established in \cite{JL25}; however, an explicit rate function of the integral form above was not derived there.
 \end{remark}

 \begin{proof}[Proof of Proposition \ref{prop-LDP-sample-paths}]
It suffices to consider $l\in\mathtt{L}^{\mu}_{n}$ with $l\leq n-\delta n$, since the event under consideration has probability zero for $l>n-\delta n$. Indeed,
  if $l > (1 - \delta)n$, then for any $k$,
\[   | Z_k | \leq k < \frac{k}{n} l + \delta n \quad \text{ and }\quad   | Z_k  | \geq  | Z_n  | - ( n - k)  > \frac{k}{n} l - \delta n \]
We now divide the remainder of the proof into three steps according to the value of $l$.

\smallskip
\noindent
\underline{\textit{Step 1}}. We first consider the range $\frac{\delta}{2} n \leq l \leq (1-\delta)n$.
Set $\delta' :=  \delta^2 /6 $.
 We claim that
   \begin{equation}\label{eq-union-bound-1}
    \Bigl\{ \exists \, k \text{ s.t. } \big|  | Z_k | - \frac{k}{n}l   \big| > \delta n ,  |Z_n| = l  \Bigr\}
  \subset
  \bigcup_{m=1}^{n-1} \Bigl\{ \big|  | Z_m  | - \frac{m}{n}l   \big| > \delta' n,  Z_m \prec_{\mathbb{F}} Z_n,  |Z_n| = l \Bigr\}.
  \end{equation}

  To verify \eqref{eq-union-bound-1}, suppose that  $ |  | Z_k | - \frac{l}{n}k  | > \delta n$. We choose  the index $m$  to be either $k_1$ or $k_2$, with $k_1 \leq k \leq k_2$ satisfying $Z_{k_1} = Z_{k_2} = Z_k \wedge Z_n  $, to ensure
$ |  | Z_m  | - \frac{l}{n}m  | >  \delta' n.$
\begin{itemize}
  \item If $ | Z_k  | < \frac{l}{n}k - \delta n$, then $| Z_{k_2} | \leq | Z_k | < \frac{l}{n} k_2 - \delta' n$ and hence the claim follows.
  \item If instead $| Z_k | > \frac{k}{n}l + \delta n$.  We  show  that either $| Z_{k_1} | > \frac{l}{n} k_1 + \delta' n$ or $| Z_{k_2} | < \frac{l}{n} k_2 - \delta' n$.
  Otherwise, $\frac{l}{n} k_2 - \delta' n \leq |Z_{k_2}|=|Z_{k_{1}}| \leq  \frac{l}{n} k_1 + \delta' n$, which yields$\frac{l}{n} ( k_2 - k_1 ) \leq 2 \delta' n$.
Thus,
  \[
  \left| Z_k \right| \leq \left| Z_{k_1} \right| + k_2 - k_1 \leq \frac{l}{n} k_1 + \delta' n+ \frac{2n}{l} \delta'n \leq   \frac{l}{n} k+   \frac{3n}{l} \delta'n  \leq \frac{l}{n} k+ \delta n,
  \]
  which contradicts the initial assumption. This proves \eqref{eq-union-bound-1}.
\end{itemize}

Fix $l\in\mathtt{L}^{\mu}_{n}$ with $\delta n/2\leq l\leq(1-\delta)n$. Observe that $Z_k\prec_{\mathbb{F}}Z_n$ implies $|Z_k^{-1}Z_n|=|Z_n|-|Z_k|$.
 Applying the Markov property and Proposition \ref{prop-LDP-|Zn|}, and noting that $l\in\mathsf{L}_n^\mu$ and $m\in\mathsf{L}_k^\mu$ imply $l-m\in\mathsf{L}_{n-k}^\mu$, we obtain, for all $1\leq k\leq n-1$,
  \begin{align}
  &  \mathbf{P} \Bigl(  \Bigl| | Z_k | - \frac{k}{n} l \Bigr| > \delta' n, Z_k \prec Z_n \,\Bigm|\, | Z_n | = l \Bigr)  \\
   & \quad \leq
 \sum_{\substack{  m\in\mathtt{L}^{\mu}_{k}, 0 \leq  l-m \leq n-k \\   | m - \frac{k}{n}l   | > \delta' n}} \frac{ \mathbf{P} \bigl( |Z_k| = m \bigr) \mathbf{P} \bigl( | Z_{k}^{-1} Z_{n}| = l-m  \bigr)}{\mathbf{P} \bigl( |Z_{n}|=l \bigr)} \\
       &\quad \leq  \sum_{\substack{m\in\mathtt{L}^{\mu}_{k},\ 0\leq l-m\leq n-k \\ |m-\frac{k}{n}l|>\delta'n}} \exp\Bigl\{-n\Bigl[\frac{k}{n}L^{*}\Bigl(\frac{m}{k}\Bigr)+\frac{n-k}{n}L^{*}\Bigl(\frac{l-m}{n-k}\Bigr)-L^{*}\Bigl(\frac{l}{n}\Bigr)\Bigr]+3\hyperref[def-widetilde-B-n]{\widetilde{B}}_{n}\Bigr\} . \label{e:ZkyZnx}
     \end{align}
  Note that  $ \frac{k}{n} \cdot  \frac{m}{k}   + \frac{n-k}{n}\cdot \frac{l-m}{n-k} = \frac{l}{n}$,
  \[   \Bigl|  \frac{l}{n} - \frac{m}{k} \Bigr|     = \frac{1}{k}   \Bigl| m- \frac{k}{n}l\Bigr| \geq \delta' \quad \text{ and } \quad   \Bigl| \frac{l}{n} - \frac{l-m}{n-k} \Bigr|     = \frac{1}{n-k}    \Bigl| m- \frac{k}{n}l\Bigr|  \geq \delta'  .  \]
  To bound \eqref{e:ZkyZnx}, we provide a uniform lower bound for the expression inside the square brackets.
For any   $1 \leq k \leq n-1$,  from the convexity of $L^{*}$, we get
    \begin{align}
      &\inf_{\substack{  m \leq k, 0 \leq  l-m \leq n-k \\   \left| m - \frac{k}{n}l   \right| > \delta' n}}   \Bigl\{ \frac{k}{n}L^*\Bigl(\frac{m}{k}\Bigr)  +  \frac{n-k}{n}L^*\Bigl(\frac{l-m}{n-k}\Bigr)- L^*\Bigl( \frac{l}{n}\Bigr)\Bigr\} \\
     &  \geq \inf_{q \in [\delta/2, 1-\delta]}  \Bigl\{ \frac{L^*\bigl(q-\delta'\bigr)+L^*\bigl(q+\delta'\bigr)}{2}-L^*\bigl(q\bigr)  \Bigr\}   =: C_{\refeq{eq-allpy-concavity}}(\delta)>0 .  \label{eq-allpy-concavity}
      \end{align}
Above, the positivity of $C_{\refeq{eq-allpy-concavity}}(\delta)$ follows by applying   Taylor's formula with mean-value forms of the remainder, and using the formula \eqref{eq-L*prime} for $(L^{*})''$.

Combining \eqref{eq-union-bound-1}, \eqref{e:ZkyZnx} and \eqref{eq-allpy-concavity}, we conclude that
\begin{equation}
     \max_{  \epsilon n\leq l\leq(1-\delta)n ,l\in\mathtt{L}^{\mu}_{n}} \mathbf{P} \Bigl( \Bigl| |Z_k|-\frac{k}{n}|Z_n| \Bigr|>\delta n \text{ for some } k\leq n \,\Bigm|\, |Z_n|=l \Bigr) \leq n^{2}e^{-nC_{\refeq{eq-allpy-concavity}}(\delta)+3\hyperref[def-widetilde-B-n]{\widetilde{B}}_{n}} \, .
\end{equation}

 \smallskip
 \noindent\underline{\textit{Step 2}}. Next we treat the case $|Z_n|=0$ when $0\in\mathtt{L}^{\mu}_{n}$. The union bound and  Markov's property yield that $\mathbf{P}(\exists k\text{ s.t. }|Z_k|>\delta n\mid Z_n=e)$ is bounded above by
 \begin{equation}
   \sum_{k=\delta n}^{n-1} \mathbf{P}(   |Z_{k}| =  \delta n \mid Z_{n}= e)  \le \sum_{k=\delta n}^{n-1}  \sum_{|x|= \delta n}  \frac{\mathbf{P}( Z_{k} = x) \mathbf{P}( Z_{n-k}= x^{-1} ) }{\P (Z_{n}= e)} \, .
\end{equation}
Since $|x|\in\mathtt{L}^{\mu}_{k} $ implies $|x^{-1}| \in \mathtt{L}^{\mu}_{n-k}$, using  Proposition~\ref{T:LLT}   and also  the concavity of $\hyperref[eq-def-Psi*]{\Psi^*}$,   the right-hand side is bounded by
 \begin{align}
  &  \sum_{k=\delta n}^{n-1} \ \sum_{|x|= \delta n, |x|\in\mathtt{L}^{\mu}_{k} } \!
 \exp \left\{  k \Psi^{*} \left( \xi(k,x) \right) + (n-k) \Psi^{*} \left(  \xi(n-k,x^{-1})\right)  -n \Psi^{*}(0)+ 3 \hyperref[def-widetilde-beta-n]{\widetilde{\beta}}_{n}   \right\} \\
      & \leq  n\sum_{|x|= \delta n}   \exp  \Bigl\{   n \Psi^{*}  \Bigl(  \frac{\Xi(x)+ \Xi(x^{-1})}{n}\Bigr)  -n \Psi^{*}(0) + 3 \hyperref[def-widetilde-beta-n]{\widetilde{\beta}}_{n}  \Bigr\}.
    \end{align}
Since $\nabla \Psi^*$ is strictly decreasing by Lemma~\ref{l:concave}, we define
   \begin{equation}\label{eq-psi*-zero}
     C_{\refeq{eq-psi*-zero}}(\delta):= \min _{\| \xi \|_{1} \geq \delta }   \left\{  \Psi^{*}(0) +   \langle \nabla \Psi^{*}(0) , \xi \rangle  - \Psi^{*}(\xi)    \right\} > 0.
   \end{equation}
   Noting $\|\Xi(x) + \Xi(x^{-1})\|_1 = 2 \delta n$ for $|x| = \delta n$, it follows that   $  \mathbf{P}( \exists k \text{ s.t. }  |Z_{k}| >   \delta n \mid Z_{n}= e)  $ is bounded above by
   \begin{equation}
     n e^{      - n  C_{\refeq{eq-psi*-zero}}(\delta)    + 2\hyperref[def-widetilde-beta-n]{\widetilde{\beta}}_{n} } \sum_{|x|= \delta n} \exp \left\{    \langle  \nabla \Psi^{*}( 0 ) , \Xi(x)+ \Xi(x^{-1}) \rangle \right\}   .
   \end{equation}

By Lemma \ref{l:concave}, we have  $\nabla \Psi^{*}( 0 ) =\psi(\ln R)= (\psi_{a}( \ln  R))_{a \in \mathcal{A}} $. Together with $F_{a}(R)=F_{a^{-1}}(R)$ and $\Xi_{a}(x)=\Xi_{a^{-1}}(x^{-1})$, we obtain
\begin{equation}
 D_{m}:= \sum_{x \in \mathbb{F}_{m}} \exp  \bigl\{ \langle  \nabla \Psi^{*}(0) , \Xi(x)+ \Xi(x^{-1}) \rangle  \bigr\}   =   \sum_{x \in \mathbb{F}_{m}} \exp \Bigl\{   2 \sum_{a \in \mathcal{A}} \Xi_a(x) \psi_a(\ln R) \Bigr\}  .
\end{equation}
This gives $D_{m} = \mathsf{Z}_{m}(2\psi(\ln R))$, adopting the notation from \S\ref{sec-asymp-F}.  Lemmas \ref{lem-eigenvalue} and \ref{lem-rho2psi} yields
\begin{equation}
   \lim_{m \to \infty }\frac{1}{m}  \ln D_{m} =  \varrho(2 \psi(\ln R)) =0  .
\end{equation}
\phantomsection\label{def-widetilde-D-n}
Now  let $\hyperref[def-widetilde-D-n]{\widetilde{D}}_{n} := \sup_{1 \leq k \leq n} |\ln D_{k}|$. We therefore conclude that $ \hyperref[def-widetilde-D-n]{\widetilde{D}}_{n} =o(n)$, and
\begin{equation}\label{eq-ldp-zn=e}
  \mathbf{P} \bigl(  \exists k \text{ s.t. }  |Z_{k}| >   \delta n \mid  Z_{n}= e  \bigr) \leq  n \exp \bigl\{   -n  C_{\refeq{eq-psi*-zero}}(\delta)  + 3\hyperref[def-widetilde-beta-n]{\widetilde{\beta}}_{n} + \hyperref[def-widetilde-D-n]{\widetilde{D}}_{n}  \bigr\}.
\end{equation}

\smallskip
\noindent
 \underline{\textit{Step 3}}. It remains to consider $|Z_n|=l$ with $l\in\mathtt{L}^{\mu}_{n}\cap[1,\frac{\delta}{2}n]$.
Note that $\left| |Z_k| - \frac{k}{n}l \right| > \delta n$ implies $|Z_k| > \delta n + \frac{k}{n}l \geq \delta n$ (since $|Z_k| < \frac{k}{n}l - \delta n < 0$ is impossible). Thus, there exist $k_1 < k < k_2$ such that $Z_{k_1} = Z_{k_2} = Z_k \wedge Z_n$, satisfying
   \begin{equation}
     |Z_{k}|- |Z_{k} \wedge Z_{n}|   > \delta n - l \geq \delta n /2 .
   \end{equation}
which in turn implies $k_2-k_1\geq\delta n$. Hence, for any $x\in\mathbb{F}$ with $|x|\in\mathtt{L}^{\mu}_{n}$ and $|x|\leq\frac12\delta n$,
\begin{align}
  & \mathbf{P}\Bigl( \exists \, k \leq n \text{ s.t. } \Bigl| | Z_k | - \frac{k}{n} |x|  \Bigr| > \delta n ,   Z_{n}=x   \Bigr)   \\
&\qquad  \leq
\sum_{\substack{0 \leq k_1 < k_2 \leq n \\ k_2 - k_1 > \delta n}} \ \sum_{y \prec_{\mathbb{F}} x}  \mathbf{P}\Bigl( Z_{k_{1}} =  y ,  \exists  k_1<k<k_2    \text{ s.t. }  |y^{-1} Z_{k}| > \delta n /2,  Z_{k_{2}}= y,  Z_{n}=x  \Bigr)  \label{eq-inclusion-relation-0}
  \end{align}
If $0\notin\mathsf{L}_{k_2-k_1}^\mu$, the corresponding inner sum in \eqref{eq-inclusion-relation-0} is zero.
For any fixed $k_1,k_2$ with $0 \in \mathtt{L}^{\mu}_{k_2-k_1}$    the Markov property    yields the inner sum over $y$ is bounded above by
 \begin{align}
   &   \sum_{y \prec_{\mathbb{F}} x} \mathbf{P}\bigl( Z_{k_{1}} =  y \bigr) \mathbf{P}\bigl( Z_{k_{2}-k_{1}}= e  \bigr) \mathbf{P}\bigl( Z_{n-k_{2}}= y^{-1}x \bigr) \\
   & \qquad \qquad \qquad \times  \mathbf{P}\Bigl( \exists   k <k_{2}-k_{1}  \text{ s.t. }  | Z_{k}| > \delta n /2 \,\Bigm|\, Z_{k_{2}-k_{1}}= e  \Bigr)   \\
  &   \leq  \mathbf{P}\bigl(Z_{n}=x\bigr) \times
  n \exp \bigl\{  - (k_{2}-k_{1}) C_{\refeq{eq-psi*-zero}}(\delta/2)  + 3\hyperref[def-widetilde-beta-n]{\widetilde{\beta}}_{n} + \hyperref[def-widetilde-D-n]{\widetilde{D}}_{n}    \bigr\}.  \label{eq-bd-Zn-deltan}
 \end{align}
Above,   we   used \eqref{eq-ldp-zn=e}  and  the fact   $ \sum_{y \prec_{\mathbb{F}} x} \mathbf{P} ( Z_{k_{1}} =  y )    \mathbf{P} ( Z_{k_{2}-k_{1}}= e  )   \mathbf{P} ( Z_{n-k_{2}}= y^{-1}x ) \le \P(Z_{n}=x)$.
Combining \eqref{eq-inclusion-relation-0}  and \eqref{eq-bd-Zn-deltan} yields
 for any $x\in\mathbb{F}$ with $|x|\in\mathtt{L}^{\mu}_{n}$ and $|x|\leq\delta n/2$,
\begin{equation}
   \mathbf{P} \Bigl( \Bigl| | Z_k | - \frac{k}{n} |Z_{n}|  \Bigr| > \delta n  \text{ for some } k \leq n   \,\Bigm|\, Z_{n}=x \Bigr) \leq
   n^{4} e^{ - n \delta C_{\refeq{eq-psi*-zero}}(\delta/2) + 3 \hyperref[def-widetilde-beta-n]{\widetilde{\beta}}_{n}  + \hyperref[def-widetilde-D-n]{\widetilde{D}}_{n} } \label{eq-cond-Zn-spldp} .
\end{equation}

\smallskip
\noindent
\underline{\textit{Conclusion of the proof}}.
Integrating the results from Steps 1, 2, and 3, together with Propositions~\ref{prop-LDP-|Zn|} and \ref{T:LLT} (which ensure that quantities $\hyperref[def-widetilde-B-n]{\widetilde{B}}_{n}$, $\hyperref[def-widetilde-beta-n]{\widetilde{\beta}}_{n}$, and $n^{-1}\hyperref[def-widetilde-D-n]{\widetilde{D}}_{n}$ are $o(n)$),  we  take
\begin{equation}\label{eq-def-C_1-varepsilon-n}
 C_{\refeq{e:ZkZnl}}(\delta) := \min  \Bigl\{ C_{\refeq{eq-allpy-concavity}}\! \bigl(  \frac{\delta}{2}, \delta \bigr)   , \, \delta C_{\refeq{eq-psi*-zero}}\! \bigl( \frac{\delta}{2} \bigr)  \Bigr\}
  \ ; \
  \hyperref[eq-def-C_1-varepsilon-n]{\varepsilon}_n := \frac{10}{n}  \bigl(  \E[|\mathcal{T}_1|^2]   +  \ln n + \hyperref[def-widetilde-beta-n]{\widetilde{\beta}}_{n} + \hyperref[def-widetilde-B-n]{\widetilde{B}}_{n} + \hyperref[def-widetilde-D-n]{\widetilde{D}}_{n} \bigr) .
\end{equation}
Here, the term $\E[|\mathcal{T}_1|^2]$ is included in the definition of $\hyperref[eq-def-C_1-varepsilon-n]{\varepsilon}_n$ solely for convenience in the next subsection. With these choices, \eqref{e:ZkZnl} holds, completing the proof.
 \end{proof}

 \subsection{Truncated particle counts}
 \label{sec:2nd-moment}
 For each $u\in\mathcal{T}_n$, let $u_k$ denote its ancestor at
generation $k\leq n$. We define, for $\delta>0$ and $g,x\in \mathbb{F}$,
\begin{equation}\label{eq-def-N-x-delta}
\mathcal{N}_{n,x,\delta}(g)
:=
\sum_{u\in\mathcal{T}_n}
\ind{
V(u)=x,\,
\lvert |g^{-1}V(u_k)|
-\frac{k}{n}|g^{-1}x|\rvert
\leq\delta n
\text{ for all }k\leq n
}.
\end{equation}
Translation invariance gives
$
(\mathcal{N}_{n,x,\delta}(g),\P _g)
\overset{\mathrm{law}}{=}
(\mathcal{N}_{n,g^{-1}x,\delta},\P )$.
For $0\leq m\leq n$, define
\begin{align}\label{eq-def-N-m-g-delta}
\mathcal{N}_{n,m,\delta}(g)
 :=
\sum_{x\in \mathbb{F}_m(g)}
\mathcal{N}_{n,x,\delta}(g) \quad , \quad
\mathcal{N}_{n,m,\delta}^{ \mathbb{F}}(g)
 :=
\sum_{x\in \mathbb{F}_m(g)}
\ind{\mathcal{N}_{n,x,\delta}(g)\geq1}.
\end{align}
As before, we omit $g$ when $g=e$.
The following lemma gives a lower bound on the probability that
$\mathcal{N}_{n,m,\delta}$ is at least a fixed fraction of
$\E [\mathcal{N}_{n,m}]$.

\begin{lemma}\label{lem-concentration-eta}
For every $\delta \in (0,1/4)$, let $C_{\refeq{eq-lb-Nnmdelta}}(\delta)>0$ be the constant defined in \eqref{eq-C-2-delta}; it satisfies
$C_{\refeq{eq-lb-Nnmdelta}}(\delta) \to 0$ as $\delta \to 0$. Then, for every $n$ such that $2\hyperref[eq-def-C_1-varepsilon-n]{\varepsilon}_n \leq C_{\refeq{e:ZkZnl}}(\delta)$ and every $0 \leq m \leq n$,
\begin{equation}\label{eq-lb-Nnmdelta}
   \mathbf{P} \Bigl( \mathcal{N}_{n,m,\delta} \geq \frac{4}{5} \mathbf{E}\bigl[\mathcal{N}_{n,m}\bigr] \Bigr)
   \geq
   \exp \bigl\{  -n [C_{\refeq{eq-lb-Nnmdelta}}(\delta)+\hyperref[eq-def-C_1-varepsilon-n]{\varepsilon}_n] \bigr\}
   \min\bigl\{ 1, r^{n} e^{-n L^{*} (\frac{m}{n} )} \bigr\}.
\end{equation}
 \end{lemma}

 The factor $\min\{1,r^ne^{-nL^{*}(m/n)}\}$ reflects two regimes. For $m\in\mathtt{L}^{\mu}_{n}$, the many-to-one formula and Proposition~\ref{prop-LDP-|Zn|} give $\mathbf{E}[\mathcal{N}_{n,m}]=r^ne^{-nL^{*}(m/n)+o(n)}$. If $L^{*}(m/n)<\ln r-\epsilon$, this factor equals $1$ for all sufficiently large $n$, and the lower bound is determined by the concentration cost. If $L^{*}(m/n)>\ln r+\epsilon$, then $\mathbf{P}(\mathcal{N}_{n,m}\geq1)\leq\mathbf{E}[\mathcal{N}_{n,m}]$, which decays exponentially. Off the support, both $\mathcal{N}_{n,m}$ and its expectation vanish.

We also need an analogous estimate for the fixed-site truncated count
$\mathcal{N}_{n,x,\delta}$.

\begin{lemma}\label{lem-concentration-Nnx}
Let  $r \in(R,\infty)$. There is a constant $C_{\refeq{eq-N-x-delta-C-3}}>0 $   defined in \eqref{eq-def-C-3} such that  for any $\delta \in (0,1/4)$,     and for each $n$ satisfying $2 \hyperref[eq-def-C_1-varepsilon-n]{\varepsilon}_n  \leq C_{\refeq{e:ZkZnl}}(\delta)$ and for any $x \in \mathbb{F}$ with $ |x| \leq \delta n/2$,
 \begin{equation}\label{eq-N-x-delta-C-3}
   \mathbf{P}  \Bigl( \mathcal{N}_{n,x,\delta}   \geq  \frac{4}{5} \mathbf{E}  [\mathcal{N}_{n,x}  ]   \Bigr)
   \geq \exp \bigl\{ -n[  C_{\refeq{eq-N-x-delta-C-3}} \delta + \hyperref[eq-def-C_1-varepsilon-n]{\varepsilon}_n ]   \bigr\}  .
 \end{equation}
  \end{lemma}

  The rest of this subsection is devoted to establish Lemmas \ref{lem-concentration-eta} and \ref{lem-concentration-Nnx}.

\begin{proof}[Proof of Lemma \ref{lem-concentration-eta}]
If $m\notin\mathtt{L}^{\mu}_{n}$, then $\mathcal{N}_{n,m,\delta}=\mathcal{N}_{n,m}=0$ almost surely, so \eqref{eq-lb-Nnmdelta} is trivial. Hence, in the rest of the proof we assume $m\in\mathtt{L}^{\mu}_{n}$.
We use the classical second moment argument: For any non-negative random variable $Y$ satisfying $\E Y>0$ and $\E[Y^2]\in(0,\infty)$,
\begin{equation}\label{eq-Cauchy-ineq}
 \P(Y \geq \lambda \E Y) \geq(1- \lambda)^2 \frac{[\E Y]^2}{\E[Y^2 ]}\,, \   \forall \ \lambda \in (0,1).
\end{equation}

 \smallskip
\noindent\underline{\textit{Step 1}}.
For each $u\in\mathcal{T}_n$, let $A_u=A_u(n,m)$ be the event that $|V(u)|=m$ and
$| |V(u_k)|-\frac{m}{n}k|\leq\delta n$ for all $k\leq n$. Then, it holds
\[ \mathcal{N}_{n,m,\delta}=\sum_{u\in\mathcal{T}_n}\ind{A_u}. \]
To apply \eqref{eq-Cauchy-ineq} to $\mathcal{N}_{n,m,\delta}$, we first estimate its expectation. The many-to-one lemma, together with Propositions~\ref{prop-LDP-|Zn|} and \ref{prop-LDP-sample-paths}, gives
\begin{align}
    \mathbf{E} [  \mathcal{N}_{n,m,\delta} ] &=  r^{n}  \mathbf{P}(|Z_{n}|= m )   \Bigl[   1   -\mathbf{P}  \Bigl(   \Bigl| |Z_{k}|- \frac{m}{n} k  \Bigr|  > \delta n  \text{ for some } k \leq n  \,\Bigm|\,  |Z_{n}|= m \Bigr)  \Bigr]  \\
    & \geq    r^{n} e^{-nL^{*}(m/n)-\hyperref[def-widetilde-B-n]{\widetilde{B}}_{n}}   \Bigl(   1-  e^{-n[C_{\refeq{e:ZkZnl}}(\delta) - \hyperref[eq-def-C_1-varepsilon-n]{\varepsilon}_n  ]}  \Bigr) .  \label{eq-Nnm-exp-lb}
\end{align}
We next derive an upper bound on $\E[\mathcal{N}_{n,m,\delta}^2]$. Conditioning on the genealogy  $\mathcal{T}$, we obtain
\begin{equation}\label{eq-conditional-2nd-moment}
   \mathbf{E}[  \mathcal{N}_{n,m,\delta}   ^2 \mid \mathcal{T} ]
   \leq  \sum_{u \in \mathcal{T}_{n} } \sum_{ s=1}^{n}  \, \sum_{  v \in \mathcal{T}_{n}:  |v \wedge u|=n-s }     \mathbf{P}(A_{u} \cap  A_{v}  \mid \mathcal{T})
   +  \mathbf{E}[  \mathcal{N}_{n,m,\delta}   \mid \mathcal{T}   ].
\end{equation}

 \smallskip
\noindent\underline{\textit{Step 2}}.
We claim that, whenever $A_{v}$, $v \in \mathcal{T}_{n}$ occurs, it  holds  $\left| |V(v_{k})^{-1} V(v)| - \frac{m}{n}(n-k) \right| \leq  3\delta n $ for each $k \leq n$. To see this, fix $k\leq n$ and choose $k_1\leq k\leq k_2$ such that
$V(v_{k_1})=V(v_{k_2})=V(v_k)\wedge V(v)$. For the upper bound,
\begin{align}
  |V(v_{k})^{-1} V(v)| &=| |V(v_{k_{2}})|- |V(v_{k})| | + | |V(v)|- |V(v_{k_{2}})| | \\
  & \leq \frac{m}{n} k_{2} + \delta n -\Bigl(\frac{m}{n}k -\delta n \Bigr) + m -\Bigl(\frac{m}{n}k_{2}-\delta n \Bigr) = \frac{m}{n}(n-k) + 3 \delta n .
\end{align}
For the lower bound,
\begin{equation}
  |V(v_{k})^{-1} V(v)|   \geq   |V(v)|- |V(v_{k_{1}})|  \geq  m-  \Bigl(    \frac{m}{n} k_{1} + \delta n  \Bigr)  \geq \frac{m}{n}(n-k) - 3 \delta n .
\end{equation}

 \smallskip
\noindent\underline{\textit{Step 3}}.
Now fix  $v\in\mathcal{T}_n$ such that $|v\wedge u|=n-s$. By the branching property, conditionally on $\mathcal{T}$, the process
$\left(V(v_{n-s})^{-1}V(v_{n-s+j})\right)_{1\leq j\leq s}$
is independent of $A_u$. The preceding claim implies
\begin{equation}
  \label{eq:A-u-A-v}
   \mathbf{P} ( A_u \cap A_v  \mid  \mathcal{T}  ) \leq \mathbf{P}  ( A_u \mid \mathcal{T}  ) \, \mathbf{P}  \Bigl(  \frac{m}{n}s - 3\delta n \leq  | Z_s  | \leq \frac{m}{n}s + 3\delta n  \Bigr).
\end{equation}
For $s\geq1$, the expected number of ordered pairs $(u,v)\in
\mathcal{T}_{n}^{2}$ satisfying $|u\wedge v|=n-s$ is bounded as follows:
\begin{equation}\label{eq:pair-d-s}
    r^{n-s} \times \mathbf{E} [|\mathcal{T}_{1}|(|\mathcal{T}_{1}|-1) ]\times
 r^{2(s-1)} \le \, \mathbf{E} [|\mathcal{T}_{1}|^2   ]  r^{n+s}.
\end{equation}
Substituting  \eqref{eq:pair-d-s} and  \eqref{eq:A-u-A-v} into
\eqref{eq-conditional-2nd-moment}, since $\mathbf{P}  ( A_u  )$ does not really depend on $u$,   we obtain
\begin{align}
  \mathbf{E}  \bigl[   \mathcal{N}_{n,m,\delta}^2 \bigr]
  &\leq \mathbf{E} [|\mathcal{T}_{1}|^2   ] \ \mathbf{E} [ \mathcal{N}_{n,m,\delta}  ]\,
  \sum_{s=0}^{n}r^{s}\mathbf{P} \Bigl(   \frac{m}{n}s-3\delta n\leq |Z_s|\leq
  \frac{m}{n}s+3\delta n  \Bigr) .
  \label{eq-sec-mon-Ndelta-1}
\end{align}

 \smallskip
\noindent\underline{\textit{Step 4}}.
To bound
$\mathbf{P} ( \frac{m}{n}s - 3\delta n \leq |Z_s| \leq \frac{m}{n}s + 3\delta n )$,
we distinguish two regimes. For $s\geq 3\sqrt{\delta}n$, rewriting this probability as
$\sum_{|k-\frac{m}{n}s|\leq 3\delta n}\mathbf{P}(|Z_s|=k)$
and applying Proposition~\ref{prop-LDP-|Zn|}, we obtain that it is bounded above by
\begin{align}
  n \exp \Bigl\{
  -\inf_{|h|\leq \frac{3\delta n}{s}}
  sL^{*} (\frac{m}{n}+h)
  +\hyperref[def-widetilde-B-n]{\widetilde{B}}_n\Bigr\}
  \leq
  n\exp  \bigl\{
  -sL^{*} (\frac{m}{n})
  +n\omega_{L^{*}}(\sqrt{\delta})
  +\hyperref[def-widetilde-B-n]{\widetilde{B}}_n
  \bigr\}.
\end{align}
Here, we used
$\frac{3\delta n}{s}\leq\sqrt{\delta}$.
For $s\leq 3\sqrt{\delta}n$, we simply use the trivial upper bound $1$.
Plugging these two bounds into \eqref{eq-sec-mon-Ndelta-1}, the sum on the right-hand side of \eqref{eq-sec-mon-Ndelta-1} is bounded above by
\begin{align}
 & n \bigg(  \sum_{s =1}^{3 \sqrt{\delta} n}  r^{s} + \sum_{s =3 \sqrt{\delta} n}^{n} r^{s} e^{-s L^{*} ( m/n )} e^{   n\omega_{L^{*}}(\sqrt{\delta})  + \hyperref[def-widetilde-B-n]{\widetilde{B}}_n  }  \bigg) \\
&\leq    n   e^{n C_{\ref{eq-lb-Nnmdelta}}(\delta) +  \hyperref[def-widetilde-B-n]{\widetilde{B}}_n}     \sum_{s = 0}^{n} r^{s} e^{-s L^{*} (m/n) }  \leq   n^2   e^{n C_{\ref{eq-lb-Nnmdelta}}(\delta) +  \hyperref[def-widetilde-B-n]{\widetilde{B}}_n}     \max \bigl\{  r^{n} e^{-n L^{*} (m/n) }, 1 \bigr\} .
\end{align}
where
\begin{equation}\label{eq-C-2-delta}
C_{\ref{eq-lb-Nnmdelta}}(\delta) := 3\sqrt{\delta} \max_{q  \in [0,1]} L^{*}(q)+ \omega_{L^{*}}(\sqrt{\delta}).
\end{equation}

Combining the this estimate with the first-moment bound in
\eqref{eq-Nnm-exp-lb}, we finally obtain
\begin{align}
\mathbf{E} [ \mathcal{N}_{n,m,\delta}^2 ]
\leq  \mathbf{E} [|\mathcal{T}_{1}|^2   ] \,  n^2  e^{n C_{\ref{eq-lb-Nnmdelta}}(\delta) + 2 \hyperref[def-widetilde-B-n]{\widetilde{B}}_n} \
\frac{ \max \{ r^{-n} e^{n L^{*} ( \frac{m}{n} )}, 1 \} }{
  1 - e^{-n[C_{\refeq{e:ZkZnl}}(\delta) - \hyperref[eq-def-C_1-varepsilon-n]{\varepsilon}_n  ]} } \
 ( \mathbf{E} [ \mathcal{N}_{n,m,\delta} ] )^2.
\label{eq-Nnm-2nd-moment}
\end{align}
Recall that $ \hyperref[eq-def-C_1-varepsilon-n]{\varepsilon}_n$ is defined by~\eqref{eq-def-C_1-varepsilon-n}. In particular,
$  \mathbf{E} [|\mathcal{T}_{1}|^2   ]   n^2 e^{2 \hyperref[def-widetilde-B-n]{\widetilde{B}}_n} \leq e^{n \hyperref[eq-def-C_1-varepsilon-n]{\varepsilon}_n}$.
Suppose that $n\geq1$ and
$2 \hyperref[eq-def-C_1-varepsilon-n]{\varepsilon}_n \leq C_{\refeq{e:ZkZnl}}(\delta)$.
Then
$n [ C_{\refeq{e:ZkZnl}}(\delta) - \hyperref[eq-def-C_1-varepsilon-n]{\varepsilon}_n ]
\geq n \hyperref[eq-def-C_1-varepsilon-n]{\varepsilon}_n \geq 10$.
Applying \eqref{eq-Cauchy-ineq} to $\mathcal{N}_{n,m,\delta}$ with
$\lambda=9/10$ therefore yields
\[
\mathbf{P} \Bigl( \mathcal{N}_{n, m, \delta} \geq \frac{9}{10} \mathbf{E} [ \mathcal{N}_{n, m, \delta} ]  \Bigr)
> e^{-n [ C_{\ref{eq-lb-Nnmdelta}}(\delta) + \hyperref[eq-def-C_1-varepsilon-n]{\varepsilon}_n ]} \min \bigl\{ r^n e^{- n L^{*} ( \frac{m}{n} )}, 1 \bigr\} , \qquad \forall\,m\in\mathtt{L}^{\mu}_{n}.
\]
Finally,
$\mathbf{E} [  \mathcal{N}_{n, m, \delta} ] \geq \mathbf{E} [  \mathcal{N}_{n, m} ] ( 1 - e^{ -n [ C_{\refeq{e:ZkZnl}}(\delta) - \hyperref[eq-def-C_1-varepsilon-n]{\varepsilon}_n   ]} ) \geq \frac{8}{9} \mathbf{E} [ \mathcal{N}_{n, m} ]$.
Since $(9/10)(8/9)=4/5$, the desired result follows.
\end{proof}

\begin{proof}[Proof of Lemma \ref{lem-concentration-Nnx}]
If $|x|\notin\mathtt{L}^{\mu}_{n}$, then $\mathcal{N}_{n,x,\delta}=\mathcal{N}_{n,x}=0$ almost surely, and the claimed inequality is immediate. Hence we assume $|x|\in\mathtt{L}^{\mu}_{n}$.
We follow the proof of Lemma~\ref{lem-concentration-eta} and record only the modifications needed for a fixed site $x$. By the many-to-one formula, Proposition~\ref{T:LLT}, and \eqref{eq-cond-Zn-spldp},
\begin{align}
  \mathbf{E} [ \mathcal{N}_{n,x,\delta} ]
  &= r^{n}\mathbf{P}(Z_{n}=x)
  \Bigl[  1-\mathbf{P}\Bigl( | |Z_{k}|-\frac{k}{n}|Z_n| \,\Bigm|\,
  >\delta n \text{ for some } k\leq n
  \mid Z_{n}=x   \Bigr) \Bigr]  \\
  &\geq r^{n}e^{
    n\Psi^{*}( \xi(n,x) )
    -\hyperref[def-widetilde-beta-n]{\widetilde{\beta}}_{n}}
  \bigl(  1-e^{-n[C_{\refeq{e:ZkZnl}}(\delta)
    -\hyperref[eq-def-C_1-varepsilon-n]{\varepsilon}_n]} \bigr) \\
  &\geq r^{n}e^{n\Psi^{*}(0)}
  e^{-n\delta\max_{|\xi\|\leq 1/2}|\nabla\Psi^{*}(\xi)|
    -\hyperref[def-widetilde-beta-n]{\widetilde{\beta}}_{n}}
  \bigl(  1-e^{-n[C_{\refeq{e:ZkZnl}}(\delta)
    -\hyperref[eq-def-C_1-varepsilon-n]{\varepsilon}_n]}  \bigr).
  \label{eq-1st-N-x-delta}
\end{align}
Here, in the last inequality we used
$\max_{\|\xi\|_{1}\leq\delta}
|\Psi^{*}(\xi)-\Psi^{*}(0)|
\leq
\max_{\|\xi\|_{1}\leq1/2}
|\nabla\Psi^{*}(\xi)|\delta$.

The same genealogical decomposition used to derive
\eqref{eq-sec-mon-Ndelta-1}
gives
\begin{equation}\label{eq-2nd-N-x-delta}
  \mathbf{E} [\mathcal{N}_{n,x,\delta}^{2}]
  \leq
  \mathbf{E} [|\mathcal{T}_{1}|^2   ]\ \mathbf{E} [\mathcal{N}_{n,x,\delta}]\
  \sum_{s=0}^{n}r^{s}
  \sum_{\substack{
    y: | |y|-\frac{n-s}{n}|x| |\leq3\delta n \\
    |y|\leq n-s,\ |y^{-1}x|\leq s}}
  \mathbf{P}(Z_s=y^{-1}x).
\end{equation}
By our assumption $|x|\leq\delta n/2$, the constraint in the inner sum implies
$|y|\leq4\delta n$.
 Terms with $|y^{-1}x|\notin\mathtt{L}^{\mu}_{s}$ vanish. Applying Proposition~\ref{T:LLT} to the remaining terms,   and using $\Psi^{*}(\xi)\leq\Psi^{*}(0)$ from Lemma~\ref{l:concave} and the crude bound $\# \{y \in \mathbb{F}: |y| \le 4 \delta n\} \le (2d)^{4 \delta n}$, we obtain
\begin{equation}
\sum_{ |y|\leq4\delta n }
\mathbf{P}(Z_s=y^{-1}x)
\le
\sum_{ |y|\leq4\delta n } e^{
    s\Psi^{*}( \xi(s,y^{-1}x) )
    +\hyperref[def-widetilde-beta-n]{\widetilde{\beta}}_n }
\leq e^{
s\Psi^{*}(0)
+4\delta n\ln(2d)
+\hyperref[def-widetilde-beta-n]{\widetilde{\beta}}_n } .
\end{equation}
Plugging this estimate into \eqref{eq-2nd-N-x-delta}, using
\eqref{eq-1st-N-x-delta} and the fact  $re^{\Psi^*(0)}={r}/{R}>1$, we obtain
\begin{align}
  \mathbf{E}[\mathcal{N}_{n,x,\delta}^{2}]
  &\leq
  \mathbf{E} [|\mathcal{T}_{1}|^2   ]\,\mathbf{E}[\mathcal{N}_{n,x,\delta}]
  e^{4\delta n\ln(2d)
    +\hyperref[def-widetilde-beta-n]{\widetilde{\beta}}_{n}}
  \sum_{s=0}^{n}r^{s}e^{s\Psi^{*}(0)} \\
  &\leq
  \mathbf{E} [|\mathcal{T}_{1}|^2   ]\, n
  e^{\delta nC_{\refeq{eq-N-x-delta-C-3}}
    +2\hyperref[def-widetilde-beta-n]{\widetilde{\beta}}_{n}}\,
  (1-e^{-n[C_{\refeq{e:ZkZnl}}(\delta)
    -\hyperref[eq-def-C_1-varepsilon-n]{\varepsilon}_n]})^{-1}
  \mathbf{E}[\mathcal{N}_{n,x,\delta}]^{2}
\end{align}
where
\begin{equation}\label{eq-def-C-3}
  C_{\refeq{eq-N-x-delta-C-3}}
  :=
  4\ln(2d)
  +\max_{|\xi\|\leq1/2}|\nabla\Psi^{*}(\xi)|.
\end{equation}
Finally, since
$\mathbf{E} [|\mathcal{T}_{1}|^2   ] n
e^{2\hyperref[def-widetilde-beta-n]{\widetilde{\beta}}_{n}}
\leq
e^{n\hyperref[eq-def-C_1-varepsilon-n]{\varepsilon}_n}$,
the same application of \eqref{eq-Cauchy-ineq} as in the proof of
Lemma~\ref{lem-concentration-eta} gives the desired bound.
\end{proof}

 \subsection{Bootstrap: Proof of Lemma \ref{lem-level-set-size}}
 \label{sec:bootstrap}

 Recall the truncated forward-cone counts
$\mathcal{N}_{n,m,\delta}(g)$ and
$\mathcal{N}_{n,m,\delta}^{\mathbb{F}}(g)$ defined in
\eqref{eq-def-N-m-g-delta}.
The following lemma, whose proof is postponed to the end of this
subsection, transfers the estimate at $e$ to every starting point.
Together with \eqref{eq:level-set-size}, applied with $\epsilon/2$,
it yields \eqref{eq:forward-cone-size}.  We then  prove \eqref{eq:level-set-size} by a bootstrap argument.

\begin{lemma}\label{lem-stochastical-domination}
Given $g \in \mathbb{F}$, integers $n \geq 1$, $0\leq m\leq n$, and a  real number $\delta \in (0,1)$, the following inequality holds:
  \begin{equation}
  \mathbf{P}_{g}( \mathcal{N}_{n,m,\delta}(g)  \geq  k) \geq \frac{1}{2}\mathbf{P}( \mathcal{N}_{n,m,\delta}   \geq 2k) \quad \text{ for all } \  k \geq 0.
  \end{equation}
  Similarly, we have $\mathbf{P}_{g}( \mathcal{N}_{n,m,\delta}^{\mathbb{F}}(g)  \geq  k) \geq \frac{1}{2}\mathbf{P}( \mathcal{N}_{n,m,\delta}^{\mathbb{F}}  \geq 2k)$ for all $k \ge 0$.
\end{lemma}

\begin{proof}[Proof of Lemma \ref{lem-level-set-size} for $b_n^{3/4}n^{1/4} \le m \le n(1-{\delta}_n)$.] We first introduce necessary notations.
For all   large $n$, choose $\delta^{(1)}_n=o(1)$ so that
$
2\hyperref[eq-def-C_1-varepsilon-n]{\varepsilon}_n
=C_{\refeq{e:ZkZnl}}(\delta^{(1)}_n).
$
This is possible as $C_{\refeq{e:ZkZnl}}(\cdot)$ is continuous and strictly increasing at $0$ by Proposition \ref{prop-LDP-sample-paths}. Moreover,
$ \hyperref[eq-def-C_1-varepsilon-n]{\varepsilon}_n
\lesssim C_{\refeq{e:ZkZnl}}(\delta^{(1)}_n) \lesssim \delta^{(1)}_n. $
Lemmas~\ref{lem-concentration-eta} and \ref{lem-stochastical-domination} therefore imply
\begin{equation}
    a_n
 :=\sup_{m/n\in\hyperref[def-I-r]{I(r)}}
\sup_{g\in \mathbb{F}} \,
\ln\mathbf{P}_g \Bigl(
\mathcal{N}_{n,m,\delta^{(1)}_n}(g)
\geq \frac{2}{5} \E [\mathcal{N}_{n,m}] \Bigr)^{-1}
 \leq  n[
C_{\refeq{eq-lb-Nnmdelta}}(\delta^{(1)}_n)
+\hyperref[eq-def-C_1-varepsilon-n]{\varepsilon}_n
] + \ln 2  .
\end{equation}
Set $a_n^*:=\max_{1\leq j\leq n}a_j =o(n)$, and define
\begin{align}
  b_n & :=\lceil
(na_n^*)^{1/2}
+\max_{1\leq j\leq n}(j\delta^{(1)}_j)
+n^{3/4}
\rceil \quad \text{so that} \quad  a_n^*\ll b_n\ll n \ , \
n\delta^{(1)}_n\leq b_n. \\
\qquad \   \hyperref[lem-level-set-size]{\delta}_n
& :=\sup_{j\geq n} \Bigl\{
\bigl(   \frac{b_j}{j}\bigr) ^{1/2}
+\omega_{L^*}\bigl(   \frac{2b_j}{j-b_j}\bigr) ^{1/2}
 \Bigr\}
 \quad \text{so that} \quad
  \hyperref[lem-level-set-size]{\delta}_n \downarrow 0 \ , \
b_n\ll n\hyperref[lem-level-set-size]{\delta}_n \ , \
n\hyperref[lem-level-set-size]{\delta}_n\to \infty. \label{def-of-deltan}
\end{align}

Recall the descendant BRW $V^w$ defined in \eqref{eq:sub-BRW}. Define
its truncated counts, analogously to
\eqref{eq-def-N-m-g-delta} and \eqref{eq-def-N-x-delta}, by\begin{align}
  \mathcal{N}_{n,m,\delta}^{w}(g )  &:= \sum_{u\in\mathcal{T}_n(w)}
\ind{ V(u) \in \mathbb{F}_m(g),
\lvert |g^{-1}V(u_{|w|+k})|-\frac{m}{n}k\rvert
\leq\delta n
\text{ for all } k\leq n
}; \\
\mathcal{N}_{n,x,\delta}^{w}   &:= \sum_{u\in\mathcal{T}_n(w)}
\ind{ V(w)^{-1}V(u) =x,
\lvert | V(w)^{-1}V(u_{|w|+k})|-\frac{|x|}{n}k\rvert
\leq\delta n
\text{ for all } k\leq n
}.
\end{align}
By the branching property and translation invariance,
conditionally on $\mathcal{F}^{V}_{|w|}$, these counts have the same
laws as $\mathcal{N}_{n,m,\delta}(g)$ under $\P_g$ and
$\mathcal{N}_{n,x,\delta}$ under $\P$, respectively.

 \smallskip
\noindent\underline{\textit{Step 1}}.
The restriction
$m\leq n(1-\hyperref[lem-level-set-size]{\delta}_n)$ and
$b_n=o(n\hyperref[lem-level-set-size]{\delta}_n)$ imply
$m\leq n-b_n$ for all sufficiently large $n$.
Let the BRW start at $e$. For all sufficiently large $n$ and
$m>b_n^{3/4}n^{1/4}$, we claim that
\begin{equation}\label{eq-case-I-inclusion}
\sum_{w\in\mathcal{T}_{b_n}}
\mathcal{ N}^{w}_{n-b_n,m-|V(w)|,\delta^{(1)}_{n-b_n}}
(V(w) )
\leq
\mathcal{N}_{n,m,\hyperref[def-of-deltan]{\delta}_n}.
\end{equation}
Indeed, the sets of descendants indexed by distinct
$w\in\mathcal{T}_{b_n}$ are disjoint. It thus
suffices to show that every particle $u \in \mathcal{T}_n$ counted by any summand on the
left-hand side is also counted by the random variable on the
right-hand side.

Write $\ell=|V(w)|\leq b_n$.
Since $\ell\in\mathtt{L}^{\mu}_{b_n}$,   $m\in \mathtt{L}^{\mu}_{n}$, $m>b_n^{3/4}n^{1/4}\gg b_n$ and $m\leq n-b_n$, we have $m-\ell\in\mathtt{L}^{\mu}_{n-b_n}$.
Thus every descendant count appearing in
\eqref{eq-case-I-inclusion} is well-defined. Clearly
the endpoint $ V(u)$  satisfies $ V(w) \prec_{\mathbb{F}} V(u)$,   and  hence
$|V(u)|=\ell+(m-\ell)=m$.  For $0\leq k\leq b_n$, both
$|V(u_k)|$ and $km/n$ are at most $b_n$. For $b_n\leq k\leq n$, the triangle inequality and the descendant path constraint give
\[
 \Bigl|  |V(u_k)|-\frac{km}{n}\Bigr|
 \leq \ell
+\Bigl| |V(w)^{-1}V(u_k)|
-\frac{k-b_n}{n-b_n}(m-\ell)\Bigr|
+\Bigl|  \frac{k-b_n}{n-b_n}(m-\ell)-\frac{km}{n} \Bigr|
\leq 4 b_n.
\]
Here we used $(n-b_n) \delta^{(1)}_{n-b_n} \le \max_{1\leq j\leq n}(j\delta^{(1)}_j) \le b_n$, and $ |\frac{k-b_n}{n-b_n} -\frac{k}{n}| \le \frac{b_n}{n}$.
Since $b_n=o(n\hyperref[lem-level-set-size]{\delta}_n)$, the required path constraint follows. This shows the preceding claim.

 \smallskip
\noindent\underline{\textit{Step 2}}.  We next  show that
\begin{equation}\label{eq-bootstrap-root-comparison}
 \min_{l\in\mathtt{L}^{\mu}_{b_n}}
\frac25\E [\mathcal{N}_{n-b_n,m- l }]
\geq
2[r^ne^{-nL^*(m/n)}]^{1-\epsilon}.
\end{equation}
 Indeed, for any $l\in\mathtt{L}^{\mu}_{b_n}$,
$
\lvert
\frac{m-l}{n-b_n}-\frac mn
\rvert
\leq\frac{2b_n}{n-b_n}$
and hence
\begin{equation}
  \label{eq:verify-I-r}
  \ln r-L^*(\frac{m-l}{n-b_n})
 \geq
\ln r-L^*(m/n)
-\omega_{L^*}(\frac{2b_n}{n-b_n})
 \geq
\hyperref[def-of-deltan]{\delta}_n
-\hyperref[def-of-deltan]{\delta}_n^2>0.
\end{equation}
The definition of $\hyperref[def-of-deltan]{\delta}_n$,
Proposition~\ref{prop-LDP-|Zn|}, and the hypothesis
$\ln r-L^*(m/n)\geq
\hyperref[def-of-deltan]{\delta}_n$
now yield, uniformly in $l\in\mathtt{L}^{\mu}_{b_n}$,
 \begin{align}
  \ln\E [\mathcal{N}_{n-b_n,m-l}]
&\geq
(n-b_n)[\ln r-L^*(m/n)]
-n \, \omega_{L^*}(\frac{2b_n}{n-b_n})
-\hyperref[def-widetilde-B-n]{\widetilde B}_n\\
&=n[\ln r-L^*(m/n)]
-o(n\hyperref[lem-level-set-size]{\delta}_n).
 \end{align}
Here we used
$ b_n=o(n\hyperref[lem-level-set-size]{\delta}_n)$
and
$ \omega_{L^*}(\frac{2b_n}{n-b_n})
\leq\hyperref[lem-level-set-size]{\delta}_n^2
=o(\hyperref[lem-level-set-size]{\delta}_n), $
and
$\hyperref[def-widetilde-B-n]{\widetilde B}_n/n
\leq\hyperref[eq-def-C_1-varepsilon-n]{\varepsilon}_n
\lesssim\delta^{(1)}_n\leq b_n/n$.

 \smallskip
\noindent\underline{\textit{Step 3}}. The inequalities \eqref{eq-case-I-inclusion}, \eqref{eq-bootstrap-root-comparison} yield that  whenever $\mathcal{N}_{n,m,\hyperref[lem-level-set-size]{\delta}_n}
\leq[r^ne^{-nL^*(m/n)}]^{1-\epsilon}$ each summand on the
left-hand side in \eqref{eq-case-I-inclusion} mus tbe less than $\frac{2}{5}$ fraction of its expectation.
Thus
conditionally on $\mathcal{F}^{V}_{b_n}$, using the branching property (see \eqref{eq:sub-BRW}),  the definition
of $a_n^*$, and the fact $1-t \le e^{-t}$ for all $t \in (0,1)$ we obtain
\[
 \mathbf{P}\bigl(
\mathcal{N}_{n,m,\hyperref[lem-level-set-size]{\delta}_n}
\leq[r^ne^{-nL^*(m/n)}]^{1-\epsilon}
\,\bigm|\,  \mathcal{F}^{V}_{b_n}
 \bigr)
\leq\exp\{-|\mathcal{T}_{b_n}|e^{-a_n^*}\}.
\]
 Here $\frac{m-l}{n-b_n}\in I(r)$ for $l\in\mathtt{L}^{\mu}_{b_n}$ by \eqref{eq:verify-I-r},
so the definition of $a_{n-b_n}$ applies.
Consequently, taking expectation on $\{|\mathcal{T}_{b_n}| \ge r^{b_n/2} \}$ we get
\[
 \mathbf{P}(
\mathcal{N}_{n,m,\hyperref[lem-level-set-size]{\delta}_n}
\leq[r^ne^{-nL^*(m/n)}]^{1-\epsilon}
)
\le
\mathbf{P}(|\mathcal{T}_{b_n}|\leq r^{b_n/2})
+\exp\{-r^{b_n/2}e^{-a_n^*}\}.
\]
Results on GW processes (see \cite[Lemma 1.5.1]{Mallein15}) yields  $   \mathbf{P}  (  | \mathcal{T}_{ b_{n}} |  \leq  r^{\frac{1}{2} b_{n} }  ) \leq e^{-\kappa_{\mathcal{T}} \, b_{n}}$, where  $\kappa_{\mathcal{T}}>0$ is a constant depending only on the law of $\mathcal{T}$.
Together with
$b_n\gg a_n^*$, this gives
\begin{equation}\label{eq-rooted-level-set-large-m}
\mathbf{P}(
\mathcal{N}_{n,m,\hyperref[lem-level-set-size]{\delta}_n}
\leq[r^ne^{-nL^*(m/n)}]^{1-\epsilon}
)
\leq  2\, e^{-c \, b_n}.
\end{equation}
Recall that $b_n \gg \sqrt{n}$.
 This completes the proof for the case where  $m \ge b_{n}^{3/4} n^{1/4} $.
\end{proof}

 \begin{proof}[Proof of Lemma \ref{lem-level-set-size} for $0 \le m \le b_n^{3/4}n^{1/4}  $.]

If $r\leq R$, this range contains no admissible $m$ for all
sufficiently large $n$. Indeed, the finiteness of $(L^*)'(0)$ and
$L^*(0)=\ln R$ give
$
\ln r-\ln R
 \geq
\hyperref[def-of-deltan]{\delta}_n+L^*(m/n)-L^*(0)
 \geq
(\frac{b_n}{n})^{1/2}
-O((\frac{b_n}{n})^{3/4})>0,
$
contradicting $\ln r\leq\ln R$. Thus only the case $r>R$ remains. We set
\[
c_n^-:=\lfloor b_n^{3/4}n^{1/4}\rfloor \ ,
\quad
 c_n:=\langle b_n^{5/8}n^{3/8}\rangle_{c_n^-,\mu}\ ,
\quad
\delta_n^{(2)}
:=\frac{4c_n^-}{n-c_n}+\delta_{n-c_n}^{(1)}.
\]
Thus by \eqref{eq:solve-periodicity-2} $c_n-c_n^-\in\mathtt{h}_\mu\mathbb{Z}$, and $|c_n - b_n^{5/8}n^{3/8}| \le 1$.
These scales satisfy
\begin{equation}\label{eq-step4-scales}
b_n\ll c_n^-\asymp(n-c_n)\delta_n^{(2)}\ll c_n\ll n
\quad\text{and}\quad
c_n+(n-c_n)\delta_n^{(2)}
=o(n\hyperref[def-of-deltan]{\delta}_n).
\end{equation}
Indeed, since $(b_n)$ is increasing,
$
(n-c_n)\delta_n^{(2)}
=4c_n^-+(n-c_n)\delta_{n-c_n}^{(1)}
\leq4c_n^-+b_n
\asymp c_n^-,$
while the reverse bound $(n-c_n)\delta_n^{(2)}\geq4c_n^-$ is immediate.
Moreover,
$
\frac{b_n}{c_n^-}
\asymp(\frac{b_n}{n})^{1/4}$,
$
\frac{c_n}{\sqrt{nb_n}}
\asymp(\frac{b_n}{n})^{1/8}$,
$
\frac{c_n^-}{c_n}
\asymp(\frac{b_n}{n})^{1/8}.
$
All three quantities tend to $0$. Since
$n\hyperref[def-of-deltan]{\delta}_n\geq\sqrt{nb_n}$,
these estimates prove \eqref{eq-step4-scales}.

 \medskip
\noindent\underline{\textit{Step 1}}.
Let $\mathcal{W}_n$ be the set of particles counted by
$\mathcal{N}_{c_n,c_n^-,
\hyperref[def-of-deltan]{\delta}_{c_n}}$.
For $w\in\mathcal{W}_n$, let $(V(w))_m$ be the prefix of $V(w)$ of
length $m$ and set
\[
x_w:=V(w)^{-1}(V(w))_m  \quad  \text{ so that } \ |V(w) x_w| = m  \ , \  |x_w|=c_n^--m .
\]
Since $m\in\mathtt{L}^{\mu}_{n}$ and $c_n^-\in\mathtt{L}^{\mu}_{c_n}$, it follows that $|x_w|=c_n^{-} -m\in\mathtt{L}^{\mu}_{n-c_n}$.
We claim that
\begin{equation}\label{eq-step4-inclusion}
\sum_{w\in\mathcal{W}_n}
\mathcal{N}^{w}_{n-c_n,x_w,\delta_n^{(2)}}
\leq
\mathcal{N}_{n,m,\hyperref[def-of-deltan]{\delta}_n}.
\end{equation}
It suffices to show that every particle $u$ counted by the summand indexed by $w$ on the
left-hand side is also counted by
$\mathcal{N}_{n,m,\hyperref[def-of-deltan]{\delta}_n}$. Clearly,
$|V(u)|=|V(w) x_w|$, and this common value is $m$. For $k\leq c_n$, the
distance between $|V(u_k)|$ and $km/n$ is at most
$c_n+c_n^-=o(n\hyperref[def-of-deltan]{\delta}_n)$ by
\eqref{eq-step4-scales}. For $k\geq c_n$, using
\eqref{eq-step4-scales} again, we have
\begin{align}
  \Bigl| |V(u_k)|-\frac{km}{n}\Bigr|
&\leq
|V(w)|+|V(w)^{-1}V(u_k)|+m\\
&\leq
c_n^-+|x_w|+(n-c_n)\delta_n^{(2)}+m  \le 2c_n^-+(n-c_n)\delta_n^{(2)} = o(n\hyperref[def-of-deltan]{\delta}_n) .
\end{align}
Since distinct $w\in\mathcal{W}_n$ have disjoint descendant sets,
this proves \eqref{eq-step4-inclusion}.

 \smallskip
\noindent\underline{\textit{Step 2}}. Next we show that, for all large $n$ and every $m\in\mathtt{L}^{\mu}_{n}$ with $m\leq b_n^{3/4}n^{1/4}$,
  \begin{equation}\label{eq-fixed-site-comparison}
[r^ne^{-nL^*(m/n)}]^{1-\epsilon}
\leq
 \frac{4}{5} \min_{   |x|\in [0,c_{n}^{-}] \cap  \mathtt{L}^{\mu}_{n-c_n}}
\E [\mathcal{N}_{n-c_n,x}]  .
\end{equation}
Indeed,
Proposition~\ref{T:LLT}, the continuity of $\Psi^*$ at $0$, and
$\Psi^*(0)=-L^*(0)=-\ln R$ give, uniformly for
$m\leq c_n^-$ and $|x|\leq c_n^-$ with $|x|\in\mathtt{L}^{\mu}_{n-c_n}$,
\begin{align}
  [r^ne^{-nL^*(m/n)}]^{1-\epsilon}
  &=\exp\{(1-\epsilon)n[\ln r-\ln R+o_{n}(1)]\},\\
  \E [\mathcal{N}_{n-c_n,x}]
  &=\exp\{n[\ln r-\ln R+o_{n}(1)]\}.
\end{align}
Since $\ln r-\ln R>0$,  \eqref{eq-fixed-site-comparison} follows.

Moreover, for $|x|\notin\mathtt{L}^{\mu}_{n-c_n}$ the probability in the following infimum equals $1$, whereas for supported $|x|\in\mathtt{L}^{\mu}_{n-c_n}$  Lemma~\ref{lem-concentration-Nnx} applies. Hence, we have
\begin{equation}
q_n
 :=\inf_{|x|\leq c_n^-}
\mathbf{P} \Bigl(
\mathcal{N}_{n-c_n,x,\delta_n^{(2)}}
\geq\frac{4}{5} \E [\mathcal{N}_{n-c_n,x}]
  \Bigr)  \geq
e^{
-(n-c_n)[
C_{\refeq{eq-N-x-delta-C-3}}\delta_n^{(2)}
+\hyperref[eq-def-C_1-varepsilon-n]{\varepsilon}_{n-c_n}
]}= e^{-o(c_n)} \label{eq-fixed-site-small-m}
\end{equation}
Here the  equality follows from
$
(n-c_n)\delta_n^{(2)}
+(n-c_n)\hyperref[eq-def-C_1-varepsilon-n]{\varepsilon}_{n-c_n}
\lesssim
c_n^-+b_n=o(c_n)$ by \eqref{eq-step4-scales}.
We now check the hypothesis in Lemma~\ref{lem-concentration-Nnx}. Since $\delta_n^{(2)}
\geq\delta_{n-c_n}^{(1)}$, and $C_{\refeq{e:ZkZnl}}(\cdot)$ is increasing,
we have $2 \hyperref[eq-def-C_1-varepsilon-n]{\varepsilon}_{n-c_n} =  C_{\refeq{e:ZkZnl}}(\delta_{n-c_n}^{(1)} ) \le C_{\refeq{e:ZkZnl}}(\delta_{n}^{(2)} )$.  Moreover,
$\delta_n^{(2)}\to0$, and \eqref{eq-step4-scales} gives
$|x| \le
c_n^-\leq {(n-c_n)\delta_n^{(2)}}/{2}
$
for all sufficiently large $n$.

 \medskip
\noindent\underline{\textit{Step 3}}. Recall \eqref{eq:sub-BRW}. Conditionally  on $\mathcal{F}^{V}_{c_n}$, the descendant BRWs
$\{V^w: w \in \mathcal{T}_{c_n}\}$ are independent. Combining \eqref{eq-step4-inclusion}  \eqref{eq-fixed-site-comparison} and
\eqref{eq-fixed-site-small-m} we obtain
\[
 \mathbf{P}\bigl(
\mathcal{N}_{n,m,\hyperref[lem-level-set-size]{\delta}_n}
\leq[r^ne^{-nL^*(m/n)}]^{1-\epsilon}
\,\bigm|\,  \mathcal{F}^{V}_{c_n}
 \bigr)   \le (1-q_n)^{|\mathcal{W}_n|} \,.
\]

Recall that $|\mathcal{W}_n|=\mathcal{N}_{c_n,c_n^-,
\hyperref[def-of-deltan]{\delta}_{c_n}}$.
We claim that, for all sufficiently large $n$, the large-level estimate
\eqref{eq-rooted-level-set-large-m}, applied at scale $c_n$ with
$\epsilon=1/2$, yields
\begin{equation}\label{eq-intermediate-level}
\mathbf{P} \Bigl(
\mathcal{N}_{c_n,c_n^-,
\hyperref[def-of-deltan]{\delta}_{c_n}}
\leq
 e^{ \frac{1}{3} c_n \ln (r/R) }
 \Bigr)
\leq e^{- \kappa_{\mathcal{T}} \, b_{c_n}} .
\end{equation}
The hypotheses of \eqref{eq-rooted-level-set-large-m} are satisfied.
Indeed, since $(b_n)$ is increasing and $c_n=o(n)$, we have
$ c_n^-\gg b_{c_n}^{3/4}c_n^{1/4} $
and clearly
$ c_n^-\leq
c_n(1-\hyperref[def-of-deltan]{\delta}_{c_n}). $
Moreover, since
$\ln r-L^*(c_n^-/c_n)\to\ln(r/R)>0$, for all sufficiently large $n$,
$ [r^{c_n}e^{-c_nL^*(c_n^-/c_n)}]^{1/2}
\geq e^{\frac{1}{3}c_n\ln(r/R)} $.

Combining   \eqref{eq-fixed-site-small-m} \eqref{eq-intermediate-level} with \eqref{eq-intermediate-level}
 proves, uniformly for $m\in\mathtt{L}^{\mu}_{n}$ with $m\leq c_n^-$,
\begin{equation}
  \mathbf{P}(
  \mathcal{N}_{n,m,\hyperref[def-of-deltan]{\delta}_n}
  \leq[r^ne^{-nL^*(m/n)}]^{1-\epsilon}
  )
  \leq \exp\{- \kappa_{\mathcal{T}} \, b_{c_n}\}+ \exp \bigl\{ -   e^{ \frac{1}{3} c_n [\ln (r/R) -o_n(1)] }  \bigr\}  . \label{eq-rooted-level-set-small-m}
\end{equation}
Finally, $b_n\geq n^{3/4}$ implies $b_{c_n}\gg\sqrt n$. This proves  \eqref{eq:level-set-size} for the remaining range.
\end{proof}

\begin{proof}[Proof of Lemma \ref{lem-stochastical-domination}] 
  Let $  \mathcal{N}_{n,m,\delta}^{\mathrm{all}}(g)
:=
\sum_{|g^{-1}x|=m}\mathcal{N}_{n,x,\delta}(g) $ so that   $(\mathcal{N}^{\mathrm{all}}_{n,m,\delta}(g),  \mathbf{P}_{g}) $ has the same law as $(\mathcal{N}_{n,m,\delta},  \mathbf{P})$.  Then, we have
  \begin{equation}
    \label{eq:N+-N-dif-1}
    \mathcal{N}^{\mathrm{all}}_{n,m,\delta}(g)  - \mathcal{N}_{n,m,\delta}(g) =  \sum_{ |g^{-1}x|=m} \ind{x\notin\mathbb{F}_m(g)}
\mathcal{N}_{n,x,\delta}(g).
  \end{equation}
  Let $\vartheta: \mathbb{F}\to \mathbb{F}$ be the unique group automorphism
determined by $\vartheta(a_i)=a_i^{-1}$ for $1\leq i\leq d$.
Equivalently, if $x=a_{i_1}^{\epsilon_1}\cdots a_{i_m}^{\epsilon_m}$ with $1 \leq i_j \leq d$ and $\varepsilon_j = \pm 1$ is a reduced word, then
$  \vartheta(x)=a_{i_1}^{- \epsilon_1}\cdots a_{i_m}^{- \epsilon_m}$.  Note that this map is different from the ordinary inversion map, which reverses the order of multiplication and is an anti-automorphism. 
Clearly, we have $\vartheta^2=\operatorname{id}_{ \mathbb{F}}$ and
$|\vartheta(x)|=|x|$ for every $x\in \mathbb{F}$. Using the
 the symmetric  assumption on $\mu$, we obtain
  \begin{equation}
    \label{eq:same-law}
    \bigl(   \{\widetilde{V}(u):= g \vartheta (g^{-1}V(u) ) , u \in \mathcal{T}\}, \P_{g} \bigr) \overset{\mathrm{law}}{=}  ( \{V(u),u \in \mathcal{T}\}  ,\P_{g}).
  \end{equation}
  Let $(z_{k})_{k \leq n}$ be
  any path  in $\mathbb{F}$ such that $z_0=g$, $ |g^{-1}z_{n}|=m  $, $g \not\prec_{\mathbb{F}} z_{n}$.
  We claim that the ``reflected" path $(\widetilde{z}_{k})_{k \leq n}:=(g \vartheta(g^{-1}z_{k}))_{k \leq n}$, satisfies
\begin{equation}
    \widetilde{z}_0=g \  ;  \quad   |g^{-1}\widetilde{z}_{k}|=|g^{-1}z_{k}|  ,\, \forall  k \le n \ ; \quad  g \prec_{\mathbb{F}} \widetilde{z} _{n}.  \label{eq:z-tilde-z}
\end{equation}
  Indeed  if $|g^{-1}z_n|=m$ and $g\not\prec_{ \mathbb{F}}z_n$, the reduced word
representing $g^{-1}z_n$ starts with the inverse of the last letter of
$g$. Hence the reduced word representing
$\vartheta(g^{-1}z_n)$ starts with the last letter of $g$, and this forces
$g\prec_{ \mathbb{F}}\widetilde z_n$.

Let  $ \widetilde{\mathcal{N}}_{n,m,\delta}^{\mathrm{all}}(g )$ and $ \widetilde{\mathcal{N}}_{n,m,\delta}(g )$ be defined   analogously to $ \mathcal{N}_{n,m,\delta}^{\mathrm{all}}(g)$ and  $ \mathcal{N}_{n,m,\delta}(g)$ respectively, but with respect to the process $ \widetilde{V} $.  From
  \eqref{eq:N+-N-dif-1} and \eqref{eq:z-tilde-z}, it follows
  \begin{equation}
    \label{eq:stochastic-domaince}
    \widetilde{\mathcal{N}}_{n,m,\delta}^{\mathrm{all}}(g ) - \widetilde{\mathcal{N}}_{n,m,\delta}(g ) \leq     \mathcal{N}_{n,m,\delta}(g)
  \end{equation}
  Finally, applying the pigeonhole principle, \eqref{eq:stochastic-domaince} and  \eqref{eq:same-law},  we  conclude
  \begin{align}
  \mathbf{P}( \mathcal{N}_{n,m,\delta}   \geq 2k)
  & =  \mathbf{P}_{g}( \mathcal{N}_{n,m,\delta}^{\mathrm{all}}(g)  \geq 2k)  \\
  & \leq \mathbf{P}_{g}( \mathcal{N}_{n,m,\delta}(g)  \geq k)   + \mathbf{P}_{g}( \mathcal{N}_{n,m,\delta}^{\mathrm{all}}(g) - \mathcal{N}_{n,m,\delta}(g)  \geq k) \\
  & \leq  2 \mathbf{P}_{g}( \mathcal{N}_{n,m,\delta}(g)  \geq k) .
  \end{align}
The same reflection and pigeonhole argument, applied to the distinct sites occupied by truncated particles, proves the second inequality and completes the proof.
 \end{proof}

\subsection{Proof of Lemma \ref{lem-level-set-size-F}}
\label{sec:Nnx-is-small}
The proof of Lemma~\ref{lem-level-set-size-F} is based on controlling
the maximal number of truncated particles occupying a single site.
Lemma~\ref{lem-level-set-size} shows that, with high probability,
$\mathcal{N}_{n,m,\hyperref[lem-level-set-size]{\delta}_n}
\geq
[r^ne^{-nL^*(m/n)}]^{1-\epsilon}$.
On the other hand, we shall show in Lemma~\ref{lem-Nnx-is-small}  a that the
corresponding number of particles at each fixed site is
subexponential. Consequently,
\begin{equation}
\mathcal{N}_{n,m}^{\mathbb{F}}
\geq
\frac{
\mathcal{N}_{n,m,\hyperref[lem-level-set-size]{\delta}_n}
}{
\max_{|x|=m}
\mathcal{N}_{n,x,\hyperref[lem-level-set-size]{\delta}_n}
}
\geq
\exp\{(1-\epsilon)n[\ln r-L^*(m/n)]-o(n)\}.
\end{equation}
Choosing the admissible range so that the subexponential error is
negligible compared with $n[\ln r-L^*(m/n)]$ then yields
Lemma~\ref{lem-level-set-size-F}.

\begin{lemma}\label{lem-Nnx-is-small}
Let  $r \in (1,R]$. For the sequence $ \hyperref[lem-level-set-size]{\delta}_{n}$  in Lemma \ref{lem-level-set-size}, we can find  a sequence $ \alpha_{n} $ such that $\sqrt{n} \ll \alpha_{n} \ll n$ and  for sufficiently large $n$
\begin{equation}
\sup_{|x| \leq n} \mathbf{P}( \mathcal{N}_{n,x, \hyperref[lem-level-set-size]{\delta}_{n}} \geq e^{\alpha_{n}}  ) \leq \exp\{- e^{ \alpha_{n}/11}\}.
\end{equation}
\end{lemma}

 We prove Lemma~\ref{lem-level-set-size-F} here and postpone the proof of Lemma~\ref{lem-Nnx-is-small} to the end of this section.

 \begin{proof}[Proof of Lemma~\ref{lem-level-set-size-F}]
Let $\alpha_n$ be defined in Lemma~\ref{lem-Nnx-is-small}. Define
\begin{equation}\label{def-of-tilde-deltan}
\hyperref[lem-level-set-size-F]{\tilde{\delta}}_n
:=
\sup_{j\geq n} \Bigl\{
\hyperref[lem-level-set-size]{\delta}_j
+ (\alpha_j/j)^{1/2}
+\omega_{L^*}(1/j)
   \Bigr\} \quad
\text{so that} \quad
\alpha_n\ll
n\hyperref[lem-level-set-size-F]{\tilde{\delta}}_n
\ll n.
\end{equation}
Since
$\hyperref[lem-level-set-size-F]{\tilde{\delta}}_n
\geq\hyperref[lem-level-set-size]{\delta}_n$, every admissible $m$ also satisfies
$m\leq n(1-\hyperref[lem-level-set-size]{\delta}_n)$. Hence Lemma~\ref{lem-level-set-size} (more precisely \eqref{eq-rooted-level-set-large-m} and \eqref{eq-rooted-level-set-small-m}) yields that for all large $n$,
\begin{equation}
  \label{eq:level-set-large-1}
\mathbf{P} \bigl(
\mathcal{N}_{n,m,\hyperref[lem-level-set-size]{\delta}_n}
\leq
[r^ne^{-nL^*(m/n)}]^{1-\epsilon/2} \bigr) \le e^{-2\sqrt{n}} .
\end{equation}
Whenever
$\mathcal{N}_{n,m,\hyperref[lem-level-set-size]{\delta}_n}
\geq
[r^ne^{-nL^*(m/n)}]^{1-\epsilon/2}$
and
$\max_{|x|\le n}\mathcal{N}_{n,x,\hyperref[lem-level-set-size]{\delta}_n}
\leq e^{\alpha_n}$  we have
\begin{equation}
\mathcal{N}_{n,m}^{\mathbb{F}}
\geq
\frac{
\mathcal{N}_{n,m,\hyperref[lem-level-set-size]{\delta}_n}
}{
\max_{|x|=m}
\mathcal{N}_{n,x,\hyperref[lem-level-set-size]{\delta}_n}
}
\geq
[r^ne^{-nL^*(m/n)}]^{1-\epsilon/2}e^{-\alpha_n} \geq
[r^ne^{-nL^*(m/n)}]^{1-\epsilon}.
\end{equation}
Here we used
$n[\ln r-L^*(m/n)]
\geq
n\hyperref[lem-level-set-size-F]{\tilde{\delta}}_n
\gg\alpha_n$.
Therefore, combining the preceding two inequalities, we derive that   $\P(\mathcal{N}_{n,m}^{\mathbb{F}}
\leq
[r^ne^{-nL^*(m/n)}]^{1-\epsilon})$ is bounded above by
\begin{equation}
   e^{-2\sqrt{n}}     +
  \sum_{|x|\leq n}
  \mathbf{P} (
  \mathcal{N}_{n,x,\hyperref[lem-level-set-size]{\delta}_n}
  \geq e^{\alpha_n} )
   \leq e^{-2\sqrt n}
  + (2d)^n\exp\{-e^{\alpha_n/11}\}
  \leq  e^{-\sqrt n}.
\end{equation}
Here, the first inequality follows from Lemma \ref{lem-Nnx-is-small}.
This proves \eqref{eq:level-set-size-F}.

Finally, applying the second inequality in Lemma~\ref{lem-stochastical-domination} to \eqref{eq:level-set-size-F}   with $\epsilon/2$ yields \eqref{eq:forward-cone-size-F}.
% It remains to prove \eqref{eq:forward-cone-size-F}. The reflection and
% pigeonhole argument used in the proof of
% Lemma~\ref{lem-stochastical-domination}, applied to distinct occupied
% sites instead of particles, gives
% \begin{align}
% &\mathbf{P}_g\left(
% \mathcal{N}_{n,m}^{\mathbb{F}}(g)
% \geq
% [r^ne^{-nL^*(m/n)}]^{1-\epsilon}
% \right) \geq
% \frac12\mathbf{P}\left(
% \mathcal{N}_{n,m}^{\mathbb{F}}
% \geq
% 2[r^ne^{-nL^*(m/n)}]^{1-\epsilon}
% \right).
% \end{align}
% Since
% $n[\ln r-L^*(m/n)]
% \geq
% n\hyperref[lem-level-set-size-F]{\tilde{\delta}}_n
% \to\infty$, for all  large $n$,
% $ 2[r^ne^{-nL^*(m/n)}]^{1-\epsilon}
% \leq
% [r^ne^{-nL^*(m/n)}]^{1-\epsilon/2}$.
% Applying the estimate already proved with $\epsilon/2$, uniformly in
% $g$, we obtain
% \[
% \mathbf{P}_g\left(
% \mathcal{N}_{n,m}^{\mathbb{F}}(g)
% \geq
% [r^ne^{-nL^*(m/n)}]^{1-\epsilon}
% \right)
% \geq
% \frac12(1- O( e^{-\sqrt n}))
% \geq\frac13.
% \]
% This proves \eqref{eq:forward-cone-size-F}.
\end{proof}

The proof of Lemma \ref{lem-Nnx-is-small} follows the framework outlined in \cite[Theorem 1.1]{AHS19}. We need the following   inequality for inhomogeneous Galton–Watson processes.

\begin{lemma}[{\cite[Proposition 2.1]{AHS19}}]\label{lem-Prop2.1-AHS19}
Let $ ( \Gamma_n  )_{n \ge 0}$ be   defined by
$\Gamma_{n+1}=\sum_{k=1}^{\Gamma_n} \nu_n^{(k)} $,
where
$\{\nu_n^{(k)}:n\geq0,\,k\geq1\}$ are mutually independent and
independent of $\Gamma_0$, and, for each $n$, $(\nu_n^{(k)})_{k\geq1}$ have the common law $\nu_n$ with  $ \mathbf{E} [\nu_n]\in (0,\infty)$.
 Let $\alpha>1$. For all $0 \leq i<n$, we assume the existence of $\lambda_i>0$ such that
$  \mathbf{E}  [ \exp(\lambda_i \nu_i) ]  \leq  \exp( \alpha \lambda_i \mathbf{E} [\nu_i])$.

Set  $\underline{\lambda}_n=\min _{0 \leq i<n} \lambda_i $ and $\bar{\lambda}_n = \max _{0 \leq i<n} \lambda_i$.
Then for all $\delta>0$ and all integer $\ell \geq 1$,
\begin{equation}
\mathbf{P} \Bigl(  \Gamma_{n} \geq \ell\max \Big\{1,(\alpha+\delta)^n   \max _{0 \leq i<n} \prod_{j=i}^{n-1}  \mathbf{E} [\nu_j] \Big\}  \,\Bigm|\, \Gamma_0=\ell  \Bigr) \leq n \, e^{ -\frac{\delta}{\alpha+\delta}  \,\ell\  \underline{\lambda}_n + \bar{\lambda}_n }.
\end{equation}
\end{lemma}

\begin{proof}[Proof of Lemma \ref{lem-Nnx-is-small}] Choose $J_{n}$ such that $  J_{n} \ln n  \ll n $ and $ J_{n}   \hyperref[lem-level-set-size]{\delta}_{n}  \ll 1 \ll J_{n} $ as $n \to \infty$. Take $\alpha_{n}$ such that
  \begin{equation}
J_{n} \ln n \ll \alpha_{n}  \ , \     {n}/{J_{n}}    \ll  \alpha_{n}  \ , \     J_{n}   \hyperref[lem-level-set-size]{\delta}_{n}n \ll \alpha_{n} , \sqrt{n} \ll \alpha_{n}   \ll n.
  \end{equation}
We divide the  time interval $[0, n]$ into intervals of length $\frac{n}{J_{n}}$.  Set $s_i:=i \frac{n}{J_{n}}$ for $0 \leq i \leq J_{n}$. For simplicity, we treat $s_{i}$  as an integer (rigorously, the upper integer part should be used).
Fix arbitrary $n \geq 1$ and $x \in \mathbb{F}$ with $|x| \leq n$,  let $x_{k}, k \leq |x|$  represent the word consisting of the  first $k$ letters of $x $.  Define
\[ \mathsf{PATH} = \bigl\{   f : \{s_i, 0 \leq i \leq J_{n}\} \to \left\{x_{k} :  k \leq |x| \right\} : f(0)=e,  f(n)= x\bigr\} \]
It clear that $  \# \mathsf{PATH}   \leq   n ^{J_{n}}=\mathrm{e}^{ J_{n} \ln n }$.

\smallskip\noindent
\underline{\textit{Step 1.}}
Consider the BRW.  For $1 \leq i \leq J_{n}$, a particle $u \in \mathcal{T}_{s_{i}}$ is said to follow a path $f \in \mathsf{PATH}$ until $s_{i}$  if, for all $0 \leq j \leq i$, the ancestor of the particle $u$ at generation $s_j$ lies in the ball $ B_{ \mathbb{F}} (f(s_{j}),  4n \hyperref[lem-level-set-size]{\delta}_{n} ) $. Let
 \begin{equation}
 \Gamma_{i}(f):=\text { the number of particles in $\mathcal{T}_{s_i}$ following  the path } f   \text{ until time } s_{i}.
 \end{equation}

 We assert that for each $u \in \mathcal{T}_{n}$ counted by  $\mathcal{N}_{n,x, \hyperref[lem-level-set-size]{\delta}_{n}}$, it holds $|V(u_{k})|-| V(u_{k}) \wedge x | \leq 4 n  \hyperref[lem-level-set-size]{\delta}_{n} ,\forall\, k \leq n$. Consequently $u$ follows the path $f$ defined by $f(s_{i}) = V(u_{s_{i}}) \wedge x$.

 To see this, let  $k_{1}\leq k \leq k_{2}$ satisfy $V(u_{k_{1}})=V(u_{k_{2}})= V(u_{k}) \wedge x$. Then $ \frac{|x|}{n}k_{2}-n \hyperref[lem-level-set-size]{\delta}_{n} \leq |V(u_{k}) \wedge x| \leq \frac{|x|}{n}k_{1}+n \hyperref[lem-level-set-size]{\delta}_{n} $, which implies $ \frac{|x|}{n}(k_{2}-k_{1}) \leq 2n \hyperref[lem-level-set-size]{\delta}_{n}$. Thus we get  $|V(u_{k})|-| V(u_{k}) \wedge x | \leq \frac{|x|}{n}(k -k_{1})+ 2 n \hyperref[lem-level-set-size]{\delta}_{n} \leq 4 n \hyperref[lem-level-set-size]{\delta}_{n}$, showing our claim.  As a result, we derive
 \begin{equation}
 \mathcal{N}_{n,x, \hyperref[lem-level-set-size]{\delta}_{n}}   \leq \sum_{f \in \mathsf{PATH}} \Gamma_{J_{n}}(f) .
 \end{equation}

Since by assumption $ J_{n}\ln n \ll \alpha_{n} $, for large $n$  we have $\# \mathsf{PATH} \leq e^{J_{n} \ln n} \leq e^{\alpha_{{n}}/2} $. Using the union bound yields $\mathbf{P} (  \mathcal{N}_{n,x, \hyperref[lem-level-set-size]{\delta}_{n}} \geq   e^{\alpha_{n}}  )$ is bounded above by
\begin{equation}
   \sum_{f \in \mathsf{PATH}} \mathbf{P} \Bigl(   \Gamma_{J_{n}}(f) \geq    \frac{ e^{\alpha_{n}}}{\# \mathsf{PATH} }  \Bigr)    \leq  e^{J_{n} \ln n}  \max_{f \in \mathsf{PATH} } \mathbf{P} \Bigl( \Gamma_{J_{n}}(f) \geq     e^{\alpha_{n}/2}  \Bigr)    .
\end{equation}

\smallskip\noindent
\underline{\textit{Step 2.}}
 Let $G_n$ be the event that
$ |\mathcal{T}_{s_{i+1}-s_i}(u)|
<
\exp\{\alpha_n/10\} $
for every $0\leq i<J_n$ and every $u\in\mathcal{T}_{s_i}$. Since $s_{i+1}-s_i=\frac{n}{J_n}$, using the union bound again we get
\begin{align}
\mathbf{P} (G_{n}^c  )
 \le \sum_{i=0}^{J_{n}-1} r^{s_{i}} \mathbf{P} \bigl(
  |\mathcal{T}_{\frac{n}{J_{n}}}| \geq e^{\alpha_{n}/10}  \bigr)  \leq   n r^{n} \mathbf{P} \bigl(   |\mathcal{T}_{ \frac{n}{J_{n}}   }| \geq e^{\alpha_{n}/10} \bigr) .
\end{align}
By our hypothesis   $\mathbf{E}[e^{s | \mathcal{T}_{1}|} ]<\infty$ for some $s>0$,      \cite[Theorem 1.1]{Nagaev15} applies. By taking $y_0=1$ there we obtain : with some constant $c >0$ depending on the GW tree $\mathcal{T}$,
$  \mathbf{P}( |\mathcal{T}_{n}| \geq k) \leq  2  (   1 +  c r^{-n} )  ^{-k}$   for all $ k \geq 1$.
 Since  $ n/J_{n} \ll \alpha_{n} $,
 we conclude
\begin{equation}
  \mathbf{P} (G_{n}^c  ) \leq 2 n r^{n}   \exp \bigl\{  - e^{\alpha_{n}/10}  \ln [ 1+  {c}/{r^{n/J_n}} ]  \bigr\}  =o  \bigl(   \exp \{   - e^{\alpha_{n}/11}   \}  \bigr).
\end{equation}

\smallskip\noindent
\underline{\textit{Step 3.}} Fix an arbitrary $f \in \mathsf{PATH}$.
On the event $G_{n}$, we have  $\Gamma_{i+1}(f) \leq \Gamma_{i}(f) e^{\alpha_{n}/10}$. If in addition $ \Gamma_{J_{n}}(f) \geq     e^{\alpha_{n}/2} $, there must exist  $ 1 \leq i < J_{n}$ such that $ \Gamma_{i}(f) \in [e^{\alpha_{n}/4},e^{7\alpha_{n}/20}]$. Hence it follows that
\begin{equation}
  \label{eq:jump-at-i}
  \mathbf{P} \bigl(   \Gamma_{J_{n}}(f) \geq  e^{\alpha_{n}/2} , G_{n} \bigr)  \leq \sum_{i=1}^{J_{n}} \sum_{\ell= e^{\alpha_{n}/4} }^{e^{7\alpha_{n}/20}} \mathbf{P} \bigl( \Gamma_{J_{n}}(f) \geq  e^{\alpha_{n}/2} , \Gamma_{i}(f)=\ell,  G_{n}   \bigr).
\end{equation}
 For $0 \leq j \leq J_{n}-1$, let   $\nu_k^{(j)} $ be the number of descendants   of the $k$-th particle in $\Gamma_{j}(f)$ located in $B_{\mathbb{F}}(f(s_{j+1}),4n \hyperref[lem-level-set-size]{\delta}_{n})$
 at generation $s_{j+1}$. This gives
 $\Gamma_{j+1}(f)=\sum_{k=1}^{\Gamma_{j}(f)} \nu_k^{(j)} $.

 For each $x$, let $\widetilde{\nu}^{(j)}(x)$ be the number of
generation-$n/J_n$ particles in
$B_{\mathbb{F}}(f(s_{j+1}),4n
\hyperref[lem-level-set-size]{\delta}_{n})$
in a BRW starting from $x$, and take these random variables to be
independent over $x$.
Define
\begin{equation}
  \widetilde{\nu}^{(j)}
  :=
  \min \Big\{
  \lfloor e^{\alpha_n/10}\rfloor,
  \sum_x
  \widetilde{\nu}^{(j)} (x)
  \ind{\mathrm{dist}(x,f(s_j))\leq
  4n\hyperref[lem-level-set-size]{\delta}_n}
  \Big\}.
\end{equation}

Conditionally on a particle counted by $\Gamma_j(f)$ being at $z$, its
offspring count $\nu_k^{(j)}$ has the same law as $\widetilde{\nu}^{(j)}(z)$.  Thus
we may couple $\nu_k^{(j)}$ with a random variable
$\widetilde{\nu}_k^{(j)} \overset{\mathrm{law}}{=} \widetilde{\nu}^{(j)}$ so that
$ \nu_k^{(j)}\wedge\lfloor e^{\alpha_n/10}\rfloor
  \leq \widetilde{\nu}_k^{(j)}$ almost surely.
Conditionally on the current particle
configuration, these couplings can be chosen independently for
distinct particles. Starting from
$\widetilde{\Gamma}_i(f)=\Gamma_i(f)$, define recursively
$
  \widetilde{\Gamma}_{j+1}(f)
  :=
  \sum_{k=1}^{\widetilde{\Gamma}_j(f)}
  \widetilde{\nu}_k^{(j)}$, $
  i\leq j<J_n$,
where, for each $j$, the variables
$\widetilde{\nu}_k^{(j)}$, $k\geq1$, are independent copies of
$\widetilde{\nu}^{(j)}$. On $G_n$, every relevant offspring count is
smaller than $e^{\alpha_n/10}$. Hence the preceding coupling gives,
\begin{equation}
  \Gamma_j(f)\leq\widetilde{\Gamma}_j(f) \ , \
   i\leq j\leq J_n \quad \text{ on } \   G_n.
\end{equation}

 $\widetilde{\Gamma}(f)$ is an inhomogeneous
Galton--Watson process to which Lemma~\ref{lem-Prop2.1-AHS19} applies with $\lambda_j$ taken as $e^{-\alpha_{n}/9}$ for all $j$. Indeed,
let $\epsilon_0>0$ be a  small   constant s.t. $\mathrm{e}^y \leq 1+ 2y$ for all $y \in[0, \epsilon_0]$.   Since $\widetilde{\nu}^{(j)} \leq  e^{\alpha_{n}/10}$, it follows that  $\lambda_j \widetilde{\nu}^{(j)} \leq \epsilon_0$ for all  large $n$. Thus,
\begin{equation}
\mathbf{E}[ \mathrm{e}^{\lambda_j \widetilde{\nu}^{(j)}} ]   \leq 1+2 \lambda_j \mathbf{E}[\widetilde{\nu}^{(j)} ] \leq \mathrm{e}^{2 \lambda_j \mathbf{E}[\widetilde{\nu}^{(j)} ]} .
\end{equation}
Now applying  Lemma \ref{lem-Prop2.1-AHS19}  with $\alpha=2, \delta =1$ we get
 \begin{align}
     \mathbf{P} \Bigl(  \widetilde{\Gamma}_{J_{n}}(f) \geq \ell K_{i,n}(f)   &\,\Bigm|\,\widetilde{\Gamma} _i(f)=\ell \Bigr)
  \leq J_{n} \exp \left\{- \ell e^{-\alpha_{n}/9} /3 +  1 \right\} \  , \text{ where } \\
 K_{i,n}(f) & :=
\max \Big\{
1,\,
3^{J_n-i}
\max_{i\leq k<J_n}
\prod_{j=k}^{J_n-1}
\mathbf{E}[\widetilde{\nu}^{(j)}]
\Big\}.
 \end{align}

 \smallskip\noindent
\underline{\textit{Step 4.}}
To give an upper bound for $K_{i,n}(f)$, let $K_0\coloneq\sup_{n\geq1,\,|x|\leq n}\mathbf{E}[\mathcal{N}_{n,x}]<\infty$. Indeed, if $|x|\notin\mathtt{L}^{\mu}_{n}$, this expectation is zero. If $|x|\in\mathtt{L}^{\mu}_{n}$, Proposition \ref{T:LLT} and Lemma \ref{l:concave} give $\mathbf{E}[\mathcal{N}_{n,x}]=\beta_n(\xi(n,x))\exp\{[\ln r+\Psi^{*}(\xi(n,x))]n\}$, where $\Psi^{*}(\xi)\leq\Psi^{*}(0)=-\ln R\leq-\ln r$ and $\beta_n(\xi)$ is uniformly bounded. Thus $\mathbf{E}[\widetilde{\nu}^{(j)}]$ is bounded from above by
\begin{equation}
  \sum_{ x } \sum_{y } \mathbf{E}[ \mathcal{N}_{n/J_{n}, x^{-1}y} ] \ind{\dist(x,f(s_{j}) ]\leq 4n  \hyperref[lem-level-set-size]{\delta}_{n} } \ind{ \dist (y,f(s_{j+1}) ]\leq 4n  \hyperref[lem-level-set-size]{\delta}_{n}  } \leq K_0 (2d)^{8 n  \hyperref[lem-level-set-size]{\delta}_{n}} . \label{eq-upper-bound-mj}
\end{equation}
This gives $  K_{i,n}(f)  \leq    e^{ O( J_{n} n  \hyperref[lem-level-set-size]{\delta}_{n}) }  \leq e^{\alpha_{n}/20}$   since   $ J_{n} n  \hyperref[lem-level-set-size]{\delta}_{n} \ll \alpha_{n} $. Thus, we gt
\begin{equation}
  \mathbf{P}\bigl( \widetilde{\Gamma}_{J_{n}}(f) \geq  e^{\alpha_{n}/2}    \,\bigm|\, \widetilde{\Gamma}_i(f)=\ell  \bigr)  \leq J_{n} \exp \bigl\{  -  e^{\alpha_{n}/8}  \bigr\}    \  \text{ for all }   \ e^{\alpha_{n}/4} \leq \ell \leq e^{ 7\alpha_{n}/20}  .
\end{equation}
Substituting this estimate into \eqref{eq:jump-at-i}, we conclude
\[  \mathbf{P} \bigl(  \Gamma_{J_{n}}(f) \geq  e^{\alpha_{n}/2} , G_{n}  \bigr)
    \leq  \sum_{i=1}^{J_{n}} \sum_{\ell= e^{\alpha_{n}/4} }^{e^{7\alpha_{n}/20}}  \exp \bigl\{  -  e^{\alpha_{n}/8}  \bigr\}      =o  \bigl(  \exp \bigl\{  - e^{\alpha_{n}/11}    \bigr\}  \bigr) \]
as $ n \to \infty$. We complete the proof by combining the results from Steps 1--3.
\end{proof}

\section{Proof of the Lower Bound}
\label{sec-lower-bound-hdim}
Recall that $\hyperref[def-I-r]{I(r)} := \{ q \geq 0 : \hyperref[eq-L*-via-P]{L^{*}}(q) \leq \ln r \} $.
Let
\[
 \mathfrak{I}_r :=
 \bigl\{
       [\alpha,\beta]:
      0\leq\alpha\leq\beta\leq1,\
     [\alpha,\beta]\subset I(r)
 \bigr\}.
\]
 For each $[\alpha,\beta]\in\mathfrak{I}_r$ and
$\dagger\in\{\mathcal{T},\mathbb{F}\}$, introduce the notation
\begin{align}
   S_r^{\mathbb{F}}(\alpha,\beta)
      :=\Lambda_r(\alpha,\beta) \quad &, \quad   S_r^{\mathcal{T}}(\alpha,\beta)
      :=E_r(\alpha,\beta) ;\\
  \theta_r^{\mathbb{F}} (\alpha,\beta)  := \min_{q \in [\alpha,\beta]}  \frac{ \ln r- L^{*}(q)}{q} \quad &, \quad
  \theta^{\mathcal{T}}_r (\alpha,\beta) := \min_{q \in [\alpha,\beta]}  \{  \ln r- L^{*}(q)  \} .
\end{align}
The case $\dagger=\mathbb{F}$ is considered only when $1<r\leq R$.
At $q=0$, the first quotient is interpreted according to the
convention in Theorem~\ref{thm-Hdim-Lambda-alpha-beta}.  By the
convexity of $L^*$, both minima may equivalently be taken over
$q\in\{\alpha,\beta\}$.

% For each  $[\alpha, \beta] \subset I(r)$,  we set
% \begin{align}
%   \theta^{\mathbb{F}}_{\alpha, \beta} &:= \min_{q \in [\alpha,\beta]}  \frac{ \ln r- L^{*}(q)}{q}   = \min  \Bigl\{    \frac{\ln r - L^{*}(\alpha)}{\alpha},\frac{\ln r - L^{*}(\beta)}{\beta} \Bigr\} \ ,  \text{ and }   \\
%   \theta^{\mathcal{T}}_{\alpha, \beta} &:= \min_{q \in [\alpha,\beta]}  \{  \ln r- L^{*}(q)  \} =  \ln r- \max\{ L^{*}(\alpha),    L^{*}(\beta) \} .
% \end{align}

We will employ the energy method (Lemma \ref{energy-method}) to show  that $\dim_{\mathrm{H}} S^{\dag}_{r}(\alpha,\beta)   \geq \theta^{\dag}_r (\alpha,\beta) $  holds with positive probability. To this end, we construct families of probability measures $\nu^{\dag}_{\alpha,\beta}$   supported respectively on $S^{\dag}_{r}(\alpha,\beta) $, and verify that  the finite-energy condition required in Lemma \ref{energy-method}  is satisfied with positive probability.
To strengthen “with positive probability” to “with probability one”, we need   the following zero-one law.

\begin{lemma}
\label{lem-0-1-law}
 Let
 $  \dag\in\{\mathcal{T},\mathbb{F}\}$,
where we assume $1<r\leq R$ if $\dag=\mathbb{F}$.
Let $ \mathfrak{I} \subset \mathfrak{I}_r$ be an arbitrary family of
intervals, and let
$ \theta : \mathfrak{I} \to [0,\infty)$
be a deterministic function.
Suppose that there exists a measurable event
$B^{\dag}$, determined by the BRW,
satisfying
$  \P (B^{\dag})>0$
and on  $    B^{\dag}$,
\begin{equation}\label{eq-uniform-witness}
   S_r^{\dag}(\alpha,\beta) \text{ is nonempty }   ,  \text{ and } \
   \dim_{\mathrm{H}} S_r^{\dag}(\alpha,\beta)\geq \theta ([\alpha,\beta]) \quad \text{ for all } \  [\alpha,\beta] \in \mathfrak{I} .
\end{equation}
Then there exists a measurable event $\Omega^{\dag}$ with
$\P (\Omega^{\dag})=1$ such that
\eqref{eq-uniform-witness} holds on $\Omega^{\dag}$ simultaneously for
every $[\alpha,\beta]\in \mathfrak{I}$. 
Moreover, when $\dagger=\mathcal{T}$, the same event may be chosen
so that $E_r(\alpha,\beta)$ is dense in $\partial\mathcal{T}$
simultaneously for every $[\alpha,\beta]\in\mathfrak I$.
\end{lemma}

\begin{proof}
Set $p^{\dag}:=\P (B^{\dag})>0$. For
$v\in\mathcal{T}$, let $B_v^{\dag}$ be the event
$B^{\dag}$ applied to the re-rooted descendant BRW $  \{  V(v)^{-1}V(v u): u\in \mathcal{T} \}$.  Set
\[
   \Omega^{\dag}:=\bigcup_{k\geq0}C_k^{\dag} \ , \  \text{ where }  \  C_k^{\dag}
      :=\bigcup_{v\in\mathcal{T}_k}B_v^{\dag}.
\]
Clearly $ \Omega^{\dag}$is measurable. We show that it has probability one. Indeed,
conditionally on $\{ V(v): |v| \le k \}$, the branching property implies that events
$(B_v^{\dag}:v\in\mathcal{T}_k)$ are i.i.d. and both has probability $p^{\dag}$.  Since  $|\mathcal{T}_k|\to\infty$ almost surely we get
\[
  \P \bigl((C_k^{\dag})^c\bigr) =    \E \bigl[   (1-p^{\dag})^{|\mathcal{T}_k|} \bigr]  \to 0 \quad \text{ as } k \to \infty
\]
 This implies $ \P( \Omega^{\dag})=1$.

 It remains to show \eqref{eq-uniform-witness} holds on $\Omega^{\dag}$.
Assume that $B_v^{\dag}$ occurs for some $v \in \mathcal{T}$. Then  for every $[\alpha,\beta]\in \mathfrak{I}$, we have
$ S_{r,v}^{\dag}(\alpha,\beta)\neq \emptyset,
$ and $
  \dim_{\mathrm{H}} S_{r,v}^{\dag}(\alpha,\beta)\geq \theta([\alpha,\beta]). $ Here $S_{r,v}^{\dag}$ is defined analogously to $S_{r}^{\dag}$ but with respect to the re-rooted descendant BRW $  \{  V(v)^{-1}V(v u): u\in \mathcal{T} \}$.
Set  \[
  \iota_v^{\mathbb{F}}(x):=V(v)\omega \  \text{ for } \  \omega \in \partial \mathbb{F} \quad \text{ and }  \quad   \iota_v^{\mathcal{T}}(t):=v\, t \  \text{ for } \  t \in \partial \mathcal{T} . \]
Since $\iota_v^{\dag}$ preserves Hausdorff dimension, and
$ \iota_v^{\dag}
  ( S_{r,v}^{\dag}(\alpha,\beta) )
      \subset S_r^{\dag}(\alpha,\beta)$,    \eqref{eq-uniform-witness} follows.

  Finally, when $\dagger=\mathcal T$, the density assertion follows by applying the preceding result to the descendant BRW re-rooted at
each vertex $v\in\mathcal T$. Since the Ulam--Harris tree is countable,
almost surely, simultaneously for every $v\in\mathcal T$ and every
$[\alpha,\beta]\in\mathfrak{I}$,
$  E_{r,v}(\alpha,\beta)\neq\varnothing $.
Hence $[v]_{\mathcal T}\cap E_r(\alpha,\beta)\neq\varnothing$ for every
$v\in\mathcal T$. Since these cylinders form a basis of
$\partial\mathcal T$, the set $E_r(\alpha,\beta)$ is dense.
\end{proof}

\subsection{Lower bound in Theorem \ref{thm-Hdim-Lambda-alpha-beta}}
We begin by introducing some notation and definitions. Assume $r \in (1,R]$.
Since $  \hyperref[lem-level-set-size-F]{\tilde{\delta}}_{n} \to0$, we may fix $n_0$ sufficiently
large so that, for every $n\geq1$,
\[
 \bigl\{
 q\in[0,1-2\hyperref[lem-level-set-size-F]{\tilde{\delta}}_{n+n_0} ]:
 \ln r-L^*(q)
 \geq 2  \hyperref[lem-level-set-size-F]{\tilde{\delta}}_{n+n_0}
 \bigr\}
\]
is a nonempty closed interval, which we denote by
$[q_n^-,q_n^+]$.
As $n$ increases, the interval $[q_n^{-}, q_n^{+}]$  becomes larger as  $(\hyperref[lem-level-set-size-F]{\tilde{\delta}}_{n})$ is decreasing.

\medskip \noindent
\underline{Dyadic scales.}
Let $\mathbb{D}_n$ be the uniform grid of $[q_n^-,q_n^+]$
with $2^k+1$ points, where $k=\lfloor\log_2n\rfloor$:
  \begin{equation}
    \mathbb{D}_{n} := \Bigl\{ q^{-}_{n} + \frac{j}{2^{k}}( q^{+}_{n}- q^{-}_{n})  : 0 \leq j \leq 2^{k} \Bigr\} ,
    \qquad 2^{k} \leq n < 2^{k+1},\quad k\geq0.
  \end{equation}
  It follows that $\# \mathbb{D}_{n} \leq 4n$, and   $ \Delta( \mathbb{D}_{n} ):=  \min\{ |q-q'|: q ,   q' \in \mathbb{D}_{n},  q \neq  q' \} \in [ \frac{q^{+}_{n}-q^{-}_{n}}{n} , \frac{2}{n})$.

  When $r < R$,  it is clear that  $\min_{n} q_{n}^{-} > 0$. When $r=R$, since $(L^{*})'(0)<\infty$, $q_{n}^{-}$ has the same order as
  $\hyperref[lem-level-set-size-F]{\tilde{\delta}}_{n+n_0}$, and hence
  $n \min_{q \in \mathbb{D}_{n}} q = n q_{n}^{-} \to \infty$.  Moreover,
 \begin{equation}
  \bigcup_{n \geq 1}\mathbb{D}_{n} \  \text{ is a dense subset of } \ I(r).
 \end{equation}

 \smallskip\noindent
\underline{Space--time scales.}
Let
$m_n=(m_1+n-1)^4$ and $M_n:=\sum_{j=1}^n m_j$,
where $m_1$ will be chosen sufficiently large later.  Recall   \eqref{eq:solve-periodicity-2}. Fix
$\epsilon\in(0,1/4)$.
For $q\in\mathbb{D}_n$, set
\begin{equation}\label{eq-def-elln-An}
\ell_n(q):=\langle m_n/q\rangle_{m_n,\mu}  \quad \text{ and } \quad
 A_n(q) :=
   \Big\lfloor
       \frac{1}{6}\Bigl(
       r^{\ell_n(q)}
       \exp \Bigl\{
          -\ell_n(q)L^* \bigl(\frac{m_n}{\ell_n(q)} \bigr) \Bigr\}
       \Bigr)  ^{1-\epsilon/2}
      \Big\rfloor .
\end{equation}
Since $q\in\mathbb{D}_n$ satisfies $0<q\leq1$, \eqref{eq:solve-periodicity-2} gives
$\ell_n(q)\geq m_n$, $m_n\in\mathsf{L}_{\ell_n(q)}^\mu$, and
$|\ell_n(q)-m_n/q|\leq1$. Thus the exact level set
$\mathcal{N}^{\mathbb{F}}_{\ell_n(q),m_n}$ is non-degenerate.

By increasing $m_1$, we may assume, uniformly in $n\geq1$ and
$q\in\mathbb{D}_n$, that
$\ell_n(q)\geq n+n_0$,
$1/(m_n-1)\leq\hyperref[lem-level-set-size-F]{\tilde{\delta}}_{n+n_0}$, and
$\omega_{L^*}(1/(m_n-1))\leq\hyperref[lem-level-set-size-F]{\tilde{\delta}}_{n+n_0}$.
Since
$|m_n/\ell_n(q)-q|\leq1/(m_n-1)$,
$\ln r-L^*(m_n/\ell_n(q))
\geq\hyperref[lem-level-set-size-F]{\tilde{\delta}}_{n+n_0}$.
The definition of
$\mathbb{D}_n$ and the monotonicity of
$(\hyperref[lem-level-set-size-F]{\tilde{\delta}}_n)$ imply that
\begin{equation}\label{eq-grid-conditions-fecund}
 \frac{m_n}{\ell_n(q)}
 \leq1-\hyperref[lem-level-set-size-F]{\tilde{\delta}}_{\ell_n(q)},
 \qquad
 \ln r-L^*\Bigl(  \frac{m_n}{\ell_n(q)}  \Bigr)
 \geq\hyperref[lem-level-set-size-F]{\tilde{\delta}}_{\ell_n(q)}.
\end{equation}
Thus Lemma~\ref{lem-level-set-size-F} is applicable to $\mathcal{N}^{\mathbb{F}}_{\ell_n(q),m_n}$ for every $n \ge 1$.

Besides,  set
$\underline A_n:=\min_{q\in\mathbb{D}_n}A_n(q)$. Combining the preceding estimates with the definition of $A_n(q)$, and using $\ell_n(q)\geq m_n$ and
$m_n\hyperref[lem-level-set-size-F]{\tilde{\delta}}_{n+n_0}\to\infty$ (recall
$\hyperref[lem-level-set-size-F]{\tilde{\delta}}_j\geq\delta_j\geq j^{-1/8}$),  we may, after increasing $m_1$ if necessary, find a constant $\underline a>0$ such that, for every $n\geq1$,
\begin{equation}\label{eq-lower-An}
 \underline A_n \ge 2  \quad \text{ and } \quad  \log\underline A_n
 \geq
 \underline a\,m_n\hyperref[lem-level-set-size-F]{\tilde{\delta}}_{n+n_0}.
\end{equation}
In particular,
the sequence $(\underline A_n)$ grows faster than any polynomial in
$n$.

\medskip \noindent
\underline{Fecund particles.}
For a particle $v\in\mathcal{T}$, $n \ge 1$ and $q\in\mathbb{D}_n$, define the
set of distinct forward sites  reached after $\ell_n(q)$ generations of $v$ with average speed $q$   by
\begin{equation}\label{eq-def-raw-descendants}
 \mathbb{X}_n(v,q)
 := \bigl\{
 x\in\mathbb{F}_{m_n}(V(v)):
 \exists\,u\in\mathcal{T}_{\ell_n(q)}(v)
 \text{ s.t. }V(u)=x \bigr\} .
\end{equation}
We say  a particle $v \in \mathcal{T}$ is \emph{$q$-fecund for  stage $n$,} if
 $ \# \mathbb{X}_n(v,q) \geq 6 \, A_n(q)$.

\medskip \noindent
\underline{Representatives and continuation.}
  We fix the lexicographical order on $\mathcal{T}$, which is a total order.  (This is different from the partial order $\prec_{\mathcal{T}}$). For every $x\in \mathbb{X}_n(v,q) $ we may set
\begin{equation}\label{eq-def-representative}
 \operatorname{rep}_n(v,q,x)
 :=
 \min_{\mathrm{lex}}
 \bigl\{
 u\in\mathcal{T}_{\ell_n(q)}(v):V(u)=x  \bigr\} .
\end{equation}

 Assume that $v$ is $q$-fecund for stage $n$.
For $q'\in\mathbb{D}_{n+1}$, we call a site
$x\in\mathbb{X}_n(v,q)$ a \emph{$q'$-continuation site} if
$\operatorname{rep}_n(v,q,x)$ is $q'$-fecund for stage $n+1$.
We denote   all such sites by
\begin{equation}\label{eq-def-good-descendants}
\mathbb{X}_n(v;q\mid q')
:=
\bigl\{
x\in\mathbb{X}_n(v,q):
\operatorname{rep}_n(v,q,x)
\text{ is $q'$-fecund for stage $n+1$}
\bigr\}.
\end{equation}
We call $\{\operatorname{rep}_n(v,q,x) ,
 x\in\mathbb{X}_n(v;q\mid q')\}$,
 the \emph{$q'$-continuation successors} of $v$.

If $v$ is not $q$-fecund for stage $n$, 
 we %% declare that it has no
%% continuation sites and
set $\mathbb{X}_n(v;q\mid q') = \emptyset$ for convenience.

 \medskip \noindent
\underline{Selected cohort and fecund skeleton.}
We recursively define a family of selected cohorts. Set
\[
   \mathcal{S}_{ \mid q_1}:=\{\emptyset\} \quad \text{for every} \quad  q_1 \in \mathbb{D}_1 .
\]
Let $n \ge 1$. Suppose that  $\mathcal{S}_{q_{1},\dots,q_{n-1} \mid q_n}$ has been defined for  every
$(q_{1},\dots,q_{n})\in\prod_{j=1}^n\mathbb{D}_j$. Then we set, for every $(q_{1},\dots,q_{n+1})\in \prod_{j=1}^{n+1}\mathbb{D}_j $,
\begin{equation}\label{eq-def-fecund-tree}
 \mathcal{S}_{q_{1},\dots,q_{n} \mid   q_{n+1}}
 :=
 \bigcup_{v\in\mathcal{S}_{q_{1},\dots,q_{n-1} \mid q_n}}
 \Bigl\{
 \operatorname{rep}_n(v,q_n,x):
 x\in\mathbb{X}_{n}(v,q_n\mid q_{n+1})
 \Bigr\}.
\end{equation}
We call
$\mathcal{S}_{q_1,\ldots,q_n\mid q_{n+1}}$
a \emph{selected cohort after stage $n$}.
Thus, it consists precisely of the
$q_{n+1}$-continuation successors of particles in the preceding
cohort.

The
genealogically linked sequence
$ (\mathcal{S}_{q_1,\ldots,q_{n-1}\mid q_{n}})_{n\ge 1}  $
is called the \emph{fecund skeleton associated with $(q_n)_{n \ge 1}$}.
By construction, every such skeleton satisfies the
following properties.
\begin{enumerate}[(i)]

\item If $n\geq2$, each particle
$v\in\mathcal{S}_{q_1,\ldots,q_{n-1}\mid q_n}$, if any, satisfies
\footnote{Although the particles lie in generation
$\sum_{j=1}^{n-1}\ell_j(q_j)$, the cohort is
$\mathcal{F}^V_{\sum_{j=1}^n\ell_j(q_j)}$-measurable, since membership
requires checking $q_n$-fecundity. All subsequent conditional
estimates use this sigma-field, so no stopping-line argument is needed.}
\begin{equation}\label{eq:auto-fecund}
v\text{ is $q_n$-fecund for stage $n$}
 \quad \text{ and } \quad
|v|=\sum_{j=1}^{n-1}\ell_j(q_j).
\end{equation}

\item No spatial collisions occur within the selected cohort:
\begin{equation}\label{eq:cohort-spatial-injectivity}
u,v\in\mathcal{S}_{q_1,\ldots,q_{n-1}\mid q_n},
\quad u\neq v
\quad\Longrightarrow\quad
V(u)\neq V(v).
\end{equation}
Indeed, successors of the same particle represent distinct sites,
whereas successors of distinct particles lie in disjoint forward
cones.

\item For  $v\in\mathcal{T}$, $q\in\mathbb{D}_n$,
$q'\in\mathbb{D}_{n+1}$,
$
\#\mathbb{X}_n(v;q\mid q')
\leq \#\mathbb{X}_n(v,q)
\leq(2d)^{m_n}$.
Consequently,
\begin{equation}\label{eq:trivial-bnd-fecund-tree}
\#\mathcal{S}_{q_1,\ldots,q_{n-1}\mid q_n}
\leq(2d)^{M_{n-1}}.
\end{equation} 
\end{enumerate}

\medskip\noindent
\underline{Good event.}
Let $B_\epsilon^{\mathbb{F}}$ be the event on which:

\begin{enumerate}[(i)]

\item the root $\emptyset$ is $q_1$-fecund for stage $1$ for every
$q_1\in\mathbb{D}_1$;

\item for every $n\geq1$, every  tuple
$(q_1,\ldots,q_{n+1})   \in \prod_{j=1}^{n+1} \mathbb{D}_j$, and every
$v\in\mathcal{S}_{q_1,\ldots,q_{n-1}\mid q_n}$,
\[
\#\mathbb{X}_n(v;q_n\mid q_{n+1})
\geq A_n(q_n).
\]

\end{enumerate}

On $B_\epsilon^{\mathbb{F}}$, for every $n\geq1$,
\begin{equation}
    \# \mathcal{S}_{q_1,\dots,q_n \mid q_{n+1}}
   \geq
   \prod_{j=1}^n A_j(q_j)  \  \text{ for every } \  (q_1,\ldots,q_{n+1}) \in\prod_{j=1}^{n+1}\mathbb{D}_j  .  \label{eq:many-fecund-particles}
\end{equation}
Indeed, the successor sets corresponding to distinct particles in
$\mathcal{S}_{q_1,\ldots,q_{n-1}\mid q_n}$ are disjoint. Therefore,
on $B_\epsilon^{\mathbb{F}}$,
$
\#\mathcal{S}_{q_1,\ldots,q_n\mid q_{n+1}}
 =
\sum_{v\in\mathcal{S}_{q_1,\ldots,q_{n-1}\mid q_n}}
\#\mathbb{X}_n(v;q_n\mid q_{n+1}) \geq
A_n(q_n)\,
\#\mathcal{S}_{q_1,\ldots,q_{n-1}\mid q_n}$.
Iterating this inequality from
$\#\mathcal{S}_{\mid q_1}=1$ proves
\eqref{eq:many-fecund-particles}.

\begin{lemma}\label{lem-fecund-subtree-F}
For any $\epsilon\in(0,1/4)$, one can choose $m_1$ sufficiently
large so that
\begin{equation}\label{eq-positive-fecund-event}
   \P \bigl(B_\epsilon^{\mathbb{F}}\bigr)\geq 1/2 .
\end{equation}
\end{lemma}

\begin{proof} By \eqref{eq-grid-conditions-fecund}, provided $m_1$ is large we may apply \eqref{eq:level-set-size-F} in
Lemma~\ref{lem-level-set-size-F} which implies
\[   \P\bigl(\emptyset\text{ is not $q$-fecund at stage $1$}
  \bigr)
   \le   \mathbf{P}   \bigl(  \mathcal{N}_{\ell_1(q),m_1}^{\mathbb{F}}  \leq \bigl[ r^{\ell_1(q) } e^{-\ell_1(q)  L^{*}(\frac{m_1}{\ell_1(q)})}  \bigr]^{1- \frac{\epsilon}{4}}  \bigr)   \leq e^{-\sqrt{\ell_1(q)}}  \]
for every $ q \in \mathbb{D}_1$. A union bound then shows that, after increasing
$m_1$ if necessary, condition~(i) fails with probability less than
$1/4$.

 Let $n \ge 1$, $v \in \mathcal{T}$ and $q \in \mathbb{D}_n$. The event $\{v \text{ is $q$-fecund at stage } n \}$ is $\mathcal{F}^{V}_{|v|+\ell_n(q)}$ measurable.
 We claim that,
for every $q'\in\mathbb{D}_{n+1}$,
\begin{equation}\label{eq-nonfecund}
 \P \bigl(
 \# \mathbb{X}_{n}(v;q\mid q') <A_n(q)
  \,\bigm|\, \mathcal{F}^{V}_{|v|+\ell_n(q)}
 \bigr) \ind{   v \text{ is $q$-fecund at stage } n }
 \leq e^{- A_n(q)/4 } .
\end{equation}
 Recall \eqref{eq:sub-BRW} and the branching property. Conditionally
on $\mathcal{F}^{V}_{|v|+\ell_n(q)}$, the descendant BRWs $(V^{w} :
 w \in \{
 \operatorname{rep}_n(v,q,x):
 x \in \mathbb{X}_n(v,q)
 \} )$
are independent.
 Moreover, by  \eqref{eq:forward-cone-size-F} in
Lemma~\ref{lem-level-set-size-F}, for every
$q'\in\mathbb{D}_{n+1}$, each representative is $q'$-fecund at stage $n+1$ with conditional probability at least
$1/3$. (Here the hypotheses of Lemma~\ref{lem-level-set-size-F} are satisfied because of \eqref{eq-grid-conditions-fecund}.)
Consequently,   conditioning on $\mathcal{F}^{V}_{|v|+\ell_n(q)}$ and   on the event that $v$ is
$q$-fecund at stage $n$,
\[
\# \mathbb{X}_{n}(v;q\mid q')
 \quad\text{stochastically dominates}\quad
 \operatorname{Bin}(6A_n(q),1/3).
\]
Using Chernoff's inequality,  we obtain    \eqref{eq-nonfecund}.

Suppose that condition~(i) holds but (ii) fails in the definition of $B^{\mathbb{F}}_{\epsilon}$. Then there must exist  $n\geq1$,
$(q_1,\ldots,q_n)\in\prod_{j=1}^n\mathbb{D}_j$,
$v\in\mathcal{S}_{q_1\cdots q_{n-1} \mid q_{n}}$, and
$q_{n+1}\in\mathbb{D}_{n+1}$, such that
\[
 \#\mathbb{X}_{n}(v,q_n \mid q_{n+1})<A_n(q_n).
\]
The union bound and   \eqref{eq-nonfecund} thus yield that $\P (\text{(i) holds, but } B^{\mathbb{F}}_{\epsilon} \text{ fails} )$ is bounded above by
\begin{equation}
    \sum_{n \ge 1} \ \sum_{ q_{j} \in  \mathbb{D}_j , j \le n+1}\!\!\! \E  \Big(   \sum_{v \in \mathcal{S}_{q_1\cdots q_{n-1} \mid q_{n}} }  \P \bigl[   \text{(i) holds}, \#\mathbb{X}_{n}(v,q_n \mid q_{n+1})<A_n(q_n)    \,\bigm|\, \mathcal{F}^{V}_{\sum_{j=1}^{n} \ell_{j}(q_{j}) } \,\bigr] \Big)
\end{equation}
Note that each summand  $v$ must be $q_n$-fecund at stage $n$, because of \eqref{eq:auto-fecund} and condition~(i). Applying  the estimate \eqref{eq-nonfecund} and \eqref{eq:trivial-bnd-fecund-tree}, and using $A_n(q) \ge  \underline{A}_n  $,  we see that   the preceding sum is bounded above by
\begin{equation}
 \sum_{n\geq1}
 |\mathbb{D}_{n+1}|^{n+1}
 (2d)^{M_{n}}
 e^{-\underline A_n/4}  \leq
 \sum_{n\geq1} \,
 e^{
 M_{n}\ln(2d)+O(n\log n)- {\underline A_n}/{4}} \ .
\end{equation}
Using \eqref{eq-lower-An}, the last sum can be made smaller than $1/4$ by
increasing $m_1$. Since condition~(i) also fails with probability less
than $1/4$, we conclude \eqref{eq-positive-fecund-event}.
\end{proof}

 \begin{proof}[Proof of the lower bound in
Theorem~\ref{thm-Hdim-Lambda-alpha-beta}]
We first reduce the proof of
Theorem~\ref{thm-Hdim-Lambda-alpha-beta} to the following claim: for
every $\epsilon\in(0,1/4)$, on $B_\epsilon^{\mathbb{F}}$,
\begin{equation}\label{eq:reduction-thm-lower-bnd}
 \Lambda_r(\alpha,\beta)\neq\emptyset
 \quad\text{and}\quad
 \dim_{\mathrm{H}}\Lambda_r(\alpha,\beta)
 \geq
 (1-\epsilon)^2 \,\theta^{\mathbb{F}}_{r}(\alpha,\beta)
 \quad\text{for every }[\alpha,\beta]\subset I(r).
\end{equation}
Indeed, Lemma~\ref{lem-fecund-subtree-F} gives
$\P(B_\epsilon^{\mathbb{F}})\geq1/2$.  Hence, once
\eqref{eq:reduction-thm-lower-bnd} is established,
Lemma~\ref{lem-0-1-law} yields an event
$\Omega_\epsilon^{\mathbb{F}}$ of probability one on which
\eqref{eq:reduction-thm-lower-bnd} holds.
Now let $\epsilon_k\downarrow0$ and set
$
 \Omega^{\mathbb{F}}
 :=
 \bigcap_{k\geq1}\Omega_{\epsilon_k}^{\mathbb{F}}$.
Then $\P(\Omega^{\mathbb{F}})=1$, and letting $k\to\infty$ gives the
desired nonemptiness and dimension lower bound.
It therefore remains
to prove \eqref{eq:reduction-thm-lower-bnd} on
$B_\epsilon^{\mathbb{F}}$. We fix an arbitrary $[\alpha,\beta]\subset I(r)$.

\smallskip\noindent
\underline{\textit{Step 1.}}
We choose $\eta_j=\eta_j^{\alpha,\beta}\in\mathbb{D}_j$ so that, with
$ L_n:=\sum_{j=1}^n\ell_j(\eta_j)$ and $L_0=0$,
we have
\begin{equation}\label{eq-choose-eta-F}
 \liminf_{n\to\infty}\frac{M_n}{L_n}=\alpha\ ,
 \quad
 \limsup_{n\to\infty}\frac{M_n}{L_n}=\beta\ ,
 \quad
  \lim_{n\to\infty} \dist\bigl(\eta_n,[\alpha,\beta]\bigr)= 0.
\end{equation}

To this end, for $q\in I(r)$, let $\pi_j(q)$ be the smallest
point of $\mathbb{D}_j$ minimizing the distance to $q$.  Then
$\pi_j(q)\to q$.
If $\alpha=\beta$, set $\eta_j=\pi_j(\alpha)$ for every $j\geq1$.
Then $\frac{M_n}{L_n}\to\alpha$, and \eqref{eq-choose-eta-F} follows.
Suppose now that $\alpha<\beta$.  Whenever $\eta_1,\ldots,\eta_n$ have been chosen,
$  {M_n}/{L_n}$ is well-defined.   Set $s_0:=0$ and construct
$(\eta_j)$ recursively.  Given $s_{2k-2}$, let
$\eta_j=\pi_j(\beta)$ for $j>s_{2k-2}$ until the first index
$ s_{2k-1}
 := \inf  \{
 n>s_{2k-2}:
\frac{M_n}{L_n}\geq
 \beta- \frac{\beta-\alpha}{10^k}   \}  . $
Next, set $\eta_j=\pi_j(\alpha)$ for $j>s_{2k-1}$ until the first
index $ s_{2k}
 :=
 \inf \{
 n>s_{2k-1}:
 \frac{M_n}{L_n}\leq
 \alpha+\frac{\beta-\alpha}{10^k}
 \}$.
If the labels are eventually chosen as $\pi_j(q)$, then
$M_n/L_n\to q$. Moreover,
$\ell_{n+1}(\eta_{n+1})/L_n\to0$, so the overshoot at each switching
time tends to zero. Hence all switching times are finite and
\eqref{eq-choose-eta-F} follows.

\smallskip\noindent
\underline{\textit{Step 2.}}
Suppose $B_\epsilon^{\mathbb{F}}$ occurs. Let
$\mathcal{R}_{\alpha,\beta}$ be the set of rays
$t\in\partial\mathcal{T}$ satisfying
\begin{equation}\label{eq-rays-through-fecund-tree}
  t_{L_n}  \in
   \mathcal{S}_{\eta_1,\dots,\eta_n \mid \eta_{n+1}}
   \quad \text{for every }n\geq1.
\end{equation}
By condition~(ii) in the definition of
$B_\epsilon^{\mathbb F}$, for every $v\in\mathcal{S}_{\eta_1,\dots,\eta_{n-1}\mid \eta_n}$,
$ \# \{
 u\in\mathcal{S}_{\eta_1,\dots,\eta_n \mid \eta_{n+1}}:
 v\prec_{\mathcal{T}}u
 \}
 =
 \#\mathbb{X}_{n}(v,\eta_n\mid \eta_{n+1})
 \geq A_n(\eta_n)\geq2 $. Consequently, $\mathcal{R} _{\alpha,\beta}$ is nonempty; it is also closed,
  being a compact  subset of $\partial\mathcal{T}$.

Starting with unit mass at the root, we define a probability measure
$\nu_{\alpha,\beta}^{\mathcal{T}}$ on
$\mathcal{R}_{\alpha,\beta}$ recursively as follows.  For every
$n\geq1$ and every
$v\in\mathcal{S}_{\eta_1,\dots,\eta_{n-1}\mid \eta_n}$, distribute the mass assigned to
$v$ uniformly among all vertices
$\{
 u\in\mathcal{S}_{\eta_1,\dots,\eta_n \mid \eta_{n+1}}:
 v\prec_{\mathcal{T}}u
 \}$.
More precisely, for
$v\in\mathcal{S}_{\eta_1,\dots,\eta_n \mid \eta_{n+1}}$, let $v_{L_j}$ be its unique
ancestor in
$\mathcal{S}_{\eta_1,\dots,\eta_{j}\mid \eta_{j+1}}$, $0\leq j\leq n$, so that
$v_0=\emptyset$ and $v_{L_n}=v$.  Then by \eqref{eq:many-fecund-particles}, on   $B_\epsilon^{\mathbb{F}}$,
\begin{equation}\label{eq-cylinder-mass-fecund}
 \nu_{\alpha,\beta}^{\mathcal{T}}
 \bigl(\mathcal{R}_{\alpha,\beta}\cap[v]_{\mathcal{T}}\bigr)
 =
 \prod_{j=1}^n
 \frac{1}{
 \#\mathbb{X}_{j}
 (v_{L_{j-1}},\eta_j \mid \eta_{j+1})
 }  \leq
 \prod_{j=1}^n\frac1{A_j(\eta_j)} \,.
\end{equation}
These cylinder masses are consistent and hence determine a probability measure
$\nu_{\alpha,\beta}^{\mathcal{T}}$ on $\mathcal{R}_{\alpha,\beta}$.
Since
$A_j(\eta_j)\geq2$, the measure is non-atomic. Recall in subection  \ref{sec-BRWonF}, it has been shown that   $V|_{\partial \mathcal{T}}:\partial\mathcal{T}\to\partial\mathbb{F}$ is an injection almost surely.
 Define
\begin{equation}\label{eq-push-forward-measure-F}
   \nu_{\alpha,\beta}^{\mathbb{F}}
   :=
    \nu_{\alpha,\beta}^{\mathcal{T}} \circ (V|_{\partial \mathcal{T}})^{-1} .
\end{equation}
From the injectivity of $V|_{\partial \mathcal{T}}$ we see that $\nu_{\alpha,\beta}^{\mathbb{F}}$ is    non-atomic as well.

\smallskip\noindent
\underline{\textit{Step 3.}} We claim that $\mathcal{R}_{\alpha,\beta} \subset E_r(\alpha,\beta)$ on  $B_\epsilon^{\mathbb{F}}$.   Since $\Lambda_r(\alpha,\beta)= V(E_{r}(\alpha,\beta))$,this implies
\begin{equation}\label{eq-support-fecund-measure}
   \Lambda_r(\alpha,\beta) \text{ is nonempty }   , \text{ and }\ \nu_{\alpha,\beta}^{\mathbb{F}}
   \bigl(\Lambda_r(\alpha,\beta)\bigr)=1 \quad \text{ on } B_\epsilon^{\mathbb{F}} .
\end{equation}
Indeed for $t\in\mathcal{R}_{\alpha,\beta}$, the construction \eqref{eq-def-fecund-tree} gives $|V(t_{L_n})|=M_n$. Thus by
 \eqref{eq-choose-eta-F},  we have
 \begin{equation}
  \label{eq:L-n-speed}
  \liminf_{n \to \infty}\frac{|V(t_{L_n})|}{L_n}=\alpha \quad , \quad \limsup_{n \to \infty}\frac{|V(t_{L_n})|}{L_n}=\beta.
 \end{equation}
Moreover, we show that
\begin{equation}
   \sup_{L_n\leq k\leq L_{n+1} } \frac{ ||V(t_k)|-M_n |}{L_n}
   \leq \frac{ L_{n+1}-L_n }{L_n} = \frac{ \ell_{n+1}(\eta_{n+1})}{L_n} \xrightarrow{n \to \infty} 0 .   \label{eq:L-n-increment}
\end{equation}
Indeed, as $L_n\geq M_n$ and
$\ell_{n+1}(\eta_{n+1})\leq \frac{m_{n+1}}{q_{n+1}^-}+1$, we have
$\frac{\ell_{n+1}(\eta_{n+1})}{L_n}
 \leq \frac{m_{n+1}+q_{n+1}^- }{M_n q_{n+1}^-}$.
 Since $\frac{m_{n+1}}{M_n}=O(\frac{1}{n})$ and $M_n\to\infty$, while
$\inf_n q_n^->0$ when $r<R$ and $nq_n^-\to\infty$ when $r=R$, the assertion \eqref{eq:L-n-increment} follows. Combining this with \eqref{eq:L-n-speed}, we obtain that the lower and upper limits of $(\frac{1}{k}|V(t_k)|)$ are $\alpha$ and $\beta$, respectively.
This shows the claim $t\in E_r(\alpha,\beta)$.

\smallskip\noindent
\underline{\textit{Step 4.}} In Step 3, we have showed $\Lambda_r(\alpha,\beta)$ is nonempty on $B^{\mathbb{F}}_{\epsilon}$. To show
the second inequality in \eqref{eq:reduction-thm-lower-bnd},
 it suffices to consider the case where  $ \theta_r^{\mathbb{F}}(\alpha,\beta)>0$.
Fix any
$ 0<\theta<(1-\epsilon)^2\, \theta_r^{\mathbb{F}}(\alpha,\beta)$. We shall prove that, on $B^{\mathbb{F}}_{\epsilon}$,
\begin{equation}\label{eq-correct-energy-F}
 \mathbf{I} (
 \theta,\nu_{\alpha,\beta}^{\mathbb{F}}
 )
 := \int  \frac{\nu^{\mathbb{F}}_{\alpha,\beta} (\dif \omega ) \nu^{\mathbb{F}}_{\alpha,\beta} (\dif \omega' ) }{\mathrm{d}_{\partial \mathbb{F}}   ( \omega,  \omega'   )^{\theta}}   \,
 \leq \,
 \sum_{n=0}^{\infty}
 e^{\theta M_{n+1}}
 \prod_{j=1}^n\frac1{A_j(\eta_j)} \, .
\end{equation}
In Step 5 we will verify the (deterministic) series on the right-hand side is  finite.
Consequently, using $V(\mathcal{R}_{\alpha,\beta}) \subset \Lambda_r(\alpha,\beta)$, Lemma~\ref{energy-method} gives
\begin{equation}\label{eq-dimension-on-fecund-event}
 \dim_{\mathrm{H}}\Lambda_r(\alpha,\beta)
 \geq
  \theta
 \quad\text{on }B_\epsilon^{\mathbb{F}}.
\end{equation}
The desired result \eqref{eq:reduction-thm-lower-bnd} therefore follows by   letting
$\theta\uparrow(1-\epsilon)^2 \theta_r^{\mathbb{F}}(\alpha,\beta)$.

To prove \eqref{eq-correct-energy-F}, note first that the diagonal point
$\omega'=\omega$ has no contribution to the integral, because
$\nu_{\alpha,\beta}^{\mathbb{F}}$ has no atom. Next, we assert that  $B_\epsilon^{\mathbb{F}}$, for every $\omega = V(t)$ with $t \in \mathcal{R}_{\alpha,\beta}$
\begin{equation}
   \label{eq:w-2-t-1}
   \nu_{\alpha,\beta}^{\mathbb F}
  \bigl(
  [\omega_{M_n}]_{\mathbb F}
  \setminus[\omega_{M_{n+1}}]_{\mathbb F}
  \bigr)=
  \nu_{\alpha,\beta}^{\mathcal T}
  \bigl(
  \mathcal R_{\alpha,\beta}
  \cap
  ([t_{L_n}]_{\mathcal T}
  \setminus[t_{L_{n+1}}]_{\mathcal T})
  \bigr)
  \leq
  \prod_{j=1}^n\frac1{A_j(\eta_j)}.
\end{equation}
 Since these sets form a disjoint partition of
$V(\mathcal{R}_{\alpha,\beta})\setminus\{\omega\}$ and
$e^{\theta|\omega\wedge\omega'|}\leq e^{\theta M_{n+1}}$ on the
$n$-th set, summing the preceding mass estimates over $n $ yields
\eqref{eq-correct-energy-F}.

We know prove \eqref{eq:w-2-t-1}.
Let $t\in\mathcal{R}_{\alpha,\beta}$. By construction, we have  $
    V(t_{L_{k+1}})
    \in
    \mathbb{F}_{m_{k+1}}(V(t_{L_k}) )$ for every $k\geq0$.  Thus the vertices $(V(t_{L_k}))_{k\geq0}$ form a nested
sequence of prefixes of the boundary point $\omega=V(t)$. Since
$|V(t_{L_k})|=M_k$, we get
\begin{equation}\label{eq:block-prefix}
    \omega_{M_k}=V(t_{L_k})  \quad \text{ for all } \quad  k\geq0.
\end{equation}
Let $t'\in\mathcal{R}_{\alpha,\beta}$ and  $\omega'=V(t')$. From
\eqref{eq:block-prefix} and
\eqref{eq:cohort-spatial-injectivity}, it follows for every $n\geq0$,
\begin{equation}\label{eq:skeleton-prefix-correspondence}
t'\in[t_{L_n}]_{\mathcal T}
\quad\Longleftrightarrow\quad
\omega'\in[\omega_{M_n}]_{\mathbb F}.
\end{equation}
Indeed, the forward implication follows from
$t'_{L_n}=t_{L_n}$ and \eqref{eq:block-prefix}. Conversely, if
$\omega'_{M_n}=\omega_{M_n}$, then
$
V(t'_{L_n})=V(t_{L_n}),$
and \eqref{eq:cohort-spatial-injectivity} gives
$t'_{L_n}=t_{L_n}$.
Combining  \eqref{eq-cylinder-mass-fecund}, we derive  \eqref{eq:w-2-t-1}.

\smallskip\noindent
\underline{\textit{Step 5.}}
It remains to prove that the series in
\eqref{eq-correct-energy-F} is finite.
To this end, define
$h(q):={[\ln r-L^*(q)]}/{q}$,
where the value at $q=0$ is understood according to the convention in
Theorem~\ref{thm-Hdim-Lambda-alpha-beta}. With this convention, $h$ is
continuous on $I(r)$. Since
$|m_n/\ell_n(\eta_n)-\eta_n|\leq O({1}/{m_n})$, the continuity of $h$
yields
\[
 \ell_n(\eta_n) \Bigl[
 \ln r-L^* \bigl(  \frac{m_n}{\ell_n(\eta_n)}  \bigr)  \Bigr]
 =m_n\bigl[h(\eta_n)+o(1)\bigr].
\]
The assumption
$\theta_r^{\mathbb{F}}(\alpha,\beta)= \min_{q\in [\alpha,\beta]} h(q)>0$,
together with
$\dist(\eta_n,[\alpha,\beta])\to0$ in
\eqref{eq-choose-eta-F}, implies that
$h(\eta_n)\ge (1-\epsilon)\theta_r^{\mathbb{F}}(\alpha,\beta)$ for all
sufficiently large $n$. Choose $\tilde{\theta}$ such that
$\theta<\tilde{\theta}<(1-\epsilon)^2\theta_r^{\mathbb{F}}(\alpha,\beta)$.
Combining the preceding estimates with the definition of
$A_n(\cdot)$, we obtain
\begin{equation}
 \log A_n(\eta_n)\geq(1-\epsilon) \, h(\eta_n)\, m_n\ge  \tilde{\theta} m_n \quad \text{ for all  large } n.
\end{equation}
It follows that, for some constant $C<\infty$,
$
 \sum_{j=1}^n\log A_j(\eta_j)\geq \tilde{\theta} M_n-C$.
Finally, using $M_{n+1}=M_n+m_{n+1}$ and $m_{n+1}=o(M_n)$, we obtain,
for some $c>0$ and all   large $n$,
\[
 e^{\theta M_{n+1}}\prod_{j=1}^n\frac1{A_j(\eta_j)}
 \leq
 \exp\bigl\{-(\tilde{\theta}-\theta)M_n+\theta m_{n+1}+C\bigr\}
 \leq e^{-cM_n}.
\]
Therefore, the series in \eqref{eq-correct-energy-F} converges. This completes the proof.
\end{proof}

  \subsection{Lower bound in Theorem \ref{thm-dim-E-alpha-beta}}
The proof follows the same fecund skeleton construction and energy method
  as the lower-bound proof of
Theorem~\ref{thm-Hdim-Lambda-alpha-beta}.  Since most of the argument
carries over verbatim, we give details only for the modifications
required in the present setting.  The main difference is that we
count particles rather than distinct sites, and work directly with
the genealogical metric on $\partial\mathcal{T}$.

Let $r \in (1,\infty)$.
Fix $\epsilon\in(0,1/4)$. Let $\mathbb{D}_n $ be  constructed as in preceding subsection   but using the sequence $\hyperref[lem-level-set-size]{\delta}_n$ from
Lemma~\ref{lem-level-set-size},  as the
dyadic uniform grid  of the interval
\[
 \bigl\{
 q\in[0,1-2\hyperref[lem-level-set-size]{\delta}_{n+n_0}]:
 \ln r-L^*(q)\geq2\hyperref[lem-level-set-size]{\delta}_{n+n_0}
 \bigr\}.
\]
Thus $\#\mathbb{D}_n\leq4n$ and
$\bigcup_{n\geq1}\mathbb{D}_n$ is dense in $I(r)$.

Let
$
 k_n=(k_1+n-1)^4  $ with $
 K_n:=\sum_{j=1}^n k_j $.  For each $q\in\mathbb{D}_n$, set
\begin{equation}\label{eq-def-kn-An-T}
 d_n(q):=\langle qk_n\rangle_{k_n,\mu}
 \quad\text{and}\quad
 A_n^{\mathcal{T}}(q)
 :=
 \Bigl\lfloor
 \frac16
 \Bigl(
 r^{k_n}
 \exp\bigl\{
 -k_nL^*\bigl(\frac{d_n(q)}{k_n}\bigr)
 \bigr\}
 \Bigr)^{1-\epsilon/2}
 \Bigr\rfloor.
\end{equation}
By \eqref{eq:solve-periodicity-2}, $d_n(q)\in\mathsf{L}_{k_n}^\mu$ and
$|d_n(q)-qk_n|\leq1$. Hence
$\mathcal{N}_{k_n,d_n(q)}$ is non-degenerate.

Set $
 \underline A_n^{\mathcal{T}}
 :=\min_{q\in\mathbb{D}_n}A_n^{\mathcal{T}}(q)$.
As in \eqref{eq-grid-conditions-fecund} and \eqref{eq-lower-An}, by increasing $k_1$ we may
assume that Lemma~\ref{lem-level-set-size} is applicable uniformly
over $\mathcal{N}_{k_n,d_n(q)}$, and
\begin{equation}\label{eq-lower-An-T}
\underline{A}_n^{\mathcal{T}}\geq2  \quad , \quad  \log\underline A_n^{\mathcal{T}}
 \geq a_0k_n \hyperref[lem-level-set-size]{\delta}_{n+n_0}.
\end{equation}

\smallskip \noindent
\underline{Fecund particles.}
We construct the particle analogue of the fecund skeleton. For
$v\in\mathcal{T}$ and $q\in\mathbb{D}_n$, let
\[
\mathcal{U}_n(v,q)
 :=
 \bigl\{
 u\in\mathcal{T}_{k_n}(v):
 V(u)\in\mathbb{F}_{d_n(q)}(V(v))
 \bigr\}.
\]
We say $v$ is \emph{$q$-fecund for stage $n$} if
$\#\mathcal{U}_n(v,q)\geq6A_n^{\mathcal{T}}(q)$.

\smallskip \noindent
\underline{Continuation successors.}
Assume that $v$ is $q$-fecund for stage $n$. For
$q'\in\mathbb{D}_{n+1}$, we call a particle
$u\in\mathcal{U}_n(v,q)$ a \emph{$q'$-continuation successor} of
$v$ if $u$ is $q'$-fecund for stage $n+1$. We denote all such
successors by
\[
\mathcal{U}_n(v;q\mid q')
 :=
 \bigl\{
 u\in\mathcal{U}_n(v,q):
 u\text{ is $q'$-fecund for stage $n+1$}
 \bigr\} \ \text{ for every } q'\in\mathbb{D}_{n+1}.
\]
If $v$ is not $q$-fecund for stage $n$, we set
$\mathcal{U}_n(v;q\mid q'):=\varnothing$.

\smallskip \noindent
\underline{Selected cohorts and particle fecund skeleton.}
Set
$\mathcal{V}_{\mid q_1}:=\{\emptyset\}$ and define recursively, for
every $(q_1,\dots,q_{n+1})\in\prod_{j=1}^{n+1}\mathbb{D}_j$,
\[
\mathcal{V}_{q_1,\ldots,q_n\mid q_{n+1}}
 :=
 \bigcup_{v\in\mathcal{V}_{q_1,\ldots,q_{n-1}\mid q_n}}
 \mathcal{U}_n(v;q_n\mid q_{n+1}).
\]
We call
$\mathcal{V}_{q_1,\ldots,q_n\mid q_{n+1}}$ a
\emph{selected cohort after stage $n$}. Since its elements are
particles, no representatives are needed. The family of selected
cohorts forms the \emph{particle fecund skeleton}.

By construction, for every $n\geq2$, each particle
$v\in\mathcal{V}_{q_1,\ldots,q_{n-1}\mid q_n}$, if any, is
$q_n$-fecund for stage $n$ and satisfies $|v|=K_{n-1}$.
Moreover, the continuation-successor sets of
distinct particles in the same cohort are disjoint.

 \smallskip \noindent
\underline{Good event.}
Let $B_\epsilon^{\mathcal{T}}$ be the event on which the root $\emptyset$ is $q_1$-fecund for stage $1$ for
every $q_1\in\mathbb{D}_1$;
and for every $n\geq1$, every
$(q_1,\ldots,q_{n+1})\in\prod_{j=1}^{n+1}\mathbb{D}_j$, and every
$v\in\mathcal{V}_{q_1,\ldots,q_{n-1}\mid q_n}$,
\[
\#\mathcal{U}_n(v;q_n\mid q_{n+1})
 \geq A_n^{\mathcal{T}}(q_n) .
\]

\begin{lemma}\label{lem-positive-fecund-event-T}
For every $\epsilon\in(0,1/4)$, one can choose $k_1$
sufficiently large so that
\begin{equation}\label{eq-positive-fecund-event-T}
 \P(B_\epsilon^{\mathcal{T}})\geq\frac12.
\end{equation}
\end{lemma}

\begin{proof}
The proof of Lemma~\ref{lem-fecund-subtree-F} applies with
\eqref{eq:level-set-size} and \eqref{eq:forward-cone-size} in place of
their distinct-site counterparts.  The only change in the final
union bound is that \eqref{eq:trivial-bnd-fecund-tree} is replaced by
the following expectation bound:
\[
 \E\bigl[
\#\mathcal{V}_{q_1,\ldots,q_{n-1}\mid q_n}
 \bigr]
 \leq\E|\mathcal{T}_{K_{n-1}}|
 =r^{K_{n-1}}.
\]
For each $v\in\mathcal{T}$, the same argument as in the proof of
Lemma~\ref{lem-fecund-subtree-F} gives
\[
 \P \bigl(
 \#\mathcal{U}_n(v;q\mid q') <A_n^{\mathcal{T}}(q)
  \,\bigm|\,
 \mathcal{F}^V_{|v|+k_n}
\bigr)
 \ind{v\text{ is $q$-fecund for stage $n$}}
 \leq e^{-A_n^{\mathcal{T}}(q)/4}
\]
 for every $q\in\mathbb{D}_n$ and
$q'\in\mathbb{D}_{n+1}$. Hence a union bound gives
\[
 \P\bigl((B_\epsilon^{\mathcal{T}})^c\bigr)
 \leq
 \frac14+
 \sum_{n\geq1}
 \Bigl(\prod_{j=1}^{n+1}\#\mathbb{D}_j\Bigr)
 r^{K_{n-1}}e^{-\underline A_n^{\mathcal{T}}/4}
 \leq\frac12,
\]
 where the last inequality follows from \eqref{eq-lower-An-T}
after increasing $k_1$.
\end{proof}

\begin{proof}[Proof of the lower bound in Theorem~\ref{thm-dim-E-alpha-beta}]
 By Lemmas~\ref{lem-positive-fecund-event-T}
and~\ref{lem-0-1-law}, as in the preceding proof, it suffices to show
that for every $\epsilon\in(0,1/4)$, on
$B_\epsilon^{\mathcal{T}}$,
\begin{equation}\label{eq:reduction-2}
 E_r(\alpha,\beta)\neq\emptyset
 \quad\text{and}\quad
 \dim_{\mathrm{H}}E_r(\alpha,\beta)
 \geq
 (1-\epsilon)^2 \,\theta^{\mathcal{T}}_{r}(\alpha,\beta)
 \quad\text{for every }[\alpha,\beta]\subset I(r).
\end{equation}

\smallskip\noindent
\underline{\textit{Step 1.}}
 Fix $[\alpha,\beta]\subset I(r)$. As in Step~1 of the preceding
proof, we may choose $\eta_j\in\mathbb{D}_j$ such that, with
$
 D_n
 :=\sum_{j=1}^n d_j(\eta_j)$,
\begin{equation}\label{eq-choose-eta-fecund-T}
 \liminf_{n\to\infty}\frac{D_n}{K_n}=\alpha,
 \qquad
 \limsup_{n\to\infty}\frac{D_n}{K_n}=\beta,
 \qquad
 \dist(\eta_n,[\alpha,\beta])\longrightarrow0.
\end{equation}

\smallskip\noindent
\underline{\textit{Step 2.}}
On $B_\epsilon^{\mathcal{T}}$, let
$\mathcal{R}_{\alpha,\beta}^{\mathcal{T}}$ consist of all
$t\in\partial\mathcal{T}$ such that
\[
 t_{K_n}\in
\mathcal{V}_{\eta_1,\ldots,\eta_n\mid\eta_{n+1}}
 \qquad\text{for every }n\geq1.
\]
Starting with unit mass at the root, at stage $n$ we divide the mass
of each $v\in\mathcal{V}_{\eta_1,\ldots,\eta_{n-1}\mid\eta_n}$ uniformly
among the particles in
$\mathcal{U}_n(v;\eta_n\mid\eta_{n+1})$.
Since each such set contains at least
$A_n^{\mathcal{T}}(\eta_n)\geq2$ particles, this defines a non-atomic
probability measure $\nu_{\alpha,\beta}^{\mathcal{T}}$ on
$\mathcal{R}_{\alpha,\beta}^{\mathcal{T}}$ satisfying
\begin{equation}\label{eq-cylinder-mass-fecund-T}
 \nu_{\alpha,\beta}^{\mathcal{T}}
 \bigl(\mathcal{R}_{\alpha,\beta}^{\mathcal{T}}
 \cap[t_{K_n}]_{\mathcal{T}}\bigr)
 =
 \prod_{j=1}^n
 \frac{1}{
\#\mathcal{U}_j
 (t_{K_{j-1}};\eta_j\mid\eta_{j+1})}
 \leq
 \prod_{j=1}^n\frac1{A_j^{\mathcal{T}}(\eta_j)}
\end{equation}
for every $t\in\mathcal{R}_{\alpha,\beta}^{\mathcal{T}}$.

\smallskip\noindent
\underline{\textit{Step 3.}}
For each ray
$t\in\mathcal{R}_{\alpha,\beta}^{\mathcal{T}}$, the forward-cone
construction gives
$
 |V(t_{K_n})|=D_n$ and $
 \sup_{K_n\leq k\leq K_{n+1}}
 \frac{ ||V(t_k)|-D_n |}{K_n}
 \leq\frac{k_{n+1}}{K_n}\longrightarrow0$.
Together with \eqref{eq-choose-eta-fecund-T}, this implies
\begin{equation}\label{eq-support-fecund-T}
 \mathcal{R}_{\alpha,\beta}^{\mathcal{T}}
 \subset E_r(\alpha,\beta).
\end{equation}
In particular, $E_r(\alpha,\beta)$ is nonempty on
$B_\epsilon^{\mathcal{T}}$.

\smallskip\noindent
\underline{\textit{Step 4.}}
It remains to show the second inequality in \eqref{eq:reduction-2}.  We may assume
$\theta_r^{\mathcal{T}}(\alpha,\beta)>0$ and show that on $B_\epsilon^{\mathcal{T}}$ for any
$
  0<\theta<
 (1-\epsilon)^2\theta_r^{\mathcal{T}}(\alpha,\beta)$,
$
  \mathbf{I} (
 \theta,\nu_{\alpha,\beta}^{\mathcal{T}}
 ) <\infty  $
which implies $   \dim_{\mathrm{H}}E_r(\alpha,\beta)\geq\theta $ by Lemma~\ref{energy-method}.

For $t'\neq t$ in $\mathcal{R}_{\alpha,\beta}^{\mathcal{T}}$,
partition according to the unique $n$ such that
$t'\in[t_{K_n}]_{\mathcal{T}}\setminus[t_{K_{n+1}}]_{\mathcal{T}}$.
Since this set is contained in $[t_{K_n}]_{\mathcal{T}}$,
\eqref{eq-cylinder-mass-fecund-T} gives
\begin{equation}\label{eq-correct-energy-T}
 \mathbf{I}\bigl(
 \theta,\nu_{\alpha,\beta}^{\mathcal{T}}
 \bigr)
 \leq
 \sum_{n=0}^{\infty}
 e^{\theta K_{n+1}}
 \prod_{j=1}^n\frac1{A_j^{\mathcal{T}}(\eta_j)}
\end{equation}

\smallskip\noindent
\underline{\textit{Step 5.}}
The convergence argument is the same as in Step~5 of the preceding
proof, with $g(q):=\ln r-L^*(q)$ in place of $h(q)$.  More precisely,
\eqref{eq-choose-eta-fecund-T} and the definition of
$A_n^{\mathcal{T}}$ imply that, for every
$
 \theta<\widetilde\theta<
 (1-\epsilon)^2\theta_r^{\mathcal{T}}(\alpha,\beta)$
and all   large $n$,
$
 \log A_n^{\mathcal{T}}(\eta_n)
 \geq\widetilde\theta k_n.$
Consequently, for some $C,c>0$ and all sufficiently large $n$,
\[
 e^{\theta K_{n+1}}
 \prod_{j=1}^n\frac1{A_j^{\mathcal{T}}(\eta_j)}
 \leq
 \exp\bigl\{
 -(\widetilde\theta-\theta)K_n+\theta k_{n+1}+C
 \bigr\}
 \leq e^{-cK_n}.
\]
Thus the series in \eqref{eq-correct-energy-T} is finite. This completes the proof.
\end{proof}

\section{More questions}
In addition to Questions \ref{ques1}, \ref{ques2},  \ref{ques3}, and Remark \ref{rmk-rw-LDP},
we conclude by listing several intriguing questions related to our proofs, which could serve as potential topics for further research.

\begin{enumerate}[(i)]
  \item  Characterize the symmetric probability measures $\mu$ on
$\mathcal A$ for which the pressure comparison condition
\eqref{eq-pressure-comparison} holds. Determine whether this condition
is necessary and sufficient for \eqref{eq-hdim-Lambda-alph} to hold
throughout the transient regime.

 \item Prove Lemma \ref{lem-Nnx-is-small} under the weaker assumption that the offspring distribution  of the BRW has  only finite variance.  If successful, we can eliminate
 the assumption that the offspring distribution has exponential moments from Theorems \ref{thm-Hdim-Lambda-alpha} and \ref{thm-Hdim-Lambda-alpha-beta}. Indeed this assumption  is used exclusively in the proof of Lemma \ref{lem-Nnx-is-small}, and we believe it may be unnecessary.

  \item Motivated by Section \ref{sec-pathLD}, we are interested in
the following  question:
\begin{equation}
   \text{What is the typical value of $Z_{k}$ conditioned on $Z_{n}=x$?}
\end{equation}
There seems to be no reason to expect that  a result analogous  to \eqref{eq-smaple-path-RW}  would hold for the random walk $Z_n$ on $\mathbb{F}$.
Based on the proof of Proposition \ref{prop-LDP-sample-paths},  we believe the following assertion holds, which
suggests a partial answer to the previous question.
  Let $\Phi(\xi):= (\psi_{a}'(s(\xi)))_{a \in \mathcal{A}}$  appeared in Lemma \ref{l:concave}.  For any $\delta>0$ there exists a constant $C_{\delta}> 0$ such that   for large  $n $
  \begin{equation}
  \max_{x \in \mathbb{F}, |x| \leq n }\mathbf{P} \Bigl( \exists\, k \leq n,  \big|  \big\langle \Xi(Z_{k})- \frac{k}{n} \Xi(x),  \Phi\big( \frac{\Xi(x)}{n} \big) \big\rangle  \big|  > \delta n     \,\big|\,  Z_n = x \Bigr) \leq e^{  - C_{\delta} n}.
  \end{equation}

  \item In Lemma \ref{lem-concentration-Nnx} we  prove that for $x \in \mathbb{F}$ with relative small word length, there is a non-negligible probability that  $ \mathcal{N}_{n,x}$ is of the same order as its expectation.  We guess that,  however, in the recurrent regime with large offspring mean $r$,
  \begin{equation}
   \limsup_{n \to \infty } \max_{x \in \mathbb{F},|x| \leq n } \Big|\frac{1}{n}\ln \mathcal{N}_{n,x} - \frac{1}{n}\ln  \mathbf{E}[\mathcal{N}_{n,x}] \Big| >0.
  \end{equation}
 To address this question.   a good understanding for  question (iii) above is essential.

\end{enumerate}

\section*{Acknowledgement}
We extend our sincere gratitude to an anonymous referee for carefully reading the manuscript, identifying several mistakes, and providing helpful and constructive comments.
We would also like to thank all the organizers of the 2nd International Conference of Young Probability Scholars in China, where this work was initiated. HM thanks his PhD advisor, Yan-Xia Ren, for her encouragement. LW is supported by the National Natural Science Foundation of China under Grant No. 12171252.

%%%%%%%%%%%%%%%%%%%%%%%%%%%%%%%%%%%
%%%%%%%%%%%%%%%%%%%%%%%%%%%%%%%%%%%
%%                               %%
%%           Reference           %%
%%                               %%
%%%%%%%%%%%%%%%%%%%%%%%%%%%%%%%%%%%
%%%%%%%%%%%%%%%%%%%%%%%%%%%%%%%%%%%
 \bibliographystyle{alpha}
  \bibliography{msbrw}

  %%%%%%%%%%%%%%%%%%%%%%%%%%%%%%%%%%%%%%%%%%%%%%%%
  %%%%%%%%%%%%%%%%%%%%%%%%%%%%%%%%%%%%%%%%%%%%%%%%
  %%%%                                        %%%%
  %%%%             APPENDIX                   %%%%
  %%%%                                        %%%%
  %%%%%%%%%%%%%%%%%%%%%%%%%%%%%%%%%%%%%%%%%%%%%%%%
  %%%%%%%%%%%%%%%%%%%%%%%%%%%%%%%%%%%%%%%%%%%%%%%%

  %%%%%%%%%%%%%%%%%%%%%%%%%%%%%%%%%%%%%%%%%%%%%%%%
%%%%%%%%%%%%%%%%%%%%%%%%%%%%%%%%%%%%%%%%%%%%%%%%
%%%%                                        %%%%
%%%%             APPENDIX                   %%%%
%%%%                                        %%%%
%%%%%%%%%%%%%%%%%%%%%%%%%%%%%%%%%%%%%%%%%%%%%%%%
%%%%%%%%%%%%%%%%%%%%%%%%%%%%%%%%%%%%%%%%%%%%%%%%

\appendix

\section{Proof of Lemmas in Section \ref{sec:ldp}}

\subsection{Concavity of the rate function $\Psi^{*}$: Proof of Lemma \ref{l:concave}}
\label{app-pf-concave}

\begin{proof}[Proof of Lemma \ref{l:concave}]
For each fixed $s\leq \ln R$, the map
$\xi\mapsto \langle \xi,\psi(s)\rangle-s$ is affine. Hence, by
\eqref{eq-def-Psi*}, $\Psi^{*}$ is concave on $\Omega$ as the pointwise
infimum of affine functions.

Now let $0<\|\xi\|_1<1$. The infimum in \eqref{eq-def-Psi*} is attained
at the unique $s=s(\xi)$ satisfying \eqref{eq-s-xi}. Hence 
$\Psi^{*}(\xi)=\langle\xi,\psi(s(\xi))\rangle-s(\xi)$. For every $b\in\mathcal{A}$, taking partial derivative  with respect to $\xi_b$ and using
\eqref{eq-s-xi}, we obtain, 
\[
\frac{\partial\Psi^{*}}{\partial\xi_b}(\xi)
=\psi_b(s(\xi))
+\bigl(\langle\xi,\psi'(s(\xi))\rangle-1\bigr)
  \frac{\partial s}{\partial\xi_b}(\xi)
=\psi_b(s(\xi)).
\]
In other words, $\nabla\Psi^{*}(\xi)=\psi(s(\xi))$ for $0<\|\xi\|_1<1$.
Moreover, \eqref{eq-s-xi} implies that $s(\xi)\uparrow\ln R$ as
$\xi\to0$ and that $s(\lambda\xi)\in[s(\xi),\ln R)$ for $0<\lambda\leq1$.
Applying the mean value theorem to $\lambda\mapsto\Psi^{*}(\lambda\xi)$, we obtain,
for some $\theta_\xi\in(0,1)$,
$
\Psi^{*}(\xi)-\Psi^{*}(0)
=\langle\psi(s(\theta_\xi\xi)),\xi\rangle
=\langle\psi(\ln R),\xi\rangle+o(\|\xi\|_1)$,
which shows that $\nabla\Psi^{*}(0)=\psi(\ln R)$. This proves assertion (i).

For any $\xi\in\Omega\setminus\{\mathbf 0\}$ and
$0\leq\lambda_1<\lambda_2\leq1$, the monotonicity of $\psi_a$ and
\eqref{eq-properties-psi} give,
\begin{align*}
&\bigl[\langle\lambda_1\xi,\psi(s)\rangle-s\bigr]
-\bigl[\langle\lambda_2\xi,\psi(s)\rangle-s\bigr]=(\lambda_1-\lambda_2)\langle\xi,\psi(s)\rangle
 >0.
\end{align*}
 for every $s\leq\ln R$, using $\psi_a(s) \le \psi_a(\ln R)<0$ for all $a \in \mathcal{A}$.
Taking the infimum over $s$ in \eqref{eq-def-Psi*} yields
$\Psi^{*}(\lambda_1\xi)>\Psi^{*}(\lambda_2\xi)$, proving assertion~(2).

Finally, let $\xi$ and $h$ be as in assertion~(3). For every $t$ with $\xi+th\in\Omega$, orthogonality
gives
\[
\langle\xi+th,\psi'(s(\xi))\rangle
=\langle\xi,\psi'(s(\xi))\rangle=1.
\]
Since $\psi_a'(s(\xi))>1$ for every $a\in \mathcal{A}$, this identity implies
$0<\|\xi+th\|_1<1$. By the uniqueness in \eqref{eq-s-xi}, we have 
$s(\xi+th)=s(\xi)$, which implies 
\[
\Psi^{*}(\xi+th)
=\Psi^{*}(\xi)+t\langle h,\psi(s(\xi))\rangle.
\]
Thus $\Psi^{*}$ is linear on the stated segment.
\end{proof}
 
\subsection{Leading eigenvalues: An ingredient in the proof of Lemma \ref{lem-eigenvalue}}

\begin{lemma}\label{lem-PF-eigen-equa}
Let $\nu=(\nu_a)_{a\in\mathcal{A}}\in(0,\infty)^{\mathcal{A}}$ and define
$M=(M_{a,b})_{a,b\in\mathcal{A}}$ by
$M_{a,b}=\nu_a\ind{b\neq a^{-1}}$. The Perron--Frobenius eigenvalue of
$M$ is the largest positive solution $\rho$ of
\begin{equation}
  \sum_{a\in\mathcal{A}}
  \nu_a\frac{\rho-\nu_{a^{-1}}}{\rho^2-\nu_a\nu_{a^{-1}}}=1.
\end{equation}
If $\nu_a=\nu_{a^{-1}}$ for every $a\in\mathcal{A}$, this equation reduces to $
  \sum_{a\in\mathcal{A}}\frac{\nu_a}{\rho+\nu_a}=1 $. 
\end{lemma}

\begin{proof}
 By the Perron--Frobenius theorem, there exists an eigenvector $u=( u_a )_{a\in\mathcal{A}}$ of $M$ with eigenvalue $\rho_{\mathrm{PF}} $  such that  $\sum_{a} u_a=1$ and  $u_a>0$ for all $a\in\mathcal{A}$.  We get 
\[
\rho_{\mathrm{PF}}u_a=\nu_a(1-u_{a^{-1}}),\qquad
\rho_{\mathrm{PF}}u_{a^{-1}}=\nu_{a^{-1}}(1-u_a).
\]
Since $|\mathcal{A}|/2=d\geq2$, we have $u_a+u_{a^{-1}}<1$. Multiplying  identities above 
  therefore shows  
$\rho_{\mathrm{PF}}^2>\nu_a\nu_{a^{-1}}$. Eliminating $u_{a^{-1}}$ gives
\[
u_a=\nu_a
\frac{\rho_{\mathrm{PF}}-\nu_{a^{-1}}}
     {\rho_{\mathrm{PF}}^2-\nu_a\nu_{a^{-1}}}.
\]
Summing over $a \in \mathcal{A}$ proves that $\rho_{\mathrm{PF}}$ satisfies the stated
equation.

Conversely, let $\rho>0$ satisfy that equation and set
\[
v_a:=\nu_a\frac{\rho-\nu_{a^{-1}}}
                    {\rho^2-\nu_a\nu_{a^{-1}}}.
\]
Then $\sum_a v_a=1$, and
\[
(Mv)_a=\nu_a\sum_{b\neq a^{-1}}v_b
=\nu_a(1-v_{a^{-1}})=\rho v_a.
\]
Thus every positive solution is an eigenvalue of $M$ and hence is at
most $\rho_{\mathrm{PF}}$. This proves the first assertion. The symmetric
case follows by factoring
$\rho^2-\nu_a^2=(\rho-\nu_a)(\rho+\nu_a)$.
\end{proof}

 %%%%%%%%%%%%%%%%%%%%%%%%
%% refer to arxiv
%%%%%%%%%%%%%%%%%%%%%%%

\subsection{Large deviations for empirical measures:
Proof of Lemma \ref{lem-ldp-empirical-measure}}
\label{app-pf-empirical-measure}
The key observation is that $\xi(X_n)$ is exactly the empirical
measure of a stationary Markov chain. Let $W=(W_k)_{k\geq1}$ be the Markov
chain on $\mathcal{A}$  transition probabilities
\begin{equation}\label{Eq-transition-prob}
  p(a,b):=\P(W_{k+1}=b\mid W_k=a)
  =\frac{1}{2d-1}\ind{b\neq a^{-1}},
  \qquad a,b\in\mathcal{A},\ k\geq1.
\end{equation}
with uniform initial distribution
$\P(W_1=a)=1/(2d)$, $a\in\mathcal{A}$.    Then   $W=(W_n)_{n\geq 1}$ is a stationary   and $
  X_{n} :=W_1 \cdots W_n $
  has the uniform distribution on $\mathbb{F}_{n}$. 
  
Define the  empirical measure and the pair empirical measure   of $W$ respectively by 
\begin{equation}
	\mathcal{L}_n := \frac{1}{n} \sum_{i=1}^n \delta_{  W_i }    \quad  \text{ and } \quad  \mathcal{L}_{n}^{(2)} := \frac{1}{n} \sum_{i=1}^n \delta_{(W_i ,   W_{i+1}) } \ \text{ for } n \geq 1.
\end{equation} 
It suffices to study the large deviation probability  for  $\mathcal{L}_n$, since 
 \[   \xi(X_n)  = (\mathcal{L}_n(\{a\}))_{a \in \mathcal{A} } . \] 
 
Let $\mathcal{P}^{(2)}$ be the collection of all probability distributions on $\mathcal{A} \times \mathcal{A}$ which are absolutely continuous with respect to $p(\cdot,\cdot)$ defined in \eqref{Eq-transition-prob}. For each $\pi \in \mathcal{P}^{(2)}$,  denote by $\pi_1(\cdot)=\sum_{b\in\mathcal{A}}\pi(\cdot,  b)$ and  $\pi_2(\cdot)=\sum_{b\in\mathcal{A}}\pi(b,  \cdot)$ the two marginal distributions of $\pi$.  Define 
\begin{equation}\label{eq-def-I2-rate-function}
  I^{(2)}(\pi) = \sum_{ a,  b \in \mathcal{A} } \pi( a,   b ) \ln \frac{ \pi(a,  b)  }{\pi_1(a) p(a,  b)}      
\end{equation}
  with standard notational conventions $0\ln 0 = 0 \ln \frac{0}{0}=0$ for each $\pi \in \mathcal{P}^{(2)}$.  
  
  Let  $\mathcal{P}_{\mathrm{b}}^{(2)}$ be the set of all $\pi \in \mathcal{P}^{(2)}$ satisfying   the balance condition $\pi_{1}=\pi_{2}$, i.e.,
  \begin{equation}\label{eq-balance-cond}
    \sum_{b\in\mathcal{A}}\pi(a,  b)= \sum_{b\in\mathcal{A}}\pi(b,a) \ \text{ for all }  a \in \mathcal{A}.
  \end{equation}
  
%Let $\mathcal{G}$ be the oriented graph induced by the transition probability $p(\cdot,\cdot)$ in \eqref{Eq-transition-prob}. That is, let $\mathcal{A}$ be the set of vertices and an oriented edge $(a,b)$ presents if and only if $b \neq a^{-1}$. A path in $\mathcal{G}$ is a sequence . 

For each  probability $\pi$ on $\mathcal{A} \times \mathcal{A}$, let $T_n^{(2)}(\pi)$ denote the collection of all possible trajectories of the Markov chain that have pair empirical measure $\pi$, i.e., 
\begin{equation}
	T_n^{(2)}(\pi )=     \Bigl\{  (w_{i})_{i=1}^{n+1} \in \mathcal{A}^{n+1}:   \frac{1}{n}\sum_{i=1}^n \delta_{ ( w_i, w_{i+1} )}=\pi \ ,  \   p(w_{i},  w_{i+1}) >0 , 1 \leq i \leq n    \Bigr\}
\end{equation}
Let  $ \mathcal{P}_n^{(2)} = \{ \pi \in \mathcal{P}^{(2)} :    T_n^{(2)}(\pi) \neq \emptyset \}$ be the set of all possible values that $ \mathcal{L}_{n}^{(2)}$ can take. 
Note that for any $\pi \in \mathcal{P}_n^{(2)}$, $\pi(a,b)>0$ only if $p(a,b)>0$.

For the pair empirical measure $\mathcal{L}_n^{(2)}$, Csiszar, Cover, and Choi \cite{CCC87} provided the following useful Large deviation estimate.

\begin{lemma}[{\cite[Lemma 3]{CCC87}}]\label{eq-LDP-pair-empirical-measure} 
  For every $\pi \in \mathcal{P}_n^{(2)}$, we have 
	\begin{equation}
   \frac{1}{2d} (n+1)^{-(4d^2+2d)} \exp\{-n  I^{(2)}(\pi)\} \leq \mathbf{P} (\mathcal{L}_{n}^{(2)}=\pi) \leq \exp\{-n  I^{(2)}(\pi)\}. 
	\end{equation} 
\end{lemma}

Let $\mathcal{P}$ be the set of all probability distributions on $\mathcal{A}$. For each $n \geq 1$, let  
\[ \mathcal{P}_{n} = \{ \nu \in\mathcal{P}:   n\nu(a)  \in  \mathbb{N}_0  \text{ and } [\nu(a) +\nu(a^{-1})]\ind{ \nu(a) \nu(a^{-1}) > 0 } < 1 , \forall a\in\mathcal{A}   \}. \] 
 Note that $\mathcal{P}$ can be identified  with $ \Omega^{1}$ and  $\mathcal{P}_{n}$ with $ \Omega^{1}_n$ by  viewing each distribution $\nu$ on $\mathcal{A}$ as a vector. The second condition encodes reducedness: a reduced word containing both $a$ and $a^{-1}$ must also contain another letter.

\begin{corollary}\label{cor-Ln-In}
For $n\geq1$ and $\nu\in\mathcal{P}_n$, define
  \begin{equation}
    I_n(\nu):=\min\bigl\{I^{(2)}(\pi):
    \pi\in\mathcal{P}_n^{(2)},\ \pi_1=\nu\bigr\}.
  \end{equation}
Then, uniformly in $\nu\in\mathcal{P}_n$,
$\mathbf P(\mathcal{L}_n=\nu)=\exp\{-nI_n(\nu)+O(\ln n)\}$ as
$n\to\infty$.
\end{corollary}

\begin{proof}
Since $\mathcal{L}_n$ is the first marginal of $\mathcal{L}_n^{(2)}$,
Lemma~\ref{eq-LDP-pair-empirical-measure} and
$\#\mathcal{P}_n^{(2)}\leq(n+1)^{4d^2}$ give
\[
\frac{1}{2d}(n+1)^{-(4d^2+2d)}e^{-nI_n(\nu)}
\leq \mathbf P(\mathcal{L}_n=\nu)
\leq (n+1)^{4d^2}e^{-nI_n(\nu)}.
\]
Taking logarithms proves the claim uniformly in $\nu\in\mathcal{P}_n$.
\end{proof}

\begin{lemma}\label{cor-In-I}
For $\nu\in\mathcal{P}$, define
\begin{equation}
	I(\nu):=\min\bigl\{I^{(2)}(\pi):
  \pi\in\mathcal{P}_{\mathrm b}^{(2)},\ \pi_1=\nu\bigr\}.
\end{equation}  
Then $I(\nu)=\varrho^*(\nu)$.
\end{lemma}

\begin{proof}
The minimum is attained because the feasible set is compact. By
\cite[Exercises~3.1.17 and 3.1.19(a)]{DZ09}, we have 
\[
I(\nu)=\sup_{\lambda\in\mathbb R^{2d}}
\bigl\{\langle\nu,\lambda\rangle-\ln\rho(P_\lambda)\bigr\},
\qquad
P_\lambda:=\bigl(p(a,b)e^{\lambda_b}\bigr)_{a,b\in\mathcal{A}},
\]
and $\rho(P_\lambda)$ is the Perron--Frobenius eigenvalue of
$P_\lambda$. Applying Lemma~\ref{lem-PF-eigen-equa} to the transpose of $P_\lambda$
and comparing with \eqref{eq-def-of-varrho} give
\[
\ln\rho(P_\lambda)=\varrho(\lambda)-\ln(2d-1).
\]
Therefore, we get 
$
I(\nu)=\sup_{\lambda\in\mathbb R^{2d}}
\bigl\{\langle\nu,\lambda\rangle-\varrho(\lambda)\bigr\}
+\ln(2d-1)=\varrho^*(\nu)$ as required.
\end{proof}

\begin{proof}[Proof of Lemma \ref{lem-ldp-empirical-measure}]
It suffices to show that 
\begin{equation}\label{eq-erreo-I-In}
  \Delta_{n} := \max_{\nu \in \mathcal{P}_{n}} | I(\nu)- I_{n}(\nu) | \lesssim_{d} \ln n /n . 
\end{equation}
Indeed, Corollary~\ref{cor-Ln-In} gives
$\ln\mathbf P(\mathcal{L}_n=\nu)=-nI_n(\nu)+O(\ln n)$ uniformly in $\nu$, while
Lemma~\ref{cor-In-I} identifies $I(\nu)=\varrho^*(\nu)$. Thus
\eqref{eq-erreo-I-In} allows us to replace $I_n(\nu)$ by $\varrho^*(\nu)$ at
an additional cost of only $O(\ln n)$ in the exponent, which proves the claim.

To estimate $\Delta_n$, we use the following combinatorial approximation
lemma, whose proof is deferred to the end of this subsection. For
$\nu\in\mathcal{P}_n$, write
$\mathcal{A}^\nu:=\{a\in\mathcal{A}:\nu(a)>0\}$ and call $\nu$ \emph{good} if,
for every $a\in\mathcal{A}$,
\[
a^{-1}\notin\mathcal{A}^\nu
\qquad\text{or}\qquad
n\bigl(\nu(a)+\nu(a^{-1})\bigr)\leq n-2.
\]
Set $\mathcal{P}_{n,\mathrm g}:=\{\nu\in\mathcal{P}_n:\nu\text{ is good}\}$.

\begin{lemma}\label{lem-approx-balance-flow}
There exists $K=K(d)>0$ such that, for all sufficiently large $n$, the
following holds. If $\nu\in\mathcal{P}_{n,\mathrm g}$ and
$\pi\in\mathcal{P}_{\mathrm b}^{(2)}\cup\mathcal{P}_n^{(2)}$ satisfies
$\pi_1=\nu$, then there exists
$\widetilde\pi\in\mathcal{P}_{\mathrm b}^{(2)}\cap\mathcal{P}_n^{(2)}$ such that
\[
\widetilde\pi_1=\nu
\qquad\text{and}\qquad
 \|\pi-\widetilde\pi\|_\infty\leq K/n .
\]
\end{lemma}

For $\pi\in\mathcal{P}^{(2)}$, plugging  the expression of $p(a,b)$ in\eqref{eq-def-I2-rate-function} shows $I^{(2)}(\pi)$ can be written as $\sum_{a,b}\pi(a,b)\ln\pi(a,b)
-\sum_a\pi_1(a)\ln\pi_1(a)+\ln(2d-1).$
Since $x\mapsto x\ln x$ has modulus of continuity
$O(\delta\ln(1/\delta))$ on $[0,1]$, we obtain  for every fixed $C $ and $\pi, \pi' \in \mathcal{P}^{(2)}$,
\begin{equation}\label{eq-continuity-I2}
\|\pi-\pi'\|_\infty\leq C/n
\quad\Longrightarrow\quad
|I^{(2)}(\pi)-I^{(2)}(\pi')| \lesssim  {\ln n}/{n}.
\end{equation}

Fix $\nu\in\mathcal{P}_{n,\mathrm g}$ and set
\[
J_n(\nu):=\min\bigl\{I^{(2)}(\pi):
\pi\in\mathcal{P}_n^{(2)}\cap\mathcal{P}_{\mathrm b}^{(2)},\
\pi_1=\nu\bigr\}.
\]
Applying Lemma~\ref{lem-approx-balance-flow} to minimizers defining
$I(\nu)$ and $I_n(\nu)$, and then using \eqref{eq-continuity-I2}, gives
\begin{equation}\label{eq-I-In-IIn}
0\leq J_n(\nu)-I(\nu)\lesssim   {\ln n}/{n},
\qquad
0\leq J_n(\nu)-I_n(\nu)\lesssim  {\ln n}/{n}.
\end{equation}
particularly, this gives $|I(\nu)-I_n(\nu)|\lesssim (\ln n)/n$ for every good
$\nu$.

It remains to consider $\nu\in\mathcal{P}_n\setminus\mathcal{P}_{n,\mathrm g}$.
Set $\pi_\nu(c,c'):=\nu(c)\ind{c'=c}$.
Since $\nu$ is not good, there exists $a\in\mathcal A$ such that
$a^{-1}\in\mathcal A^\nu$ and
\[
 n\bigl[\nu(a)+\nu(a^{-1})\bigr] \in \{n-1,n\}
\]
If this value is $n$, the definition of $\mathcal P_n$ forces
$\nu=\delta_{a^{-1}}$. In this case, $\pi_\nu$ is the unique balanced pair
measure with first marginal $\nu$, while every $\pi'\in\mathcal P_n^{(2)}$
with $\pi'_1=\nu$ satisfies $\|\pi'-\pi_\nu\|_\infty\leq1/n$. Hence
\eqref{eq-continuity-I2} gives
$|I(\nu)-I_n(\nu)|\lesssim(\ln n)/n$.

We may therefore assume that $    \nu(a)+\nu(a^{-1})  =1-1/n$. Since $n\nu$ is integer-valued, there
exists $b\in\mathcal A$ such that
\[
\nu(a)+\nu(a^{-1})=1-1/n,
\qquad \nu(b)=1/n,
\qquad \mathcal{A}^\nu\subseteq\{a,a^{-1},b\}.
\]

Assume that 
$\pi\in\mathcal{P}_{\mathrm b}^{(2)}$ has first marginal $\nu$.
Since transitions between $a$ and $a^{-1}$ are forbidden, any
off-diagonal entry $\pi(c,c')>0$ must satisfy $c=b$ or $c'=b$.
As $\pi_1(b)=\pi_2(b)=\nu(b)=1/n$, it follows that
$\|\pi-\pi_\nu\|_\infty\leq 2/n$.
Now let $\pi'\in\mathcal{P}_n^{(2)}$ have first marginal $\nu$, and choose a
trajectory $(w_i)_{i=1}^{n+1}\in T_n^{(2)}(\pi')$. After possibly interchanging $a$ and $a^{-1}$, there is unique $1<l<n$ such that  $w_l=b$,   $w_i=a$ for $1 \le i < l$, $w_j=a^{-1}$ for $l<j \le n$. This gives  $ \sum_{c'\neq c}\pi'(c,c') \le 3/n$. 
Hence $\|\pi'-\pi_\nu\|_\infty\leq 6/n$. 
We get 
$\|\pi-\pi'\|_\infty\leq 8/n$, and applying
\eqref{eq-continuity-I2} to minimizers of $I(\nu)$ and $I_n(\nu)$ proves
\eqref{eq-erreo-I-In}. 
\end{proof}

\begin{proof}[Proof of Lemma \ref{lem-approx-balance-flow}]
 We split the proof into two cases.

 \smallskip
\noindent\underline{Case 1: $\pi\in\mathcal P_{\mathrm b}^{(2)}$.}
Here $\pi$ is already balanced; the missing property is that it need not
be realizable as the pair empirical measure of a  length-$n$
trajectory.  Indeed, $n\pi$ is a  balanced flow whose in-degree and
out-degree at each $a\in\mathcal A^\nu$ are both  $n\nu(a) \in \mathbb{N}$.
A natural idea is   to replace $n\pi$ by a nearby nonnegative
integer flow $q$, still supported on allowed transitions, such that
$\|q-n\pi\|_\infty=O (1)$, its in-degrees and out-degrees remain
$n\nu$, and its support is strongly connected.  Let $G_q$ be the directed
multigraph having $q(a,b)$ edges from $a$ to $b$.  Since $q$ is balanced,
every vertex of $G_q$ has equal in-degree and out-degree; since the support
of $q$ is strongly connected, $G_q$ has an Eulerian circuit.  Moreover,
$G_q$ has $n$ edges, so this circuit realizes
$q/n\in\mathcal P_{\mathrm b}^{(2)}\cap\mathcal P_n^{(2)}$.

To this end, first we construct
a bounded balanced integer flow $c$ whose support connects all of
$\mathcal A^\nu$.  We then extract from $n\pi$ a large balanced integer
subflow $\varpi$, leaving only a uniformly bounded residual.  Finally, we
add $c$ and use self-loops to restore the prescribed degrees $n\nu$.
For any flow $f$, we write $f_1(a):=\sum_b f(a,b)$.

\smallskip
\noindent\underline{Step 1.} 
We first construct a  closed walk that visits every vertex of 
$\mathcal A^\nu$. 

If $\mathcal A^\nu=\{a\}$, take the self-loop
$a\longrightarrow a$.  Otherwise, consider the simple graph with vertex
set $\mathcal A^\nu$ and edge set
$
 \{\{a,b\}:a,b\in\mathcal A^\nu,\ b\neq a, b \neq a^{-1}  \}$. 
If $\#\mathcal A^\nu\geq4$, every vertex has degree at least
\[
 \#\mathcal A^\nu-2\geq\frac{\#\mathcal A^\nu}{2},
\]
so Dirac's theorem gives a Hamiltonian cycle. If
$\#\mathcal A^\nu=2$, its two points cannot be an inverse pair by the
definition of $\mathcal P_n$, so the resulting two-cycle is allowed.  A
three-point support admits an allowed three-cycle, unless it is of the form
$\{a,a^{-1},b\}$. In that exceptional case, use
$
 a\to b\to a^{-1}\to b
 \to a$.
Here the goodness of $\nu$ gives $n\nu(b)\geq2$. 
Thus in every case there is
an allowed closed walk $(v_0,\ldots,v_m=v_0)$ that leaves each vertex exactly once, except that it leaves $b$ twice in the exceptional case.

Let $c(a,b)$ be the number of times this walk traverses the directed edge $(a,b)$. Then
$c$ is a balanced integer flow with strongly connected support $\mathcal{A}^{\nu}$, and
\begin{equation}\label{eq:backbone-flow}
   c_1(a)\leq2 \ , \  c_1(a)\leq n\nu(a) \quad  \forall \,  a \in \mathcal{A} \quad 
   \text{and}\quad \|c\|_\infty\leq1.
\end{equation}

\smallskip
\noindent\underline{Step 2.} 
Among all balanced integer flows $f$ satisfying
\[
   0\leq f(a,b)\leq
   [n\pi(a,b)-c_1(a)]_+,
   \qquad a,b\in\mathcal A,
\]
choose one with maximal total mass and denote it by $\varpi$. The admissible family is finite and nonempty because it contains the zero flow. The residual flow $R:=n\pi-\varpi$ is balanced. We claim   
\begin{equation}
  \label{eq:residual-flow}
     \|n\pi-\varpi\|_\infty\leq3(2d)^{2d}.
\end{equation} 
Indeed, suppose that
$R(a_0,a_1)>3(2d)^{2d}$. Since $R$ is balanced, we may recursively
choose $a_{i+1}$ so that
$R(a_i,a_{i+1})\geq R(a_{i-1},a_i)/(2d)$.
Some vertex among $a_0,\ldots,a_{2d}$ must repeat, yielding a directed
cycle $C$ on which $R(a,b)\geq3$. Adding one unit to $\varpi$ along $C$
preserves balance; moreover, for every $(a,b)\in C$,
\[
  \varpi(a,b)+1
  \leq n\pi(a,b)-2
  \leq n\pi(a,b)-c_1(a),
\]
so the upper bounds remain valid. This contradicts the maximality of
$\varpi$.

Moreover, we have 
\begin{equation}\label{eq:row-reserve}
   \varpi_1(a)\leq n\nu(a)-c_1(a) \quad \forall a \in \mathcal{A}. 
\end{equation}
Indeed, this is immediate if
$n\pi(a,b)\leq c_1(a)$ for every $b$. Otherwise, choose $b_0$ with
$n\pi(a,b_0)>c_1(a)$. Then
$\varpi(a,b_0)\leq n\pi(a,b_0)-c_1(a)$, whereas
$\varpi(a,b)\leq n\pi(a,b)$ for $b\neq b_0$, and hence
\[
   \varpi_1(a)
   \leq n\pi(a,b_0)-c_1(a)
      +\sum_{b\neq b_0}n\pi(a,b)
   =n\nu(a)-c_1(a).
\]

\smallskip
\noindent\emph{Step 3.}
We now define for $a,b \in \mathcal{A}$
\begin{equation}\label{eq:def-tilde-pi}
 n\widetilde\pi(a,b):=
 \begin{cases}
   \varpi(a,b)+c(a,b),&a\neq b,\\
   \varpi(a,a)+c(a,a)
   +n\nu(a)-\varpi_1(a)-c_1(a),&a=b.
 \end{cases} 
\end{equation}
By \eqref{eq:row-reserve}, the diagonal corrections in
\eqref{eq:def-tilde-pi} are nonnegative integers. Moreover, 
$n\widetilde\pi$ is a balanced integer flow with first marginal $n\nu$.
All its transitions are allowed, and its support contains the strongly
connected support of $c$ and is therefore strongly connected.  Moreover,
for $a\neq b$, equations
\eqref{eq:backbone-flow} and \eqref{eq:residual-flow} give
$ |n\widetilde\pi(a,b)-n\pi(a,b) |
 \leq 1+3(2d)^{2d}$. 
Since $n\widetilde\pi$ and $n\pi$ have the same first marginal, we  have
$
 n\widetilde\pi(a,a)-n\pi(a,a)
 =-\sum_{b\neq a}
    [n\widetilde\pi(a,b)-n\pi(a,b) ]$.
Thus $\|n\widetilde\pi-n\pi\|_\infty\leq K $.  The directed multigraph
associated with $n\widetilde\pi$ is balanced, strongly connected, and has
$n$ edges, so it has an Eulerian circuit.  This circuit realizes
$\widetilde\pi$, as required.

\smallskip
\noindent\underline{Case 2: $\pi\in\mathcal P_n^{(2)}$.}
Here $\pi$ is already realized by an allowed length-$n$ trajectory; only
balance may be missing.  For any
$(w_i)_{i=1}^{n+1}\in T_n^{(2)}(\pi)$,
\[
 \pi_1-\pi_2=\frac{1}{n}  (\delta_{w_1}-\delta_{w_{n+1}} ).
\]
Thus it is enough to rearrange the first $n$ letters, without changing
their multiplicities, so that every consecutive transition, including the
last-to-first transition, is allowed.  This preserves the first marginal
$\nu$, while changing only a bounded number of pair counts.

Fix such a trajectory and put $x:=w_1$. If $w_n\neq x^{-1}$, we simply
replace the terminal letter
$w_{n+1}$ by $x$, thereby obtaining an allowed closed walk with the same
first marginal and changing at most one pair count.

Suppose now that $w_n=x^{-1}$. Write the maximal terminal block of
$w_1\cdots w_n$ as $(x^{-1})^r$, and let $b$ be the letter immediately
preceding it. Reducedness gives $b\notin\{x,x^{-1}\}$. If
\[
   w_1\cdots w_n=A\,b\,(x^{-1})^r
\]
and the last letter of $A$ is not $x$, replace this word by
$A\,(x^{-1})^r b$. Otherwise, write $A=B x^s$, where $x^s$ is the maximal
terminal block of $x$ in $A$. Since $\nu$ is good, at least two of
$w_1,\ldots,w_n$ lie outside $\{x,x^{-1}\}$; hence $B$ is nonempty and its
last letter lies outside $\{x,x^{-1}\}$. In this case, replace the word
$  B x^s b(x^{-1})^r$
by $B\,(x^{-1})^r b x^s$.

In both cases the new word has the same letter counts and becomes an allowed
closed walk when its last letter is joined to its first. Its edge-count
matrix is obtained from the original one by removing and adding at most four
unit edges. Hence its pair empirical measure $\widetilde\pi$ satisfies
\[
   \widetilde\pi\in
   \mathcal P_{\mathrm b}^{(2)}\cap\mathcal P_n^{(2)},\qquad
   \widetilde\pi_1=\nu,\qquad
   \|\pi-\widetilde\pi\|_\infty\leq 4/n
\]
This completes the proof.
\end{proof}
  
\end{document}